\documentclass[preprint, 11pt]{elsarticle}
\usepackage[margin=2.5cm]{geometry}

\newtheorem{theorem}{Theorem}

\newtheorem{remark}[theorem]{Remark}

\usepackage{hyperref}

\journal{a journal}

\usepackage{amsmath}
\usepackage{amsfonts}
\usepackage{amssymb}
\usepackage{bm}
\usepackage{mathrsfs}
\usepackage{physics}
\usepackage{upgreek}
\usepackage{mathtools}

\usepackage[labelformat=simple]{subcaption}

\usepackage{graphicx}

\usepackage{booktabs}
\usepackage{multirow}
\usepackage{caption}
\usepackage[table]{xcolor}
\definecolor{lightgray}{gray}{0.85}
\usepackage{tablefootnote}

\usepackage{float}
\makeatletter
\newcommand\fs@boxedtop
{\fs@boxed
	\def\@fs@mid{\vspace\abovecaptionskip\relax}%
	\let\@fs@iftopcapt\iftrue
}
\makeatother
\floatstyle{boxedtop}
\floatname{framedbox}{Box}
\newfloat{framedbox}{hbt}{lob}[section]
\usepackage{setspace}

\usepackage[textsize=tiny]{todonotes}

\begin{document}
	
	\begin{frontmatter}
		
		\title{Automatically Differentiable Model Updating (ADiMU): conventional, hybrid, and neural network material model discovery including history-dependency}
		
		\author[mymainaddress]{Bernardo P. Ferreira\corref{mycor2}}
		\ead{bernardo_ferreira@brown.edu}
		\author[mymainaddress]{Miguel A. Bessa\corref{mycor2}}
		\ead{miguel_bessa@brown.edu}
		
		\address[mymainaddress]{School of Engineering, Brown University, 184 Hope St, Providence, RI 02912, USA}
		
		\cortext[mycor2]{Corresponding author}
		
		\begin{abstract}
			We introduce the first Automatically Differentiable Model Updating (ADiMU) framework that finds any history-dependent material model from full-field displacement and global force data (global, indirect discovery) or from strain-stress data (local, direct discovery). We show that ADiMU can update conventional (physics-based), neural network (data-driven), and hybrid material models. Moreover, this framework requires no fine-tuning of hyperparameters or additional quantities beyond those inherent to the user-selected material model architecture and optimizer. The robustness and versatility of ADiMU is extensively exemplified by updating different models spanning tens to millions of parameters, in both local and global discovery settings. Relying on fully differentiable code, the algorithmic implementation leverages vectorizing maps that enable history-dependent automatic differentiation via efficient batched execution of shared computation graphs. This contribution also aims to facilitate the integration, evaluation and application of future material model architectures by openly supporting the research community. Therefore, ADiMU is released as an open-source computational tool, integrated into a carefully designed and documented software named HookeAI.
		\end{abstract}
		
		\begin{keyword}
			Material model, Model updating, Automatic differentiation, History-dependency, Recurrent neural network, Hybrid material model, ADiMU, Open-source
		\end{keyword}
		
	\end{frontmatter}
	
	
	
	\section{Introduction}
	
	At the core of computational mechanics lies the need to represent material constitutive behavior accurately \citep{lemaitre:1990a, simo:2006a, desouzaneto:2008a}. They have been a key enabler of important breakthroughs in structural design and analysis that span a wide variety of fields and applications \cite{bonet:2008a, gurtin:2010a, belytschko:2013a}, such as crashworthiness analysis in automotive design, the design of lightweight, damage-tolerant aerospace structures, the development of biocompatible medical implants, or the stability analysis of geological structures.
	
	Conventional (mechanistic or phenomenological) constitutive material models are formulated in rigorous physics-based theoretical frameworks \citep{hill:1950a, truesdell:1965a, malvern:1969a, ogden:1984a}, ensuring thermodynamical consistency and fundamental principles like frame invariance and material objectivity. Each model includes a small set of parameters, often linked to measurable physical quantities, and comprises a set of explicit constitutive equations. These are usually grounded in phenomenological or micromechanical considerations, and tailored to capture specific types of material behavior (e.g., elasto-plasticity, visco-plasticity, crystal plasticity, damage) \citep{gurson:1977a, kothari:1998a, holzapfel:2000a, anand:2003a}. Despite their accuracy in simple, idealized strain-stress loading paths, the underlying assumptions limit their performance when dealing with more challeging scenarios found in real applications that involve multi-axial, non-monotonic loading paths, complex multiphasic material architectures, multi-scale damage effects, and evolving history-dependent dependencies.
	
	In contrast, neural network models are highly expressive and can theoretically approximate any continuous function to arbitrary accuracy, as established by the universal approximation theorem \citep{chen:1995a}. Deep learning architectures that integrate appropriate activation functions and recurrent structures can effectively capture nonlinear, history-dependent behavior, as first shown in \cite{mozaffar:2019}. However, these models lack physics knowledge, relying entirely on data-driven learning \cite{bessa:2017a}. With architectures spanning tens to billions of non-physically interpretable parameters, their accuracy heavily depends on the availability of diverse, representative data sets \citep{bengio:2009a, rolnick:2017a}. Neural networks were first introduced as material models by Ghaboussi and coworkers \citep{ghaboussi:1991a} and were followed by early developments \citep{theocaris:1993a, yagawa:1996a, waszczyszyn:2001a, lefik:2002a}, but only gained traction twenty years later, largely driven by the emergence of accessible machine learning software with built-in automatic differentiation capabilities \citep{mozaffar:2019, ghavamian:2019a, wu:2020a}. A thorough review on data-driven material models was recently published by De Lorenzis and coworkers \citep{fuhg:2025a}, providing a valuable overview of the extensive research on the topic. Soon envisioned by Lefik, Unger and coworkers \citep{lefik:2009a, unger:2009a} and early explored by Le and coworkers \citep{le:2015a}, one of the most popular applications of these models involves replacing lower-scale models in hierarchical multi-scale simulations \citep{bessa:2017a, wang:2018a, liu:2019a}, which would otherwise be prohibitively computationally expensive. At the same time, considerable effort is being put forth to design thermodynamic consistent architectures and enforce fundamental constraints similar to conventional models \citep{he:2022a, eghbalian:2023a, kalina:2025a, jadoon:2025a}. In this work, we define a hybrid material model as any model that combines physics-based knowledge and/or constraints with neural network architectures.
	
	Regardless of their nature, the parameters of material models must be determined from measurable experimental quantities, such as displacements and forces. However, given that material models are often defined in terms of quantities that cannot be directly measured (e.g., strains, stresses, internal variables), discovering their parameters calls for the solution of an indirect inverse problem. Introduced by Kavanagh and Clough in the early 1970s \citep{kavanagh:1971a}, one of the most popular, longstanding approaches to tackle this problem is called Finite Element Model Updating (FEMU). The key idea is to leverage the finite element method as a forward propagation model and iteratively update the material parameters by minimizing the discrepancy with known experimental data. Although a thorough, comprehensive literature review was recently published by Chen and coworkers \citep{chen:2024a}, and is thus not attempted here, two key aspects are highlighted to contextualize our contribution.
	
	In the first place, the historical evolution of FEMU has mirrored advances in both experimental data acquisition and optimization techniques. While early contributions used mostly global force-displacement curves from mechanical tests, the emergence of Digital Image Correlation (DIC) and Digital Volume Correlation (DVC) provided access to full-field displacement  data \citep{bay:1999a, sutton:2009a, pan:2009a, buljac:2018a, wang:2024a}. Similarly, foundational work mostly explored gradient-free optimization methods or finite-differences-based sensitivities, after which adjoint methods gained traction for the efficient computation of gradients with respect to numerous parameters \citep{choi:2006a, mroz:2007a, friswell:2010a, haftka:2013a}. Adjoint methods continue to play a central role, namely when combined with automatic differentiation \citep{wengert:1964a, griewank:1989a}, which only recently gained widespread adoption in scientific computing \citep{baydin:2018a}. These advances enabled a progression from early linear models, to complex nonlinear conventional models accounting for history-dependency, and ultimately to expressive, data-driven neural network models \citep{thakolkaran:2022a, wu:2025a, akerson:2025a}. In the second place, it is remarkable that Ghaboussi and coworkers \citep{ghaboussi:1998a} pioneered the integration of a neural network model in the FEMU framework in the late 90s, soon followed by a few contributions \citep{shin:2000a, lefik:2003a}. Despite their novelty, these early works were constrained by the absence of automatic differentiation, limited computational resources, and the unavailability of rich, time-resolved full-field data. Consequently, they relied on shallow networks and lacked robust training procedures, limiting their ability to capture complex, nonlinear, and history-dependent behavior.
	
	Lastly, it is worth noting that alternative approaches to FEMU have also made substantial contributions, such as the virtual fields method \citep{grediac:1989a, pierron:2012a, lourenco:2024a, sun:2025a, kumar:2025a} and physics informed neural networks \citep{dissanayake:1994a, raissi:2019a, jeong:2024a, anton:2024a}. Furthermore, with the growing potential of neural network material models, it is crucial to understand the data diversity needed for effective learning, particularly when addressing history-dependent material behavior \citep{zhang:2024a}.

	\subsection{Our contribution}
	Building on the foundational vision of Ghaboussi and coworkers \citep{ghaboussi:1991a, ghaboussi:1998a}, and to the best of our knowledge, we introduce the first Automatically Differentiable Model Updating (ADiMU) framework that finds any history-dependent material model from full-field displacement and global force data (global, indirect discovery) or from strain-stress data (local, direct discovery). Alongside the streamlined handling of conventional (physics-based), neural network (data-driven) and hybrid material models, ADiMU requires no fine-tuning of hyperparameters or additional quantities beyond those ineherent to the user-selected material model architecture and optimizer. This is extensively demonstrated throughout this paper, covering numerous examples of history-dependent material models in both local and global discovery settings.
	
	This contribution also aims to facilitate the integration, evaluation and application of future material model architectures by openly supporting the research community. Therefore, ADiMU is released as an open-source computational tool, integrated into a carefully designed and documented software named HookeAI (see \ref{sec:hookeai}). This software, fully designed and implemented by the first author, is used to generate most of the (synthetic) data, perform all computational analyses, and to carry out the post-processing of all results shown in this paper.
	
	The paper is outlined as follows. We describe the key steps of ADiMU's material model discovery workflow in Section~\ref{sec:adimu}. Then, we demonstrate ADiMU's performance in Sections~\ref{sec:local_model_discovery} and ~\ref{sec:global_model_discovery}, discussing several local and global material model discovery examples, respectively. Lastly, we summarize the main conclusions and discuss future challenges in Section~\ref{sec:conclusions}. Several appendices are also included with complementary methods, results and details.
	
	\section{Automatically Differentiable Model Updating (ADiMU) \label{sec:adimu}}
\begin{figure}[hbt]
		\centering
		\includegraphics[width=0.95\textwidth]{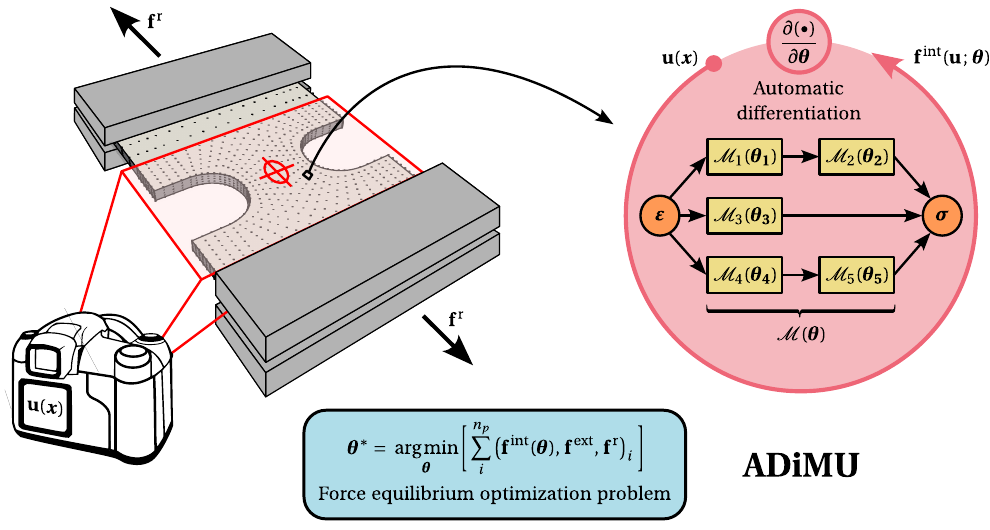}
		\caption{Automatically Differentiable Model Updating (ADiMU) framework, enabling the automatic, indirect discovery of any history-dependent parametric material model, $\mathcal{M}(\bm{\theta})$, parameterized by $\bm{\theta}$, from full-field displacement and global force data. The gradient-based solution of a force equilibrium optimization problem is unlocked with automatic differentiation, as ADiMU's workflow is fully automatically differentiable.}
		\label{fig:adimu_global_discovery}
	\end{figure}

	The Automatically Differentiable Model Updating (ADiMU) framework introduced in this paper enables the automated, indirect discovery of any history-dependent material model from measured full-field displacement and global force data, as illustrated in Figure~\ref{fig:adimu_global_discovery}. This section outlines the fundamental concepts of ADiMU, complemented by several appendices, which collectively provide the necessary background to comprehend the extensive results presented in the remainder of this paper.
	
	As a starting point, consider a specimen for which the underlying material is completely unknown. Suppose we use a Universal Testing Machine (UTM) to perform a uniaxial tensile test, at any given temperature, and collect experimental data, namely: (i) the specimen's displacement field, $\bm{u}(\bm{x})$, throughout the deformation history (e.g., with Digital Image Correlation (DIC) or Digital Volume Correlation (DVC)); and (ii) the corresponding reaction forces, $\bm{f}^{\, \mathrm{r}}$, history (e.g., Universal Testing Machine load cell). With this displacement-force data in hand, the goal is to automatically discover a parametric model, $\mathcal{M}(\bm{\theta})$, parameterized by $\bm{\theta}$, that accurately describes the (local) material behavior,
	\begin{equation}
		\bm{\sigma}(\bm{x}) = \mathcal{M}(\bm{\varepsilon}(\bm{x})\, ;\, \bm{\theta}) \, ,
	\end{equation}
	where $\bm{x}$ denotes a spatial coordinate, $\bm{\sigma}$ is the Cauchy stress tensor and $\bm{\varepsilon}$ is the infinitesimal strain tensor.\footnote{While this paper demonstrates ADiMU in the infinitesimal strains setting, the framework is the same regardless of the material model formulation under infinitesimal or finite strains. In the latter setting, the displacement field is used to compute the deformation gradient, $\bm{F}(\bm{x})$, which is subsequently employed to determine any required strain tensor and convert any stress tensor into the Cauchy stress tensor.}
	
\begin{figure}[hbt]
		\centering
		\includegraphics[width=0.95\textwidth]{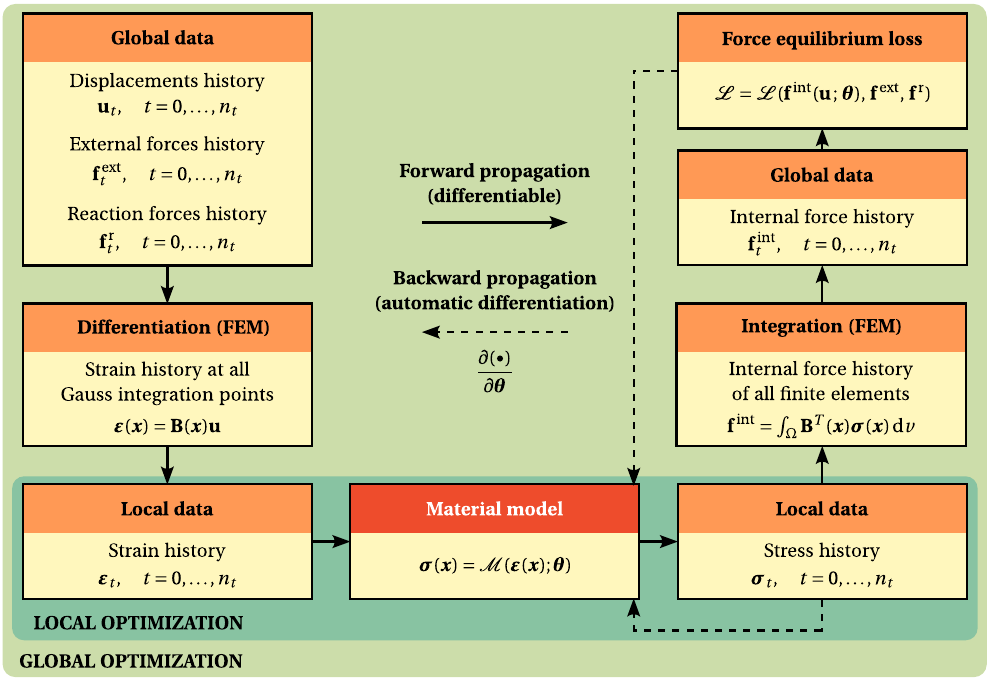}
		\caption{Material model discovery workflow with Automatically Differentiable Model Updating (ADiMU). In the local discovery setting, a general material model, $\mathcal{M}$, parameterized by $\bm{\theta}$, is discovered directly from strain-stress data, $(\bm{\varepsilon}, \, \bm{\sigma})$. In contrast, in the global discovery setting, the material model is indirectly discovered from displacement-force data, $(\bm{u}, \, \bm{f})$. Note that ADiMU's local forward propagation is shared between both local and global discovery settings.}
		\label{fig:adimu_optimization}
	\end{figure}

\begin{figure}[hbt]
		\centering
		\includegraphics[width=0.85\textwidth]{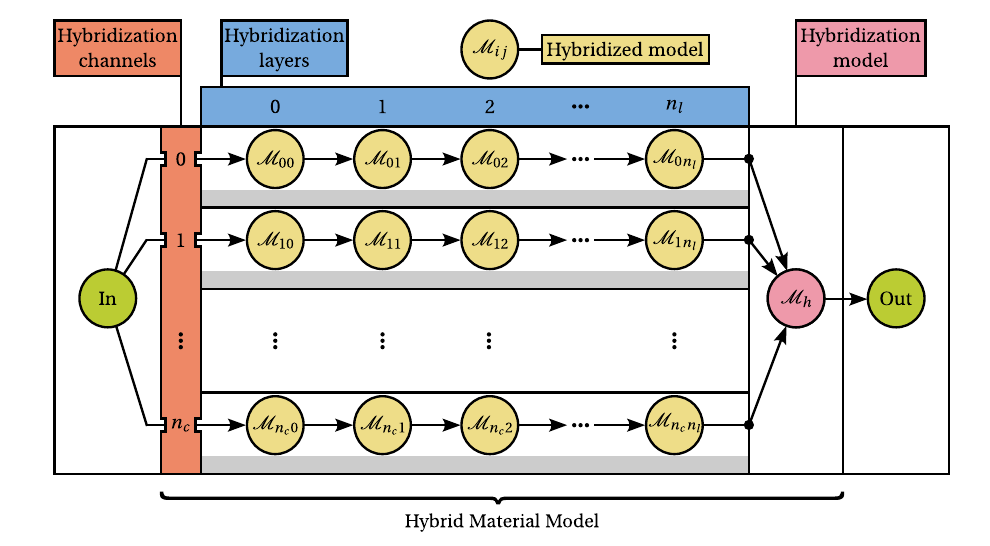}
		\caption{Hybrid material model architecture. The material input data flows into each hybridization channel, each containing one or more hybridized models (one per hybridization layer). The outputs from the hybridization channels are then processed by the hybridization model, which yields the final material output data.}
		\label{fig:hybrid_architecture}
	\end{figure}

	ADiMU's material model discovery workflow is illustrated in Figure~\ref{fig:adimu_optimization}. The overall framework involves the following key steps:
	\begin{enumerate}
		\item \textbf{Displacement-force data enconding.} The first step discretizes the specimen in a finite element mesh and encodes the experimental displacement-force history data as nodal quantities. While a point-to-point mapping is required to encode the DIC/DVC displacement field data, a standard Universal Tensile Machine (UTM) only provides the global uniaxial reaction force history. Given the homogeneous strain state often found at the specimen's shoulders (gripped by the UTM), the global uniaxial reaction force can be uniformly distributed over the corresponding nodes as a reasonable approximation;
		\item \textbf{Choice of parametric model architecture.} The second step consists in selecting a suitable parametric model, $\mathcal{M}(\bm{\theta})$, to describe the material constitutive behavior. Note that the parameters $\bm{\theta}$ can be completely or partially unknown. Akin to any finite element simulation software, $\mathcal{M}(\bm{\theta})$ can be any conventional material model (e.g., the well-known von Mises \citep{mises:1913} or Drucker-Prager \citep{drucker:1952a} elasto-plastic models), where $\bm{\theta}$ could include parameters such as the elastic constants or hardening properties. However, ADiMU can handle practically any material model, namely neural network models, where $\bm{\theta}$ often includes the so-called weights and biases parameters. In fact, our implementation also includes a general hybrid material model architecture, illustrated in Figure~\ref{fig:hybrid_architecture}. This architecture enables a flexible combination of any conventional and neural network models, each integrated as a hybridized model with a corresponding set of parameters, in a graph-like structure of hybridization channels and layers\footnote{The number of hybridized models in each hybridization channel is independent. Moreover, despite not illustrated in Figure~\ref{fig:hybrid_architecture}, the hybrid model architecture allows residual and interchannel connections compatible with the forward data flow.}. The outputs of the different hybridization channels are then collected by a hybridization model (e.g., a weighted mixture rule), itself a general parametric model, that yields the resulting material stress prediction. In this case, the parameters $\bm{\theta}$ of the hybrid material model, $\mathcal{M}(\bm{\theta})$, include all parameters of the underlying hybridized models and hybridization model. While we demonstrate a particular candidate-corrector hybrid architecture in this paper, different possibilities are discussed in the future remarks;
		\item \textbf{Selection of optimizer.} The third step entails selecting an appropriate optimizer and its associated hyperparameters. In the proposed implementation, the underlying workflow computation graph is entirely automatically differentiable, i.e., gradients with respect to the parameters, $\bm{\theta}$, are computed via automatic differentiation. Consequently, gradient-based optimization is available for the model discovery process, which is essential for handling material models with parameter counts ranging from a few to millions. Moreover, some parameters may require setting bounds (e.g., avoid non-physical parameters) or the enforcement of model specific constraints (e.g., preserve the yield surface convexity), which is particularly relevant in conventional models;
		
		\begin{remark}
			We have not changed the optimizer and its hyperparameters for all the examples considered in this paper. We did not need to fine-tune and/or compare different types of optimizers\footnote{Gradient-free optimizers can be used when the number of parameters is small, namely when considering conventional models. In this case, gradient computation via automatic differentiation should be disabled to prevent unnecessary computational costs.}, as convergence was not an issue. All model discovery processes are performed with the well-known gradient-based optimizer Adam \citep{kingma:2014a}, assuming default decay rates for first and second moment estimates, 0.9 and 0.999, respectively, and a learning rate exponential decay scheduler. The suitable range for the learning rate is dependent on factors such as data and parameters normalization, model architecture, and the resulting loss landscape. Nevertheless, a learning rate range of [$\mathit{10^{-1},\, 10^{1}}$] is employed for all conventional models, whereas the range [$\mathit{10^{-5},\, 10^{-3}}$] is adopted for all neural network and hybrid models.
		\end{remark}
		
		\item \textbf{Material model discovery.} ADiMU's model discovery workflow involves a force equilibrium optimization problem formulated using an implicit version of the Finite Element Method (FEM). Knowing that the experimental displacement field corresponds \textit{de facto} to an equilibrium state, the objective function driving the optimization problem is naturally based on the static (or quasi-static) force equilibrium of a solid structure. Akin to most FEMU contributions \citep{chen:2024a}, the commonly defined force equilibrium loss is defined as \citep{giton2006hyperelastic,avril2008overview}
		\begin{equation}
			\mathcal{L}(\bm{\theta}) = \sum_{t}^{n_{t}} \sum_{j}^{n_{p}} \sum_{k}^{n_{d}} \Big[ \mathbf{f}^{\, \mathrm{int}}_{j, \, k} (\bm{\theta}) - \mathbf{f}^{\, \mathrm{ext}}_{j, \, k} - \mathbf{f}^{\, \mathrm{r}}_{j, \, k} \Big]_{t}^{2} \, ,
		\end{equation}
		where $\mathbf{f}^{\, \mathrm{int}}$, $\mathbf{f}^{\, \mathrm{ext}}$ and $\mathbf{f}^{\, \mathrm{r}}$ denote the internal, external and reaction forces, respectively, while $n_{t}$, $n_{p}$ and $n_{d}$ denote, respectively, the number of (pseudo-)time steps\footnote{In this article we refer to load increments as pseudo-time increments, but note that in the entire article there is no time variable (all analyses are quasi-static). Yet, there is no limitation to consider time-dependent problems.}, the number of nodes and the number of degrees of freedom per node. In this paper, we propose the previous loss function to become dimensionless as
		\begin{equation}
			\mathcal{L}(\bm{\theta}) = \alpha \sum_{t}^{n_{t}} \Bigg\{ \sum_{j}^{n_{p}} \sum_{k}^{n_{d}} \Big[ \mathbf{f}^{\, \mathrm{int}}_{j, \, k} (\bm{\theta}) - \mathbf{f}^{\, \mathrm{ext}}_{j, \, k} - \mathbf{f}^{\, \mathrm{r}}_{j, \, k} \Big]_{t}^{2} \Bigg\} \Delta t \, ,
			\label{eq:force_equilibrium_loss}
		\end{equation}
		where $\Delta t$ denotes a (pseudo-)time increment and the scaling factor, $\alpha$, is defined as
		\begin{equation}
			\alpha^{-1} = E_{c}^{\, 2} \; l_{c}^{\, 4} \; t_{c} \, ,
		\end{equation}
		where $E_{c}^{2}$ is a characteristic Young modulus, $l_{c}$ is a characteristic length, and $t_{c}$ is a characteristic (pseudo-)time length. The force equilibrium optimization problem is thus postulated as
		\begin{equation}
			\bm{\theta}^{*} = \underset{\bm{\theta}}{\mathrm{argmin}} \, \mathcal{L}(\bm{\theta}) \, ,
		\end{equation}
		where $\bm{\theta}^{*}$ is the set of material model parameters that best satisfies the force equilibrium of the specimen throughout the whole deformation history. As previously mentioned, the gradient of the force equilibrium loss with respect to the parameters, $\nabla_{\bm{\theta}} \, \mathcal{L}(\bm{\theta})$, is computed via automatic differentiation;
		
		\begin{remark}
			In general, the force equilibrium optimization problem is non-convex. On the one side, the non-convexity of the loss landscape is induced by the displacement-force nature of the optimization problem, i.e., the indirect dependency of the internal forces on the material model parameters through the force equilibrium equations. On the other side, as the material model becomes more expressive and nonlinear, namely when involving neural network architectures, the optimization loss landscape becomes increasingly non-convex. Consequently, multiple local minima, saddle points and other features of non-convex loss landscapes are expected and should be considered when selecting a suitable optimizer.
		\end{remark}
		
		\item \textbf{Material model prediction performance.} After solving the optimization problem, the accuracy and reliability of the discovered material model, $\mathcal{M}(\bm{\theta}^{*})$, should be then tested on diverse multi-axial, non-monotonic loading paths.
		
	\end{enumerate}
	
\begin{figure}[hbt]
		\centering
		\includegraphics[width=0.9\textwidth]{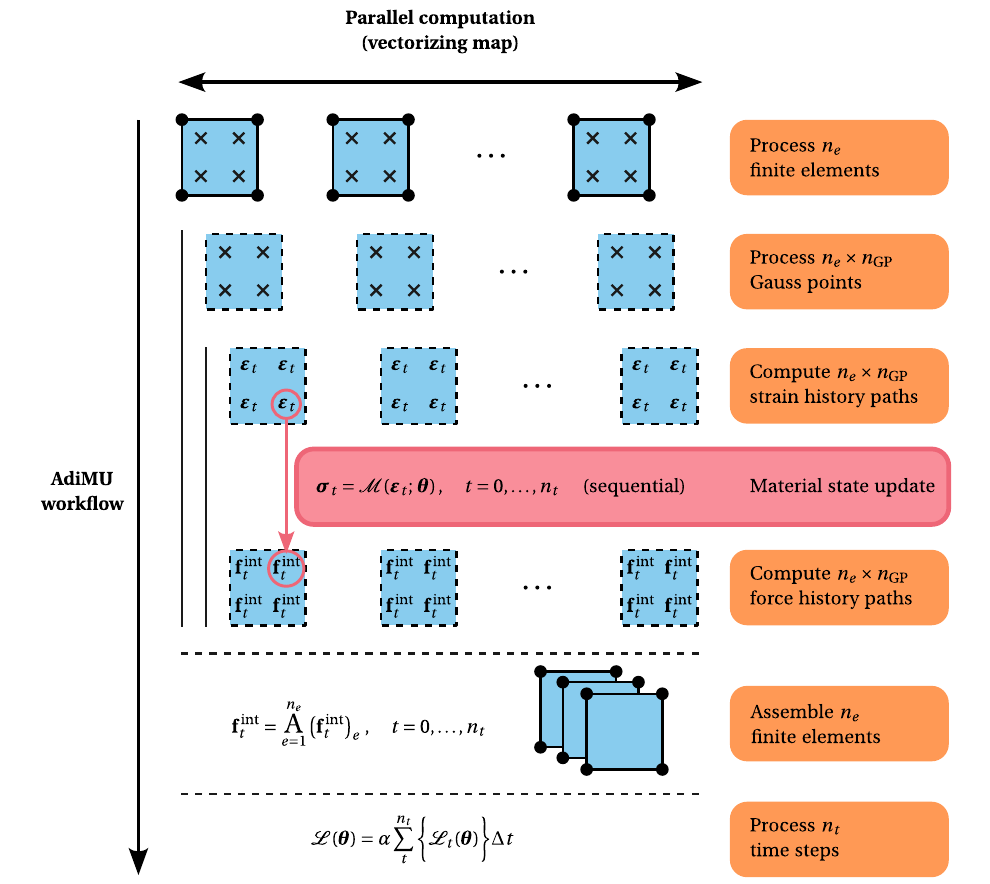}
		\caption{Algorithmic design of the Automatically Differentiable Model Updating (ADiMU) framework. Vectorizing maps significantly reduce computational costs by leveraging batched execution with shared computation graphs.}
		\label{fig:adimu_vectorizing_maps}
	\end{figure}

	Given that the full displacement field history is known from the experimental measurements, ADiMU's forward propagation does not involve the solution of the equilibrium problem like a standard implicit FEM simulation. Nevertheless, the material model state update at each integration point may still require the solution of a nonlinear system of equations, namely when a nonlinear conventional model is involved (e.g., elasto-plastic return-mapping). This motivates the proposal of a rather non-conventional algorithmic design that significantly leverages vectorizing maps in the computation of the force equilibrium loss (see Figure~\ref{fig:adimu_vectorizing_maps}). The use of vectorizing maps enables substantial efficiency gains in the evaluation of the differentiable computation graph over all elements in the finite element mesh, particularly relevant under history-dependent material behavior. By eliminating explicit loops and leveraging batched execution with shared computation graphs, vectorization significantly reduces both runtime and memory overhead. This is especially critical in the context of automatic differentiation and GPU acceleration, where vectorized operations allow more efficient memory management and parallel computation.
	
	\clearpage
	
	\begin{remark}
		While history-independent material models are commonly addressed in the literature, history-dependent models -- such as those involving internal variables or incorporating recurrence -- introduce significant computational and memory demands due to their dependency on the entire deformation history. Under such conditions, a properly vectorized implementation is not just beneficial but essential, particularly when full automatic differentiation is employed to compute gradients. Without vectorization, the accumulation of history-dependent operations and gradient tracking makes the problem computationally intractable.
	\end{remark}
	
	Last but not least, a note on ADiMU's versatily is in order. While so far we have focused on global material model discovery, where the material model is indirectly discovered from displacement-force data, it is worth highlighting that local material model discovery, where the material model is discovered directly from strain-stress data, is also available without any modifications other than the selection of a different loss function. In fact, it transpires from ADiMU's global workflow and design architecture (see Figures~\ref{fig:adimu_optimization} and \ref{fig:adimu_vectorizing_maps}) that the local strain-to-stress computations are common to both local and global discovery settings forward propagation. In particular, given a local data set of strain-stress paths (akin to the batched Gauss integration points strain-stress paths), the optimization problem loss function can be set, e.g., as the Mean Squared Error (MSE),
	\begin{equation}
		\mathcal{L}(\bm{\theta}) = \mathrm{MSE} \, (\bm{\sigma}(\bm{\theta}), \, \bar{\bm{\sigma}}) \, ,
	\end{equation}
	where $\bm{\sigma}(\bm{\theta})$ denotes the stress paths predicted by the material model, $\mathcal{M}(\bm{\theta})$, and $\bar{\bm{\sigma}}$ denotes the stress paths (numerical) ground-truth\footnote{Similar to the force equilibrium loss function (see Equation~\ref{eq:force_equilibrium_loss}), the whole deformation history is considered in the loss computation.}. In this case, $\bm{\theta}^{*}$ is the set of material model parameters that best explains the (numerical) ground-truth material stress response. Although stress data is not experimentally measurable, ADiMU's local model discovery remains valuable for several applications. An important application case is the discovery of homogenized (surrogate) material models from strain-stress data stemming from the multi-scale simulation of heterogeneous materials representative volume elements \cite{bessa:2017a,mozaffar:2019}. Additionaly, it facilitates performance analyses when testing new material model architectures and/or exploring different types of material behavior before addressing the more challenging global model discovery scenario.
	
	The remainder of this paper demonstrates the versatility, robustness, and performance of the ADiMU framework across both local and global model discovery scenarios.
	
	\section{From strain to stress: Local direct model discovery \label{sec:local_model_discovery}}
	
	Let us start by addressing ADiMU's local direct model discovery, where a given model is discovered directly from a local strain-stress data set. Several examples of the three different types of models are discussed in the following sections, namely conventional, neural network and hybrid material models.
	
	\subsection{Conventional models}
	The simplest scenario consists in finding the parameters of a given conventional model. This task has been extensively reported in the literature, often using gradient-free optimizers. Here we leverage ADiMU's automatically differentiable, gradient-based optimization instead. Two elasto-plastic conventional models are selected for demonstrative purposes: the well-known von Mises (VM) model and the Lou-Zhang-Yoon (LZY) model (see \ref{sec:lzy_model}). The `ground-truth' parameters can be found in \ref{sec:material_parameters}. A small local strain-stress data set comprising of only 8 random polynomial strain-stress paths is generated for each model (see \ref{sec:local_datasets}).
	
	\begin{remark}
		In this paper, we employ a very tight convergence criterion to terminate the local discovery of conventional models. Convergence is only achieved when the relative change of each parameter is less than $\mathit{10^{-4}}$ for five consecutive epochs. In practice, a higher convergence tolerance can be used to achieve a reasonable solution with lower computational costs, i.e., a lower number of epochs.
	\end{remark}
	
\begin{figure}[h!]
		\centering
		\begin{subfigure}[b]{0.38\textwidth}
			\centering
			\includegraphics[width=\textwidth]{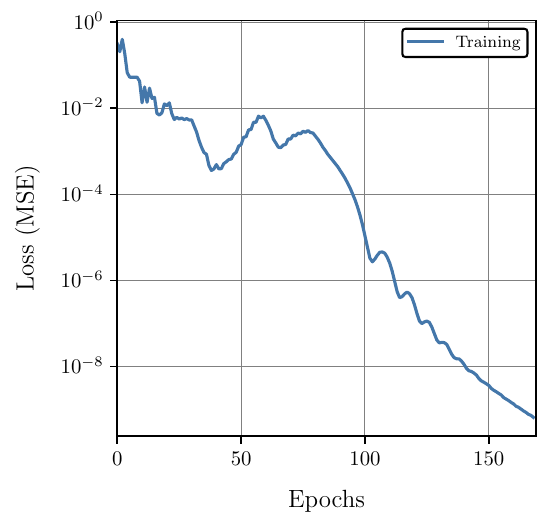}
			\caption{}
			\label{subfig:local_vm_training_loss_history}
		\end{subfigure}
		\begin{subfigure}[b]{0.38\textwidth}
			\centering
			\includegraphics[width=\textwidth]{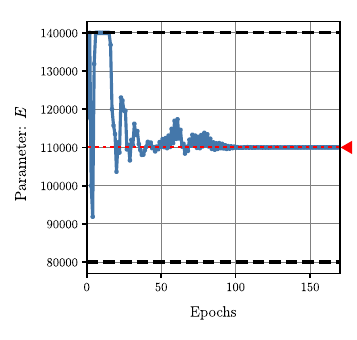}
			\caption{}
			\label{subfig:local_vm_model_parameter_history_E}
		\end{subfigure}\hfill
		\begin{subfigure}[b]{0.38\textwidth}
			\centering
			\includegraphics[width=\textwidth]{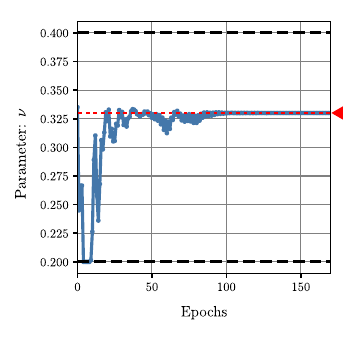}
			\caption{}
			\label{subfig:local_vm_model_parameter_history_v}
		\end{subfigure}
		\begin{subfigure}[b]{0.38\textwidth}
			\centering
			\includegraphics[width=\textwidth]{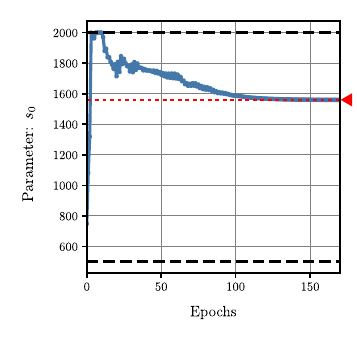}
			\caption{}
			\label{subfig:local_vm_model_parameter_history_s0}
		\end{subfigure}\hfill
		\begin{subfigure}[b]{0.38\textwidth}
			\centering
			\includegraphics[width=\textwidth]{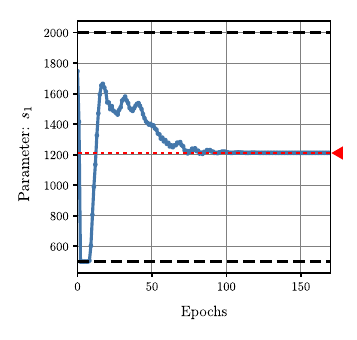}
			\caption{}
			\label{subfig:local_vm_model_parameter_history_a}
		\end{subfigure}
		\begin{subfigure}[b]{0.38\textwidth}
			\centering
			\includegraphics[width=\textwidth]{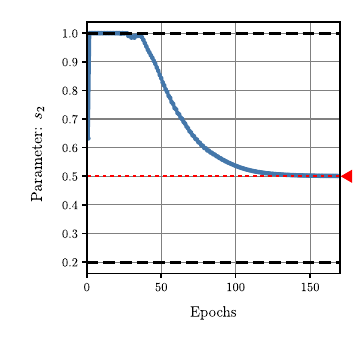}
			\caption{}
			\label{subfig:local_vm_model_parameter_history_b}
		\end{subfigure}\hfill
		\caption{VM model local discovery history from VM strain-stress data ($E=110\,$GPa, $\nu=0.33$, $s_{0}=900\sqrt{3}\,$MPa, $s_{1}=700\sqrt{3}\,$MPa, $s_{2}=0.5$) with random model initialization:
			\subref{subfig:local_vm_training_loss_history} Training loss (MSE) history throughout the discovery process;
			\subref{subfig:local_vm_model_parameter_history_E} Young modulus $E$; \subref{subfig:local_vm_model_parameter_history_v} Poisson ratio $\nu$;
			\subref{subfig:local_vm_model_parameter_history_s0} Yield parameter $s_{0}$;
			\subref{subfig:local_vm_model_parameter_history_a} Yield parameter $s_{1}$;  \subref{subfig:local_vm_model_parameter_history_b} Yield parameter $s_{2}$. Black dashed lines correspond to the optimization upper and lower bounds. Red dashed lines correspond to the parameters `ground-truth'.}
		\label{fig:vm_local_discovery}
	\end{figure}

	The first example shows the local discovery of all VM model parameters (see Figure~\ref{fig:vm_local_discovery}), namely: (i) the isotropic elastic constants ($E$, $\nu$) and (ii) the parameters of the Nadai-Ludwik isotropic strain hardening law ($s_{0}$, $s_{1}$ and $s_{2}$). Broad exploratory ranges are assigned to all parameters, particularly those associated with the hardening law that would be more difficult to find experimentally. Inspection of Figures~\ref{subfig:local_vm_model_parameter_history_E}-\ref{subfig:local_vm_model_parameter_history_b} shows that the `ground-truth' values of all parameters are successfully found after ca. 150 epochs (see Figure~\ref{subfig:local_vm_training_loss_history}). The local discovery process was repeated three times with a random parameter initialization, all yielding the same converged values despite traversing distinct paths in the loss landscape.
	
	\begin{remark}
		Despite the previous demonstration, the isotropic elastic constants ($E$, $\nu$) can be easily determined experimentally. At the same time, they play a major role in the strain-stress response of most conventional elasto-plastic material models. Therefore, to reduce the complexity of the loss landscape and improve the discovery of other parameters, the elastic constants should be excluded from the optimization problem unless necessary (or found at a stage where the specimen is only deforming elastically).
	\end{remark}

\begin{figure}[hbt]
		\centering
		\begin{subfigure}[b]{0.45\textwidth}
			\centering
			\includegraphics[width=\textwidth]{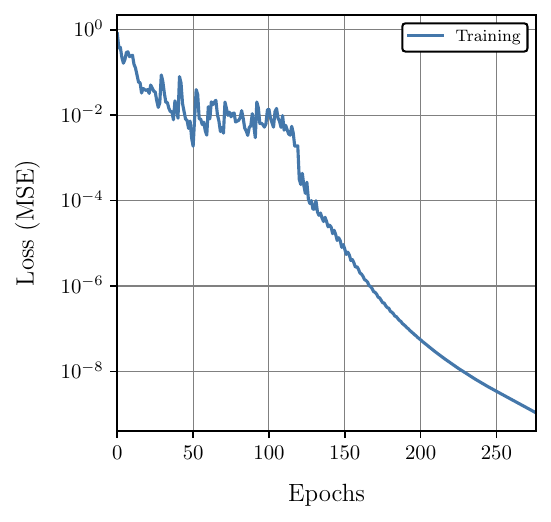}
			\caption{}
			\label{subfig:lou_local_training_loss_history}
		\end{subfigure}
		\begin{subfigure}[b]{0.447\textwidth}
			\centering
			\includegraphics[width=\textwidth]{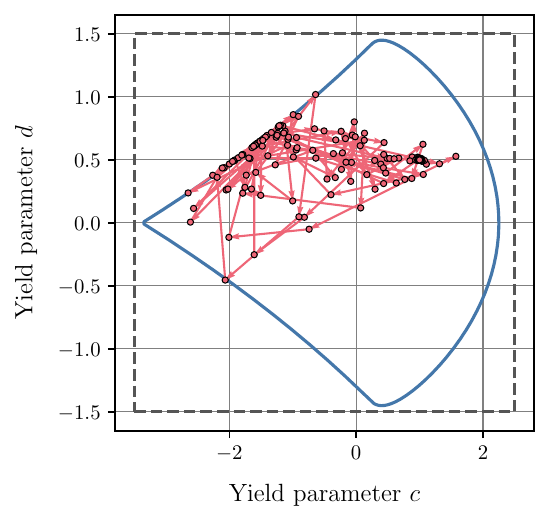}
			\caption{}
			\label{subfig:lou_model_local_convexity_history}
		\end{subfigure}\hfill
		\caption{LZY model local discovery from LZY strain-stress data ($E=110$GPa, $\nu=0.33$, $a=1.0$, $b=0.05$, $c=1.0$, $d=0.5$, $s_{0}=900\,$MPa, $s_{1}=700\,$MPa, $s_{2}=0.5$) with random model initialization:
			\subref{subfig:lou_local_training_loss_history} Training loss (MSE) history throughout the discovery process;
			\subref{subfig:lou_model_local_convexity_history} Yield surface convexity enforcement throughout the discovery process.}
		\label{fig:lou_local_discovery_loss}
	\end{figure}

\begin{figure}[h!]
		\centering
		\begin{subfigure}[b]{0.38\textwidth}
			\centering
			\includegraphics[width=\textwidth]{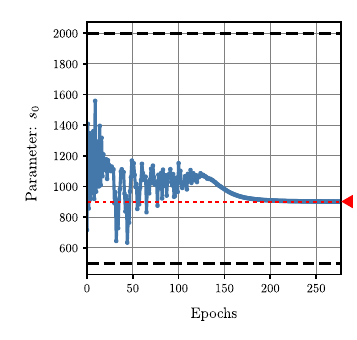}
			\caption{}
			\label{subfig:local_lou_model_parameter_history_s0}
		\end{subfigure}
		\begin{subfigure}[b]{0.38\textwidth}
			\centering
			\includegraphics[width=\textwidth]{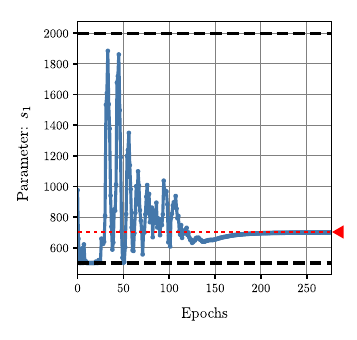}
			\caption{}
			\label{subfig:local_lou_model_parameter_history_a}
		\end{subfigure}\hfill
		\begin{subfigure}[b]{0.38\textwidth}
			\centering
			\includegraphics[width=\textwidth]{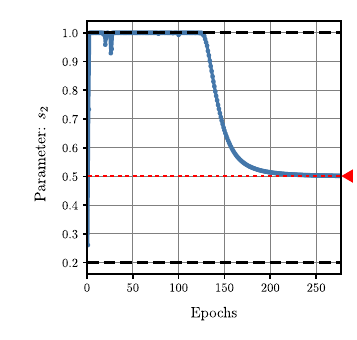}
			\caption{}
			\label{subfig:local_lou_model_parameter_history_b}
		\end{subfigure}
		\begin{subfigure}[b]{0.38\textwidth}
			\centering
			\includegraphics[width=\textwidth]{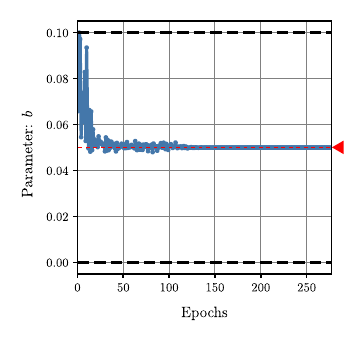}
			\caption{}
			\label{subfig:local_lou_model_parameter_history_yield_b_s0}
		\end{subfigure}\hfill
		\begin{subfigure}[b]{0.38\textwidth}
			\centering
			\includegraphics[width=\textwidth]{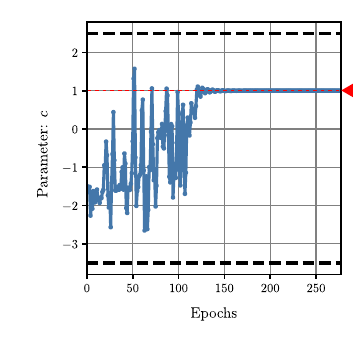}
			\caption{}
			\label{subfig:local_lou_model_parameter_history_yield_c_s0}
		\end{subfigure}
		\begin{subfigure}[b]{0.38\textwidth}
			\centering
			\includegraphics[width=\textwidth]{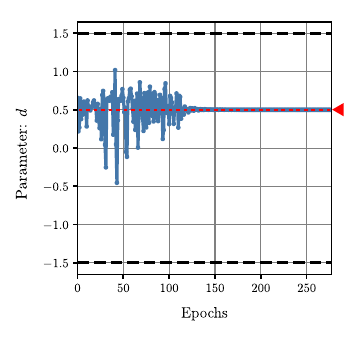}
			\caption{}
			\label{subfig:local_lou_model_parameter_history_yield_d_s0}
		\end{subfigure}\hfill
		\caption{LZY model local discovery history from LZY strain-stress data ($E=110$GPa, $\nu=0.33$, $a=1.0$, $b=0.05$, $c=1.0$, $d=0.5$, $s_{0}=900\,$MPa, $s_{1}=700\,$MPa, $s_{2}=0.5$) with random model initialization:
			\subref{subfig:local_lou_model_parameter_history_s0} Yield parameter $s_{0}$;
			\subref{subfig:local_lou_model_parameter_history_a} Yield parameter $s_{1}$;
			\subref{subfig:local_lou_model_parameter_history_b} Yield parameter $s_{2}$;
			\subref{subfig:local_lou_model_parameter_history_yield_b_s0} Yield parameter $b$;
			\subref{subfig:local_lou_model_parameter_history_yield_c_s0} Yield parameter $c$;
			\subref{subfig:local_lou_model_parameter_history_yield_d_s0} Yield parameter $d$. Black dashed lines correspond to the optimization upper and lower bounds. Red dashed lines correspond to the parameters `ground-truth'.}
		\label{fig:lou_local_discovery}
	\end{figure}

	In a second example, we turn to the more challenging local discovery of the LZY model parameters (see Figures~\ref{fig:lou_local_discovery_loss} and \ref{fig:lou_local_discovery}), namely\footnote{It is redundant to optimize the LZY yield surface parameter $a$ while simultaneously optimizing the isotropic strain hardening law parameters ($s_{0}$, $s_{1}$ and $s_{2}$). For a given value $a \neq 1$, the resulting `equivalent' parameters $s_{0}$ and $s_{1}$ are scaled accordingly (see the LZY model yield surface and hardening law in \ref{ssec:lzy_formulation}). In the particular case of the Nadai-Ludwik strain hardening law, $s_{0}$ scales linearly with $a$, whereas $s_{1}$ scales linearly with $a^{1+s_{2}}$.}: (i) the pressure-dependency ($b$), (ii) the yield surface curvature ($c$), (iii) the strength differential effect ($d$), and (iv) the parameters of the Nadai-Ludwik isotropic strain hardening law ($s_{1}$, $s_{2}$ and $s_{3}$). In addition to the hardening law parameters, a wide range of pressure dependency is taken into account (from pressure insensitivity to a friction angle of approximately $10^{\circ}$) as well as the whole yield surface convexity domain related to parameters $c$ and $d$. Figure~\ref{fig:lou_local_discovery} shows that the `ground-truth' values of all parameters are successfully found after ca. 250 epochs (see Figure~\ref{subfig:lou_local_training_loss_history}). Moreover, Figure~\ref{subfig:lou_model_local_convexity_history} shows that the yield surface convexity is enforced throughout the discovery process by means of the convexity return-mapping proposed in ~\ref{ssec:convexity_ginca}. As in the previous example, the local discovery process was conducted three times with random parameter initialization, consistently converging to the same values despite following different paths in the loss landscape.\footnote{The stability of the material model state update, with respect to both material parameters and loading conditions, is essential to perform an effective local model discovery. The former enables exploration of a broad parameter range, while the latter leverages diverse strain loading paths to enhance the conditioning of the discovery process.}
	
	Lastly, it is important to highlight that ADiMU's local model discovery of any conventional model only requires the corresponding state update algorithm used in a conventional simulation. All parameter derivatives needed for gradient-based optimization are computed via automatic differentiation, involving no additional formulation-related enhancements. Nevertheless, adhering to certain computational implementation strategies can yield significant efficiency gains, namely by leveraging implicit differentiation to handle the iterative solution of complex state update algorithms \citep{bhatia:2025a}.

	\subsection{Neural network models \label{ssec:local_ml_models}}
	
	We now focus on the case where the material model is entirely physics-agnostic, meaning it incorporates no physics-based knowledge or constraints of any kind. This implies that the local discovery process is purely data-driven, relying exclusively on the capacity of the select model architecture to learn from the local strain-stress data set. To demonstrate ADiMU's performance in this context, a multi-layer gated recurrent unit (GRU) neural network model is chosen as the material model architecture (see \ref{sec:gru_architecture}). The GRU model hyperparameters, namely the number of recurrent layers and the hidden layer size, are determined with the Tree-structure Parzen Estimator (TPE) algorithm \citep{bergstra:2011a}, yielding a model with approximately 3 million parameters.
	
	\begin{remark}
		Unlike the optimization of conventional models in the previous section, where each model involves a small set of parameters, gradient-free approaches are impracticable to discover such a GRU material model with millions of parameters. Moreover, optimizing the GRU material model does not require setting bounds on learnable parameters, as it does not involve solving a physics-based state update system of equations.
	\end{remark}
	
	The different local strain-stress data sets used in this section are generated based on a `ground-truth' LZY model, whose parameters are described in \ref{sec:material_parameters}, and consist of random polynomial strain-stress paths (see ~\ref{sec:local_datasets}). The different data set sizes used for training, validation and testing purposes are reported in Table~\ref{tab:local_gru_dataset_sizes}. Due to the physically uninterpretable parameters of the GRU architecture, all discovered GRU material models are assessed by evaluating their performance on an unseen testing data set comprising 512 random polynomial strain-stress paths.
	
\begin{figure}[h!]
		\centering
		\begin{subfigure}[b]{0.49\textwidth}
			\centering
			\includegraphics[width=\textwidth]{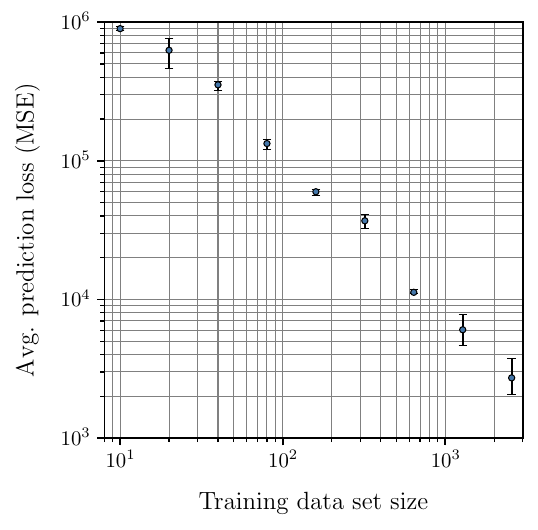}
			\caption{}
			\label{subfig:gru_conv_analysis_avg_mse}
		\end{subfigure} \vspace*{5pt} \hfill \\
		\begin{subfigure}[b]{0.49\textwidth}
			\centering
			\includegraphics[width=\textwidth]{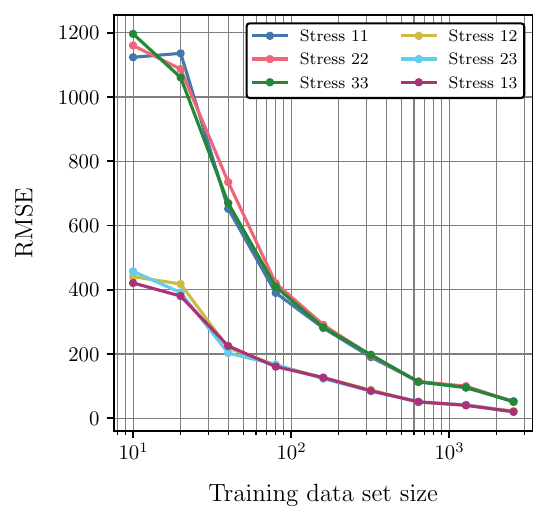}
			\caption{}
			\label{subfig:gru_conv_analysis_mean_rmse_convergence}
		\end{subfigure}
		\begin{subfigure}[b]{0.476\textwidth}
			\centering
			\includegraphics[width=\textwidth]{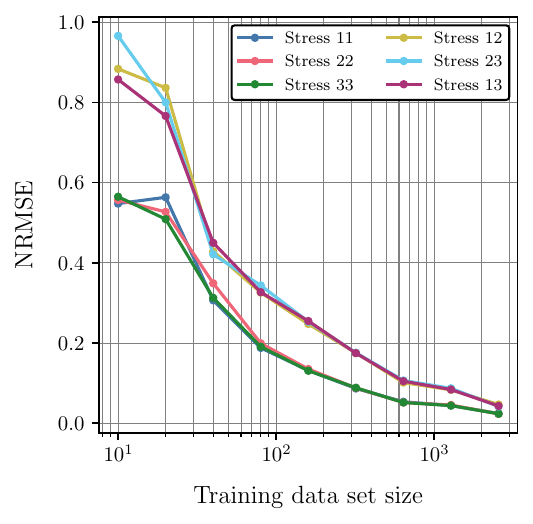}
			\caption{}
			\label{subfig:gru_conv_analysis_mean_nrmse_convergence}
		\end{subfigure}\hfill
		\caption{Performance of the GRU material model discovered from LZY strain-stress data ($E=110$GPa, $\nu=0.33$, $a=1.0$, $b=0.05$, $c=1.0$, $d=0.5$, $s_{0}=900\,$MPa, $s_{1}=700\,$MPa, $s_{2}=0.5$) and tested in a (unseen) local synthetic data set of 512 random polynomial strain-stress paths:
			\subref{subfig:gru_conv_analysis_avg_mse} Average prediction loss (MSE) with respect to the training data set size and uncertainty quantification (3 random model initializations); \subref{subfig:gru_conv_analysis_mean_rmse_convergence} Average prediction Root Mean Squared Error (RMSE) of each stress component with respect to the training data set size with random model initialization; \subref{subfig:gru_conv_analysis_mean_nrmse_convergence} Average prediction Normalized Root Mean Squared Error (NRMSE) of each stress component with respect to the training data set size with random model initialization.}
		\label{fig:gru_conv_analysis}
	\end{figure}

\begin{figure}[hbt]
		\centering
		\begin{subfigure}[b]{0.49\textwidth}
			\centering
			\includegraphics[width=\textwidth]{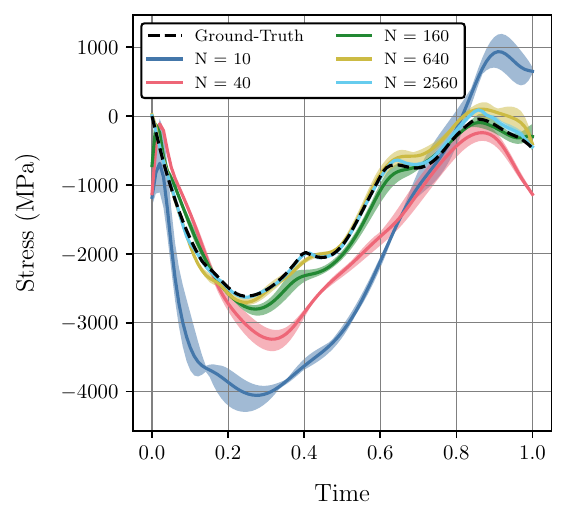}
			\caption{}
			\label{subfig:time_series_convergence_uq_stress_11_path_sample}
		\end{subfigure}
		\begin{subfigure}[b]{0.49\textwidth}
			\centering
			\includegraphics[width=\textwidth]{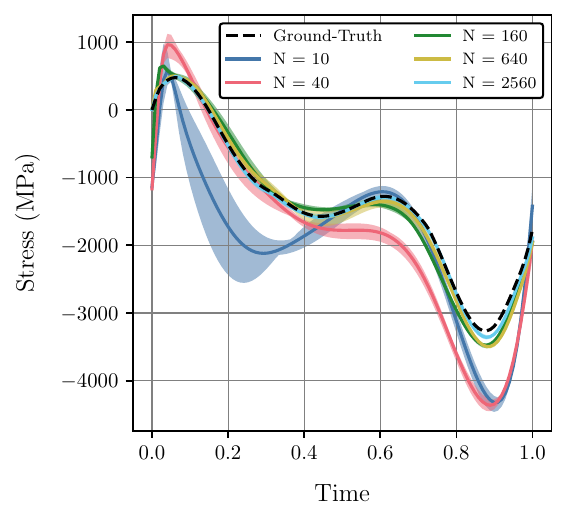}
			\caption{}
			\label{subfig:time_series_convergence_uq_stress_22_path_sample}
		\end{subfigure}\hfill
		\begin{subfigure}[b]{0.49\textwidth}
			\centering
			\includegraphics[width=\textwidth]{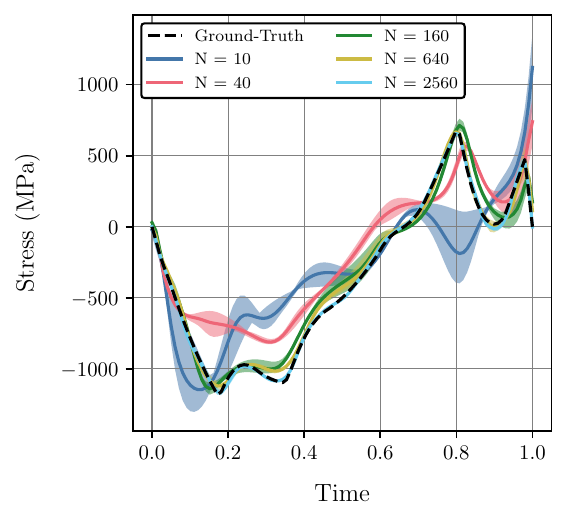}
			\caption{}
			\label{subfig:time_series_convergence_uq_stress_12_path_sample}
		\end{subfigure}
		\begin{subfigure}[b]{0.48\textwidth}
			\centering
			\includegraphics[width=\textwidth]{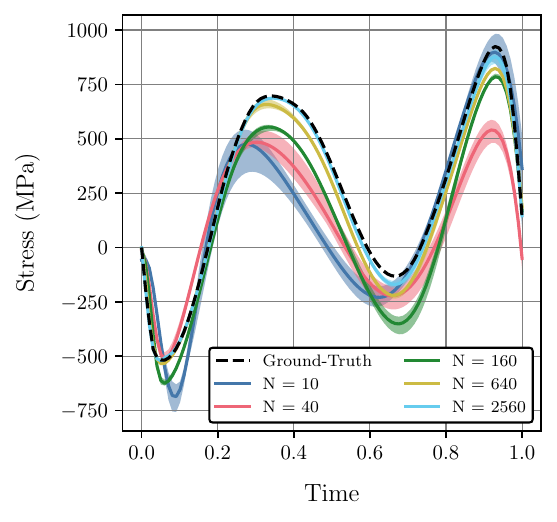}
			\caption{}
			\label{subfig:time_series_convergence_uq_stress_23_path_sample}
		\end{subfigure}\hfill
		\caption{Performance of the GRU material model discovered from LZY strain-stress data ($E=110$GPa, $\nu=0.33$, $a=1.0$, $b=0.05$, $c=1.0$, $d=0.5$, $s_{0}=900\,$MPa, $s_{1}=700\,$MPa, $s_{2}=0.5$), with uncertainty estimation (from 3 random model initializations), on randomly picked (unseen) testing random polynomial strain-stress paths and different stress components: \subref{subfig:time_series_convergence_uq_stress_11_path_sample} Normal stress 11; \subref{subfig:time_series_convergence_uq_stress_22_path_sample} Normal stress 22;
			\subref{subfig:time_series_convergence_uq_stress_12_path_sample} Shear stress 12;
			\subref{subfig:time_series_convergence_uq_stress_23_path_sample} Shear stress 23.}
		\label{fig:gru_time_series_convergence_uq_stress_paths}
	\end{figure}

	Given the data-driven nature of the GRU material model, it is instructive to conduct a convergence analysis of its performance with respect to the size of the training data set, i.e., the amount of data available to learn the material behavior. Three model realizations with a random parameter initialization are performed for each training data set size. The results are shown in Figure~\ref{fig:gru_conv_analysis} for training data sets ranging from 10 to 2560 strain-stress paths. Figure~\ref{subfig:gru_conv_analysis_avg_mse} shows the expected prediction improvement of the GRU material model with the increase of the training data set size, while also revealing a low uncertainty associated with the different model realizations. Examples of randomly picked testing samples are additionaly provided in Figure~\ref{fig:gru_time_series_convergence_uq_stress_paths} for different stress components.
	
	While the aforementioned outputs are commonly reported in the literature, we believe they offer valuable but limited insights into the effective performance of the model. First, the prediction performance is averaged over all stress components, thus hiding the actual model accuracy when predicting physically relevant components under different loading scenarios. Second, the prediction performance is often reported in stress squared (Mean Squared Error, MSE) or stress (Root Mean Squared Error, RMSE) units. This complicates the interpretation of the reported performance results, as stress magnitude is heavily influenced by both material parameters and the testing loading paths. Lastly, the magnitude of the different stress components is also highly dependent on the material behavior, thus further supporting the previous arguments. In this context, we propose and advocate for a complementary performance metric that addresses the aforementioned shortcomings. For a given stress path with $n_{t}$ time steps, the model performance when predicting each stress component is given by the Normalized Root Mean Squared Error (NRMSE) defined as
	\begin{equation}
		\mathrm{NRMSE}\, (\sigma_{ij}) = \dfrac{\mathrm{MSE}\,(\sigma_{ij})}{\mathrm{MAV}\,(\bar{\sigma}_{ij})} \, , \quad \{i, \, j\} = 1, \, \dots, \, n_\mathrm{dim} \, ,
	\end{equation}
	where $\mathrm{MSE} \, (\sigma_{ij})$ and $\mathrm{MAV} \, (\bar{\sigma}_{ij})$ denote the Mean Squared Error and the Mean Absolute Value\footnote{
		The Mean Absolute Value (MAV) is computed based on the `ground-truth' data. When the MAV of a given component is null or close to zero, a characteristic stress (e.g., the initial yield stress) can be used instead to perform the required normalization.} of the stress component $\sigma_{ij}$ over all $n_{t}$ time steps, respectively. This results in an easily interpretable, relative error for each stress component whose value does not depend on a particular set of material parameters and/or loading conditions. When reporting the performance over a given testing data set, the NRMSE of each stress component can be averaged over all stress paths.
	
	The GRU material model average prediction NRMSE for each stress component is shown in Figure~\ref{subfig:gru_conv_analysis_mean_nrmse_convergence}. The performance improvement with increasing training data set size is evident once more, but it is now apparent that the error is substantially higher in predicting the shear stress components (e.g., for a training data set size of 320 strain-stress paths, the average NRMSE of normal and shear components is around 10\% and 18\%, respectively). Note that this valuable information is hidden in Figure~\ref{subfig:gru_conv_analysis_mean_rmse_convergence}, where the RMSE of each stress component accounts for the prediction error but the normal to shear stress magnitude ratio does also depend on the particular model yield surface.

	To complete the evaluation of the GRU material model, it is useful to examine its performance when the training data contains noise. This analysis is conducted to account for the noise originating from the displacement field (experimental) measurements found in practice in the global model discovery context. Synthetic noise is thus introduced into the training data sets local random polynomial strain paths, as outlined in ~\ref{sec:local_datasets}, while the (noiseless) testing data set, consisting of 512 random polynomial strain-stress paths, remains unchanged.
	
	\begin{remark}
		The artificial noise is consistently injected into the corresponding noiseless data set of the same size. This ensures that the model performance comparison between the noiseless and noisy data scenarios is not affected by a different diversity of random strain-stress paths. Also note that in practice, after a material model is found from available data, the aim is to deploy it to perform numerical simulations. Therefore, even if the training data contains noise (e.g., experimental data), we assess the model's performance in the noiseless scenario found in numerical simulations.
	\end{remark}
	
\begin{figure}[hbt]
		\centering
		\begin{subfigure}[b]{0.45\textwidth}
			\centering
			\includegraphics[width=\textwidth]{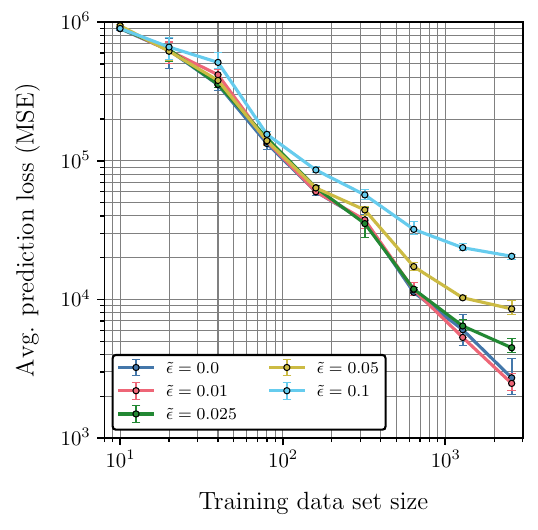}
			\caption{}
			\label{subfig:gru_conv_analysis_avg_mse_homgau_noise}
		\end{subfigure}
		\begin{subfigure}[b]{0.45\textwidth}
			\centering
			\includegraphics[width=\textwidth]{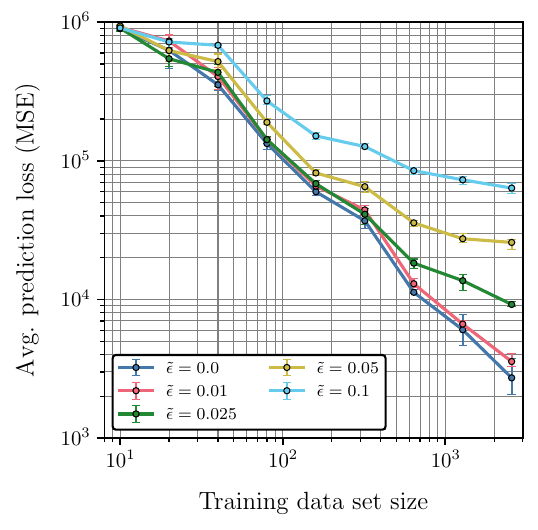}
			\caption{}
			\label{subfig:gru_conv_analysis_avg_mse_hetgau_noise}
		\end{subfigure}\hfill
		\caption{Performance of the GRU material model discovered from LZY noiseless ($\tilde{\varepsilon}=0.0$) and noisy ($\tilde{\varepsilon}>0$) strain-stress data ($E=110$GPa, $\nu=0.33$, $a=1.0$, $b=0.05$, $c=1.0$, $d=0.5$, $s_{0}=900\,$MPa, $s_{1}=700\,$MPa, $s_{2}=0.5$) and tested in a (unseen) local synthetic data set of 512 random polynomial strain-stress paths. Average prediction loss (MSE) with respect to the training data set size and uncertainty estimation (from 3 random model initializations) for different noise distributions: \subref{subfig:gru_conv_analysis_avg_mse_homgau_noise} Gaussian homoscedastic noise; \subref{subfig:gru_conv_analysis_avg_mse_hetgau_noise} Gaussian heteroscedastic noise.}
		\label{fig:gru_conv_analysis_avg_mse_hom_het_gau_noise}
	\end{figure}

	Three different noise distribution types (Gaussian, Uniform, and Spiked Gaussian), four different noise levels (see Table~\ref{tab:noise_dist_params}), and both homoscedastic and heteroscedastic noise variabilities\footnote{Homoscedastic noise means that the noise is constant for every input value. Heteroscedastic noise means that it is different at every input value.} are considered. Examples of different noisy random polynomial strain paths are shown in Figure~\ref{fig:noisy_strain_paths}. The overall performance of the GRU material model for different training data set sizes is shown in Figure~\ref{fig:gru_conv_analysis_avg_mse_hom_het_gau_noise} for both homoscedastic and heteroscedastic Gaussian noise, while the results for the remaining noise distribution types are provided in ~\ref{sec:additional_noise_local}. Three model realizations with a random parameter initialization are performed for each training data set size and noise level.
	
	Despite improving with increasing training data set size, the overall performance of the GRU material model drops in the presence of data noise, as expected. The decline in prediction accuracy worsens with higher noise levels and when transitioning from homoscedastic to heteroscedastic noise. Not surprisingly, the complex Spiked Gaussian distribution has a significantly greater impact on accuracy compared to the Gaussian and Uniform distributions. Furthermore, more challenging noise conditions disrupt the model performance at smaller training data set sizes, i.e., the noise has an increased importance when compared with the available amount of training data. At the same time, it is instructive to notice that a low noise value may actually help to regularize the model and end up improving its testing performance in comparison with the noiseless case.
	
	\clearpage
	
	\subsection{Hybrid models \label{ssec:local_hybrid_models}}
	Having addressed both conventional and neural network models, this last section on local model discovery is dedicated to hybrid models. Two types of hybrid models are addressed in what follows, namely a hybrid model architecture and a hybrid pre-trained model.
	
\begin{figure}[hbt]
		\centering
		\includegraphics[width=0.8\textwidth]{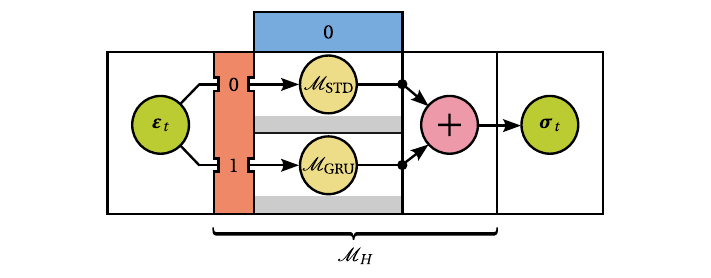}
		\caption{Hybrid model candidate-corrector type architecture. The conventional (candidate) model refers to the available physics-based model that best explains the material behavior, while the data-driven GRU (corrector) model learns the `unknown' and addresses the candidate model shortcomings.}
		\label{fig:hybrid_architecture_cc}
	\end{figure}

	A hybrid model architecture with two hybridization channels, each with a single hybridization layer, is shown in Figure~\ref{fig:hybrid_architecture_cc}. It becomes evident that this architecture is merely a particular case of ADiMU's general hybrid model architecture shown in Figure~\ref{fig:hybrid_architecture}. The first channel has a single given conventional model, such as VM or D-P models, the second channel has a single GRU model, and the hybridization model is the sum function. For a given strain path, $\bm{\varepsilon}_{t}$, $t = 0, \, \dots, \, n_{t} - 1$, the stress path, $\bm{\sigma}_{t}$, predicted by the hybrid model, $\mathcal{M}_{H}$, is thus given by
	\begin{equation}
		\bm{\sigma}_{t} = \mathcal{M}_{H}(\bm{\varepsilon}_{t}) = \mathcal{M}_{\mathrm{STD}}(\bm{\varepsilon}_{t}) + \mathcal{M}_{\mathrm{GRU}}(\bm{\varepsilon}_{t}) \, , \qquad t = 0, \, \dots, \, n_{t} - 1 \, ,
	\end{equation} 
	where $\mathcal{M}_{\mathrm{STD}}$ and $\mathcal{M}_{\mathrm{GRU}}$ denote the conventional and GRU models, respectively.
	
	This simple architecture realizes a candidate-corretor type model that can be interpreted as follows. The conventional model refers to any physics-based model deemed most suitable for describing the material behavior, thus being called candidate model. Such choice should leverage existing knowledge on the material type and take into account the model underlying constitutive assumptions. However, even if an accurate conventional model is available, it often does not fully explain the observed data. At this point, the general, data-driven GRU model comes into play with the task of learning the `unknown' and tackling the candidate model shortcomings, hence called corrector model. The performance of this hybrid model is demonstrated in what follows.
	
	\clearpage
	
\begin{figure}[h!]
		\centering
		\begin{subfigure}[b]{0.49\textwidth}
			\centering
			\includegraphics[width=\textwidth]{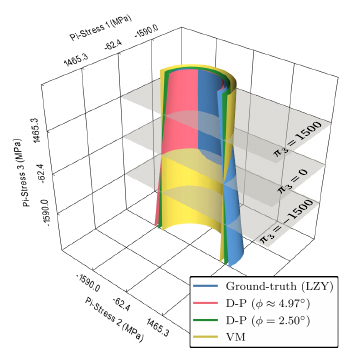}
			\caption{}
			\label{subfig:yield_surfaces_plane_cut_slices}
		\end{subfigure}
		\begin{subfigure}[b]{0.495\textwidth}
			\centering
			\includegraphics[width=\textwidth]{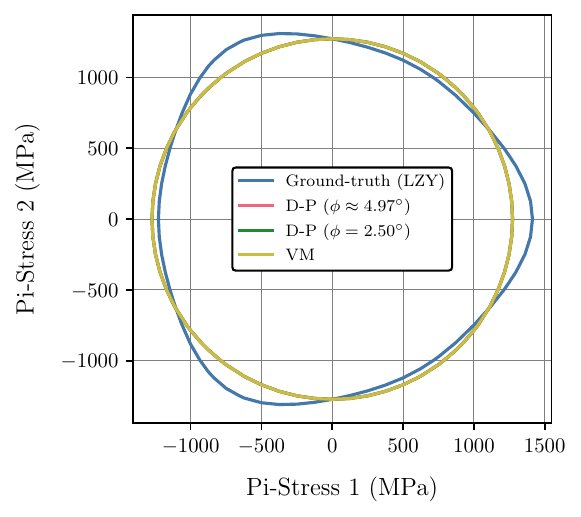}
			\caption{}
			\label{subfig:yield_surfaces_pi_plane_pressure_0}
		\end{subfigure}\vspace*{5pt} \hfill
		\begin{subfigure}[b]{0.480\textwidth}
			\centering
			\includegraphics[width=\textwidth]{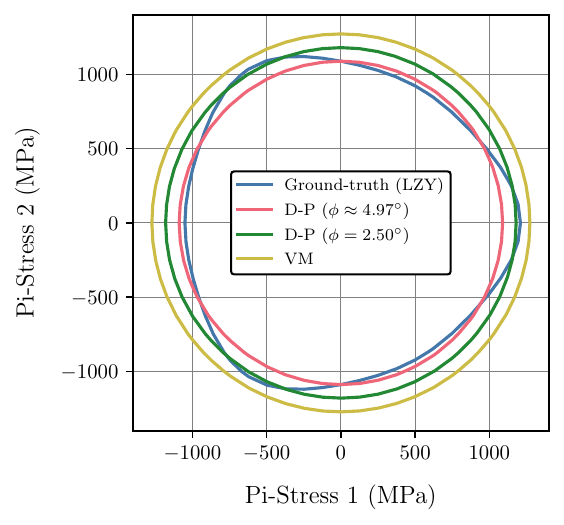}
			\caption{}
			\label{subfig:yield_surfaces_pi_plane_tensile_1500}
		\end{subfigure}
		\begin{subfigure}[b]{0.49\textwidth}
			\centering
			\includegraphics[width=\textwidth]{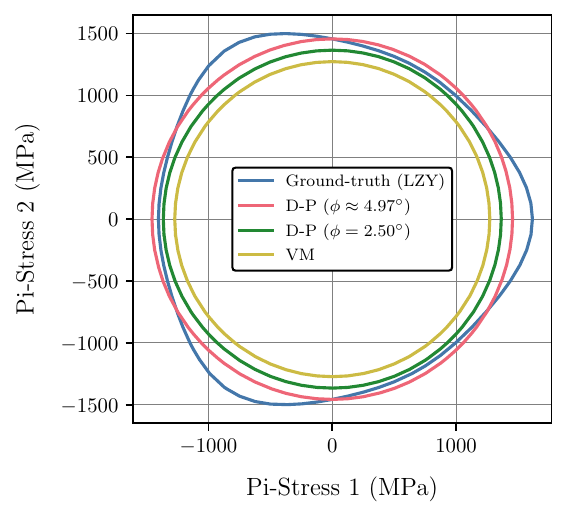}
			\caption{}
			\label{subfig:yield_surfaces_pi_plane_compression_1500}
		\end{subfigure}\hfill
		\caption{Comparison between the `ground-truth' LZY model and the different D-P and VM candidate models yield surfaces: \subref{subfig:yield_surfaces_plane_cut_slices} Symmetry-plane half-section view and three planes of constant pressure; \subref{subfig:yield_surfaces_pi_plane_pressure_0} Null pressure slice ($\pi_{3}=0$ MPa);
			\subref{subfig:yield_surfaces_pi_plane_tensile_1500} Tensile pressure slice ($\pi_{3}=1500$ MPa);
			\subref{subfig:yield_surfaces_pi_plane_compression_1500} Compressive pressure slice ($\pi_{3}=-1500$ MPa).}
		\label{fig:hybrid_yield_surfaces_comparison}
	\end{figure}

	For performance comparison purposes, consider the same LZY local strain-stress data sets as in the previous section for the purely data-driven GRU material model. In addition, let us assume that the LZY model, which is known to perfectly explain the observed data, is not available (or not known). In this context, it is instructive to illustrate the performance of the hybrid model for three different candidate models (see Figure~\ref{fig:hybrid_yield_surfaces_comparison}):
	\begin{itemize}
		\item Drucker-Prager ($\phi \approx 4.97^{\circ}$). A D-P model fitted to the observed data. The D-P model does not capture the yield surface curvature nor strength differential effect, but still manages to predict the yield pressure dependency accurately\footnote{The Drucker-Prager model parameter $\phi$ denotes the often called friction angle, associated with the pressure dependency of the yield surface.};
		\item Drucker-Prager ($\phi \approx 2.50^{\circ}$). A D-P model that does not capture the yield surface curvature nor strength differential effect, and predicts the yield pressure dependency inaccurately;
		\item von Mises. The VM model does not capture the yield surface pressure dependency, curvature nor strength differential effect.
	\end{itemize}
	These candidate models clearly show different scenarios concerning the best available candidate model, ranging from a model that has some shortcomings to a model that does not capture many important yield dependencies. Regarding the corrector model, we keep the same GRU model architecture used in the previous section. Three model realizations with a random parameter initialization are performed for each training data set size and candidate model.
	
	\begin{remark}
		For illustrative purposes, the parameters of the different candidate models are here kept fixed throughout the local discovery process of each hybrid model. Nonetheless, from a practical point of view, it seems reasonable to first discover the best candidate model that explains the data, and then find the corrector model that best addresses the remaining discrepancies. Nothing precludes, however, the simultaneous discovery of both candidate and corrector models. The comparison of these two distinct approaches is beyond the scope of this paper.
	\end{remark}

	\begin{figure}[h!]
		\centering
		\includegraphics[width=0.48\textwidth]{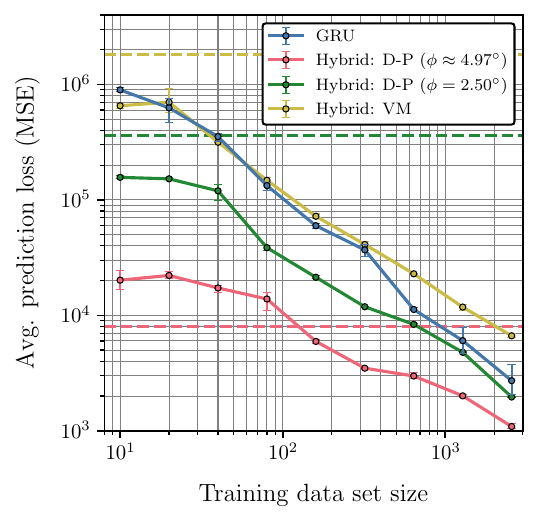}
		\caption{Performance of the GRU and three different hybrid material models discovered from LZY strain-stress data ($E=110$GPa, $\nu=0.33$, $a=1.0$, $b=0.05$, $c=1.0$, $d=0.5$, $s_{0}=900\,$MPa, $s_{1}=700\,$MPa, $s_{2}=0.5$) and tested in a (unseen) local synthetic data set of 512 random polynomial strain-stress paths. Performance is quantified by the average prediction loss (MSE) with respect to the training data set size and uncertainty quantification (3 random model initializations). Dashed reference lines correspond to the performance of the conventional candidate models.}
		\label{fig:hybrid_conv_analysis}
	\end{figure}
	
\begin{figure}[h!]
		\centering
		\begin{subfigure}[b]{0.49\textwidth}
			\centering
			\includegraphics[width=\textwidth]{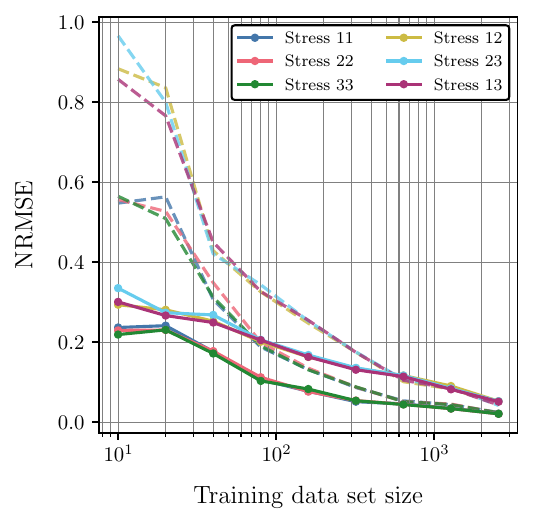}
			\caption{}
			\label{subfig:dp_2d50_plus_gru_mean_nrmse_convergence}
		\end{subfigure}
		\begin{subfigure}[b]{0.49\textwidth}
			\centering
			\includegraphics[width=\textwidth]{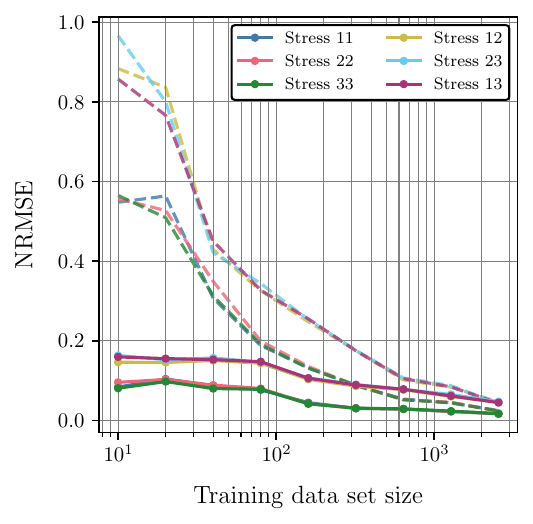}
			\caption{}
			\label{subfig:dp_4d97_plus_gru_mean_nrmse_convergence}
		\end{subfigure} \vspace*{5pt} \hfill \\
		\begin{subfigure}[b]{0.49\textwidth}
			\centering
			\includegraphics[width=\textwidth]{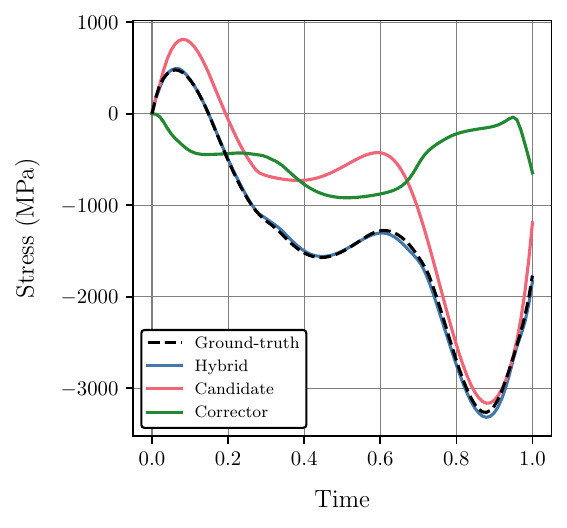}
			\caption{}
			\label{subfig:dp_2d50_plus_gru_2560_stress_22_path_sample_8}
		\end{subfigure} 
		\begin{subfigure}[b]{0.49\textwidth}
			\centering
			\includegraphics[width=\textwidth]{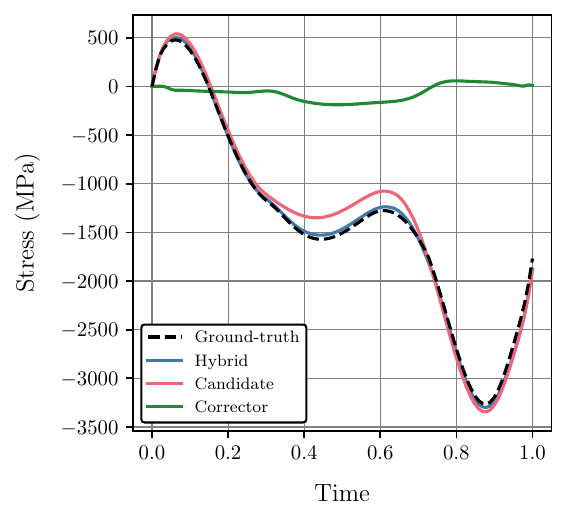}
			\caption{}
			\label{subfig:dp_4d97_plus_gru_stress_22_path_sample_8}
		\end{subfigure}\hfill
		\caption{Performance of the hybrid material models discovered from LZY strain-stress data ($E=110$GPa, $\nu=0.33$, $a=1.0$, $b=0.05$, $c=1.0$, $d=0.5$, $s_{0}=900\,$MPa, $s_{1}=700\,$MPa, $s_{2}=0.5$) and tested in a (unseen) local synthetic data set of 512 random polynomial strain-stress paths: \subref{subfig:dp_2d50_plus_gru_mean_nrmse_convergence} Hybrid model (D-P with $\phi = 2.50^{\circ}$) average prediction Normalized Root Mean Squared Error (NRMSE) of each stress component with respect to the training data set size with random model initialization;
			\subref{subfig:dp_4d97_plus_gru_mean_nrmse_convergence} Hybrid model (D-P with $\phi \approx 4.97^{\circ}$) average prediction NRMSE of each stress component with respect to the training data set size with random model initialization; \subref{subfig:dp_2d50_plus_gru_2560_stress_22_path_sample_8} Hybrid model (D-P with $\phi = 2.50^{\circ}$) prediction of randomly picked (unseen) testing sample normal stress 22; \subref{subfig:dp_4d97_plus_gru_stress_22_path_sample_8} Hybrid model (D-P with $\phi \approx 4.97^{\circ}$) prediction of randomly picked (unseen) testing sample normal stress 22. Colored dashed lines in \subref{subfig:dp_2d50_plus_gru_mean_nrmse_convergence} - \subref{subfig:dp_4d97_plus_gru_mean_nrmse_convergence} correspond to the performance of the GRU material model and are shown for comparison.}
		\label{fig:hybrid_conv_analysis_nrmse}
	\end{figure}

	Figure~\ref{fig:hybrid_conv_analysis} shows the overall performance of the three hybrid models for training data sets ranging from 10 to 2560 strain-stress paths. The performance of the previously discussed GRU material model, as well as the performance of the different conventional candidate models, are also shown for comparison. Let us start by focusing the hybrid model with the D-P candidate model ($\phi \approx 4.97^{\circ}$) that best explains the data. The performance of the hybrid model is significantly better than the pure data-driven GRU material model, regardless of the training data set size. The valuable physics-based knowledge provided by the candidate D-P model is futher highlighted in Figure~\ref{subfig:dp_4d97_plus_gru_mean_nrmse_convergence}, where the average NRMSE of the different stress components is shown. It is clearly observed that the performance improvement is inversely related to the amount of available data, i.e., the physics-based knowledge proves most valuable when the available data is scarce. In fact, it is interesting to notice that the data-driven GRU corrector model even disrupts the performance of the hybrid model in comparison with the accurate physics-based candidate model in the low end of the data size spectrum (as evidenced by the magenta dashed line that shows the average prediction loss for the D-P model only).
	
\begin{figure}[h!]
		\centering
		\begin{subfigure}[b]{0.49\textwidth}
			\centering
			\includegraphics[width=\textwidth]{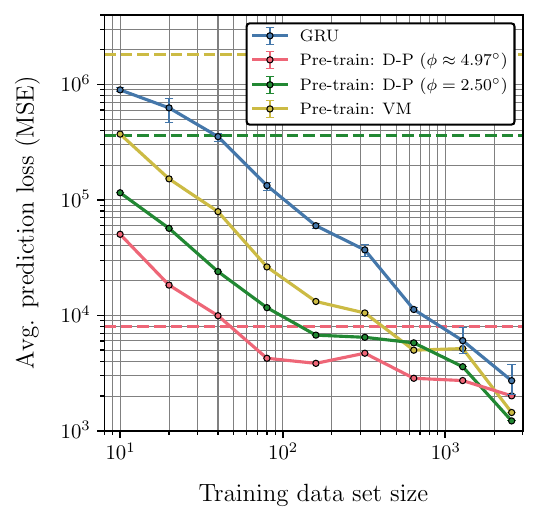}
			\caption{}
			\label{subfig:pretrain_conv_analysis_avg_mse}
		\end{subfigure}
		\begin{subfigure}[b]{0.485\textwidth}
			\centering
			\includegraphics[width=\textwidth]{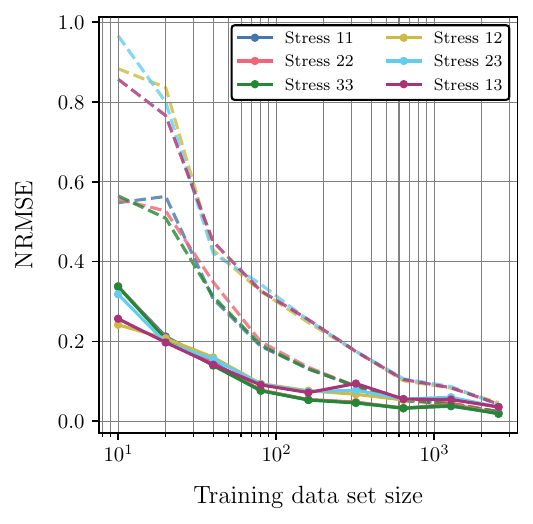}
			\caption{}
			\label{subfig:pt_gru_dp0d01_mean_nrmse_convergence}
		\end{subfigure} \vspace*{5pt} \hfill \\
		\begin{subfigure}[b]{0.485\textwidth}
			\centering
			\includegraphics[width=\textwidth]{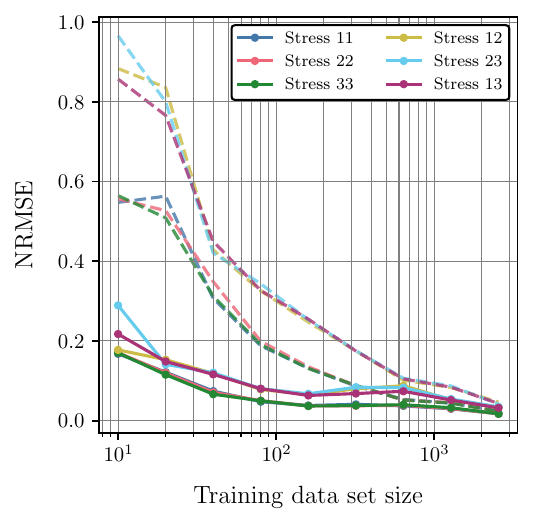}
			\caption{}
			\label{subfig:pt_gru_dp2d50_mean_nrmse_convergence}
		\end{subfigure}
		\begin{subfigure}[b]{0.485\textwidth}
			\centering
			\includegraphics[width=\textwidth]{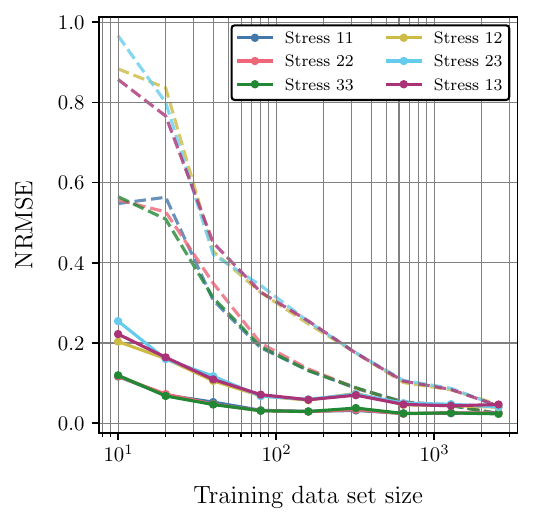}
			\caption{}
			\label{subfig:pt_gru_dp4d97_mean_nrmse_convergence}
		\end{subfigure}\hfill
		\caption{Performance of the GRU and three different hybrid pre-trained GRU material models discovered from LZY strain-stress data ($E=110$GPa, $\nu=0.33$, $a=1.0$, $b=0.05$, $c=1.0$, $d=0.5$, $s_{0}=900\,$MPa, $s_{1}=700\,$MPa, $s_{2}=0.5$) and tested in a (unseen) local synthetic data set of 512 random polynomial strain-stress paths:
			\subref{subfig:pretrain_conv_analysis_avg_mse} Average prediction loss (MSE) with respect to the training data set size and uncertainty quantification (3 random model initializations). Dashed reference lines correspond to the performance of the conventional models underlying the pre-training data; \subref{subfig:pt_gru_dp0d01_mean_nrmse_convergence} Pre-trained GRU material model (VM) average prediction Normalized Root Mean Squared Error (NRMSE) of each stress component with respect to the training data set size with random model initialization; \subref{subfig:pt_gru_dp2d50_mean_nrmse_convergence} Pre-trained GRU material model (D-P with $\phi = 2.50^{\circ}$) NRMSE of each stress component with respect to the training data set size with random model initialization; \subref{subfig:pt_gru_dp4d97_mean_nrmse_convergence} Pre-trained GRU material model (D-P with $\phi \approx 4.97^{\circ}$) average prediction NRMSE of each stress component with respect to the training data set size with random model initialization. Colored dashed lines in \subref{subfig:pt_gru_dp0d01_mean_nrmse_convergence} - \subref{subfig:pt_gru_dp4d97_mean_nrmse_convergence} correspond to the performance of the (non-pretrained) GRU material model and are shown for comparison.}
		\label{fig:pretrain_conv_analysis}
	\end{figure}

	Furthermore, we also show the additive candidate-corrector nature of this particular hybrid model architecture (with D-P considering $\phi = 4.97^{\circ}$) in Figure~\ref{subfig:pt_gru_dp4d97_mean_nrmse_convergence} for a randomly picked testing stress path. The candidate model provides a solid baseline prediction, but the corrector model further enhances accuracy where the candidate model falls short. In turn, despite the inaccurate prediction of the yield pressure dependency, the hybrid model with the D-P candidate model ($\phi = 2.5^{\circ}$) still outperforms the physics-agnostic GRU material model throughout the whole training data set size spectrum. When compared to the conventional candidate model, a modest improvement is achieved in the low data regime, but becomes very significant as more data is available. Given that the baseline prediction provided by the candidate model is worse, the hybrid model requires more data than the previous case to achieve the same performance level, as expected.
	This is evident from Figure~\ref{subfig:dp_2d50_plus_gru_2560_stress_22_path_sample_8}, which highlights that the data-driven GRU corrector model needs to provide a more difficult correction when compared to Figure~\ref{subfig:dp_4d97_plus_gru_mean_nrmse_convergence}.
	
	Lastly, observing the performance of the hybrid model with the VM candidate model leads to an interesting conclusion: if the chosen conventional (physics-based) model is too far from the ground-truth material behavior, then the correction by the GRU corrector can be more difficult than simply learning directly from the data without doing the hybridization. In other words, despite the significant improvement when compared with the conventional candidate model, the hybrid model does not outperform the pure data-driven GRU material model.
	
	Alternatively, we can also infuse knowledge by pre-training the GRU on a chosen conventional model. Rather than integrating the conventional candidate model into the hybrid model architecture, it can be used to create a large data set for pre-training the GRU material model. As a consequence, the discovery process of the actual material behavior begins with a more informed initial state of the GRU material model, commonly referred to as a `warm start', aiming to reduce the required amount of data to achieve a given prediction performance.
	
	\begin{remark}
		Despite its simplicity, this hybrid modeling approach addresses the typical scenario where high-fidelity material data is expensive and limited. Similar to multi-fidelity modeling approaches \cite{yi2024practical}, it explores low-fidelity data that can be easily acquired or generated to enhance the overall model discovery process. In this particular case, we assume that a reasonable conventional candidate model is available to efficiently generate a large low-fidelity data set.
	\end{remark}
	
	Each of the three different candidate models (see Figure~\ref{fig:hybrid_yield_surfaces_comparison}) is thus used to generate a local strain-stress data set consisting of 2560 random polynomial strain-stress paths. Accordingly, three hybrid pre-trained models are initialized by pre-training a GRU material model on each of those local data sets, with the same architecture previously described. In the following analysis, the same LZY local strain-stress data sets previously used are considered. The performance of the three hybrid pre-trained GRU material models is shown in Figure~\ref{fig:pretrain_conv_analysis}, where the GRU material model discussed in the previous section is also included for comparison, together with the performance of the different conventional candidate models. First, it is observed that all hybrid pre-trained models outperform the GRU material model, even the one pre-trained with the VM candidate model. In the second place, the performance improvement is clearly dependent on the quality of the pre-training data, increasing with the accuracy of the underlying conventional candidate model, as expect. Lastly, the pre-training effect decreases with the increase of the training data set size, as the amount of the actual material data enhances and dominates the discovery process.

	\section{From displacements to forces: Global indirect model discovery \label{sec:global_model_discovery}}
	
	We now turn to ADiMU's global indirect model discovery, where a model is indirectly discovered from a displacement-force data set. As in the previous section, several examples demonstrate ADiMU's capability of handling conventional, neural network and hybrid material models.
	
	\subsection{Conventional models}
	
	The first scenario involves discovering the parameters of an elasto-plastic conventional model using displacement-force data collected from the well-known uniaxial tensile test of a dogbone specimen. The specimen geometry and loading conditions, illustrated in Figure~\ref{fig:dogbone_specimen}, are detailed in ~\ref{ssec:tensile_dogbone}. Two elasto-plastic conventional models are selected for demonstrative purposes: the von Mises (VM) model and the Drucker-Prager (D-P) model. The `ground-truth' parameters can be found in \ref{sec:material_parameters}.
	
\begin{figure}[h!]
		\centering
		\begin{subfigure}[b]{0.45\textwidth}
			\centering
			\includegraphics[width=\textwidth]{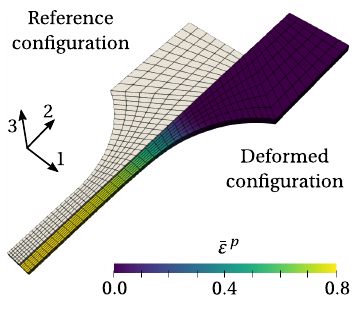}
			\caption{}
			\label{subfig:vm_tensile_dogbone_acc_p_strain}
		\end{subfigure}
		\begin{subfigure}[b]{0.42\textwidth}
			\centering
			\includegraphics[width=\textwidth]{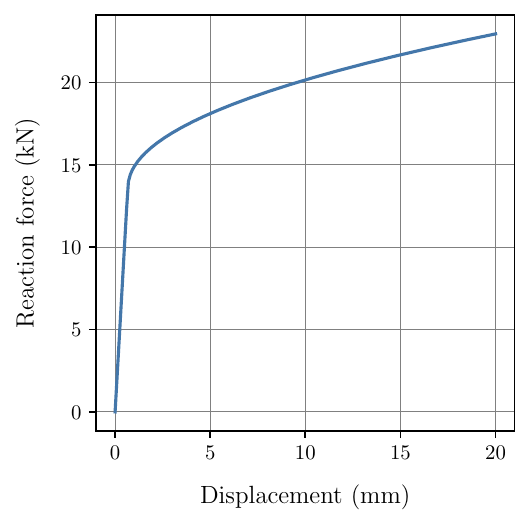}
			\caption{}
			\label{subfig:vm_tensile_dogbone_force_displacement}
		\end{subfigure}\hfill
		\begin{subfigure}[b]{0.45\textwidth}
			\centering
			\includegraphics[width=\textwidth]{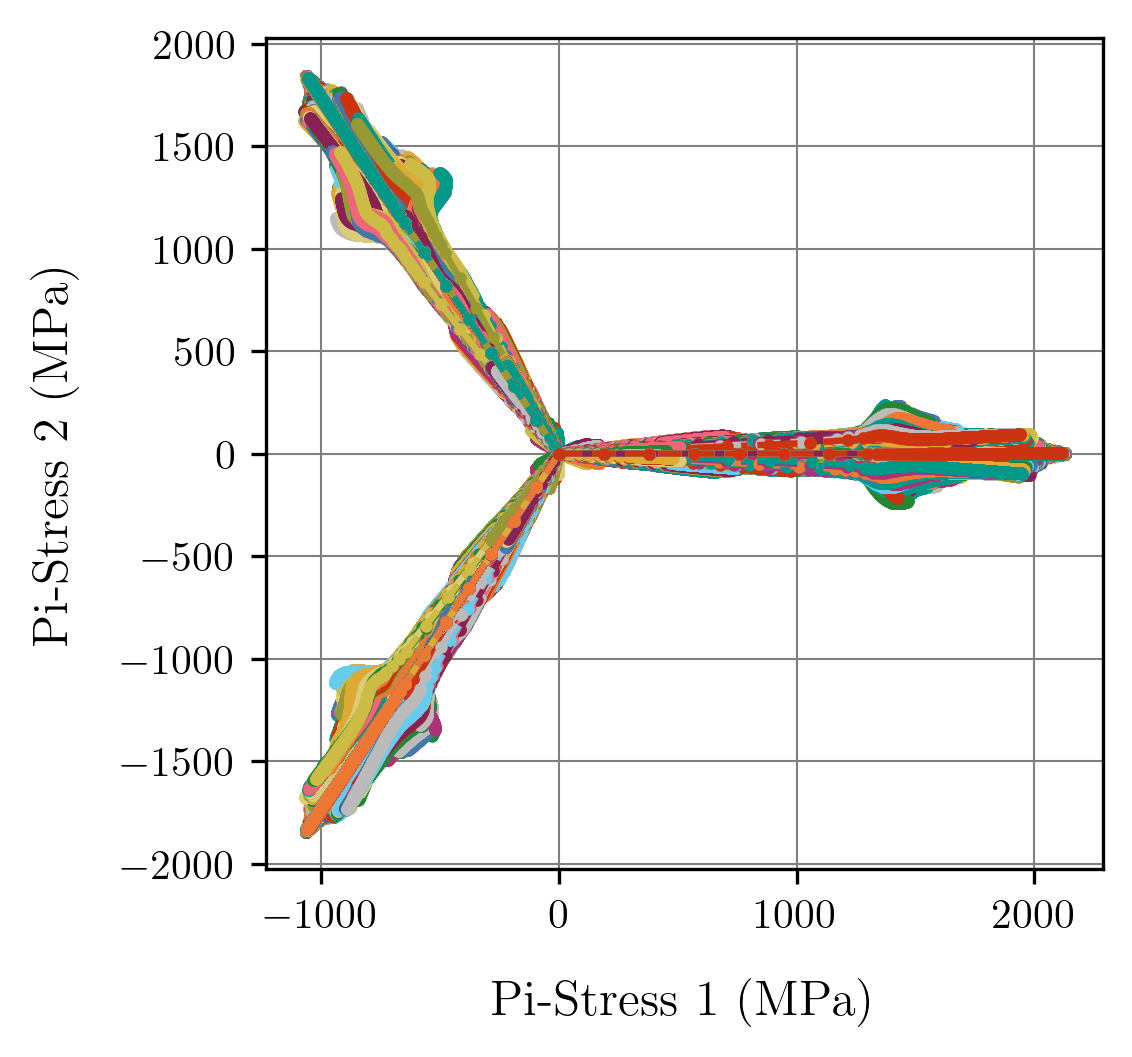}
			\caption{}
			\label{subfig:vm_tensile_dogbone_stress_path_pi_stress_1v2}
		\end{subfigure}
		\begin{subfigure}[b]{0.465\textwidth}
			\centering
			\includegraphics[width=\textwidth]{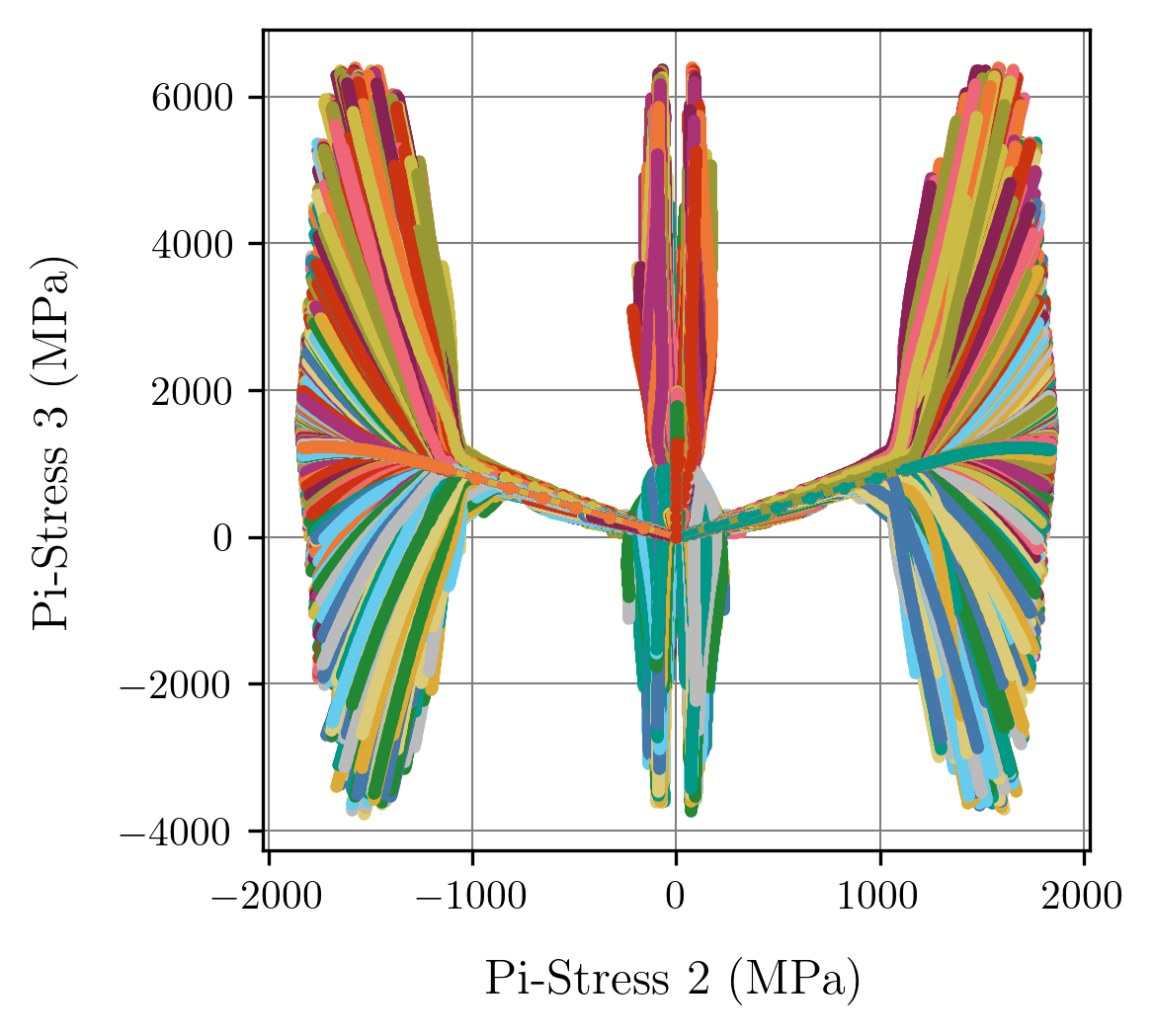}
			\caption{}
			\label{subfig:vm_tensile_dogbone_stress_path_pi_stress_2v3}
		\end{subfigure}\hfill
		\caption{Tensile dogbone VM data ($E=110\,$GPa, $\nu=0.33$, $s_{0}=900\sqrt{3}\,$MPa, $s_{1}=700\sqrt{3}\,$MPa, $s_{2}=0.5$): \subref{subfig:vm_tensile_dogbone_acc_p_strain} Accumulated plastic strain field resulting from the total prescribed displacement; \subref{subfig:vm_tensile_dogbone_force_displacement} Uniaxial force-displacement response;
			\subref{subfig:vm_tensile_dogbone_stress_path_pi_stress_1v2} $\pi_{1}-\pi_{2}$ stress projection view of the induced local strain-stress paths;
			\subref{subfig:vm_tensile_dogbone_stress_path_pi_stress_2v3} $\pi_{2}-\pi_{3}$ stress projection view of the induced local strain-stress paths.}
		\label{fig:vm_tensile_dogbone_training}
	\end{figure}

	In the first example, the deformed configuration and uniaxial force-displacement response\footnote{Note that the displacement shown in the specimen's force-displacement plot corresponds to the uniaxial displacement imposed by the tensile testing machine. ADiMU requires the full displacement field to perform the global material model discovery.} of the tensile dogbone specimen with a VM material are shown in Figures~\ref{subfig:vm_tensile_dogbone_acc_p_strain} and \ref{subfig:vm_tensile_dogbone_force_displacement}, respectively. Further insights are provided in Figures~\ref{subfig:vm_tensile_dogbone_stress_path_pi_stress_1v2} and ~\ref{subfig:vm_tensile_dogbone_stress_path_pi_stress_2v3}, where the induced local strain-stress paths are shown from two different stress projection views. The `ground-truth' VM material model cylindrical yield surface (see Figure~\ref{subfig:vm_gt_yield_surface} for reference) can be clearly identified under uniaxial stress conditions, as expected.
	
	\begin{remark}
		It is emphasized that the specimen's local strain-stress paths are solely shown for illustrative purposes, as such data is not available in the global indirect model discovery. These visualizations offer a valuable insight into the diversity of the induced local strain-stress paths that ultimately condition the indirect model discovery through the resulting forces.
	\end{remark}
	
\begin{figure}[h!]
		\centering
		\begin{subfigure}[b]{0.38\textwidth}
			\centering
			\includegraphics[width=\textwidth]{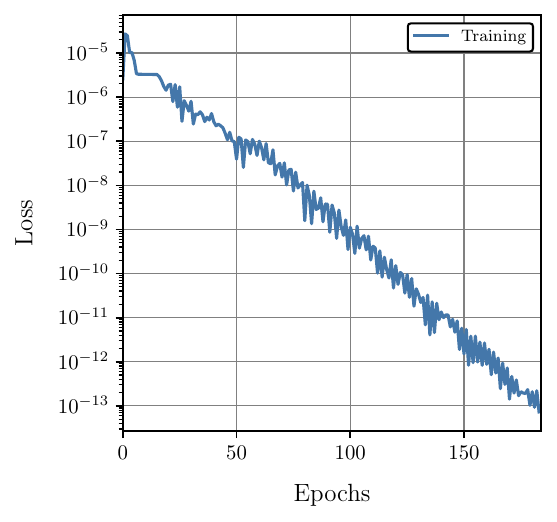}
			\caption{}
			\label{subfig:vm_tensile_dogbone_training_loss_history}
		\end{subfigure}
		\begin{subfigure}[b]{0.38\textwidth}
			\centering
			\includegraphics[width=\textwidth]{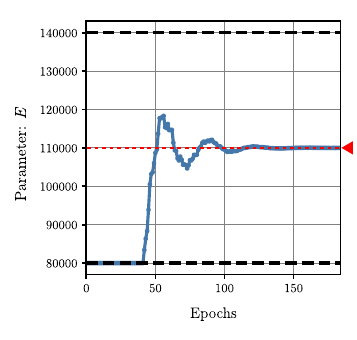}
			\caption{}
			\label{subfig:vm_model_parameter_history_E}
		\end{subfigure}\hfill
		\begin{subfigure}[b]{0.35\textwidth}
			\centering
			\includegraphics[width=\textwidth]{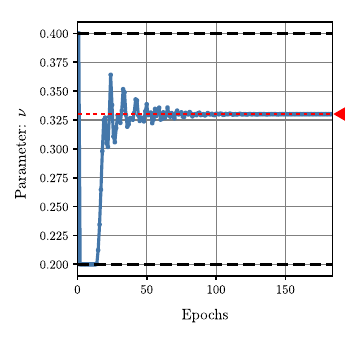}
			\caption{}
			\label{subfig:vm_model_parameter_history_v}
		\end{subfigure}
		\begin{subfigure}[b]{0.38\textwidth}
			\centering
			\includegraphics[width=\textwidth]{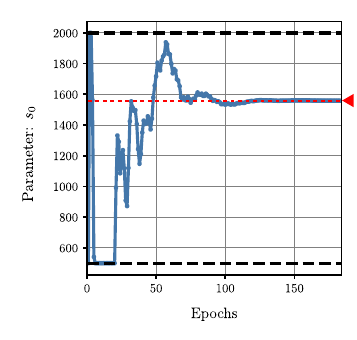}
			\caption{}
			\label{subfig:vm_model_parameter_history_s0}
		\end{subfigure}\hfill
		\begin{subfigure}[b]{0.38\textwidth}
			\centering
			\includegraphics[width=\textwidth]{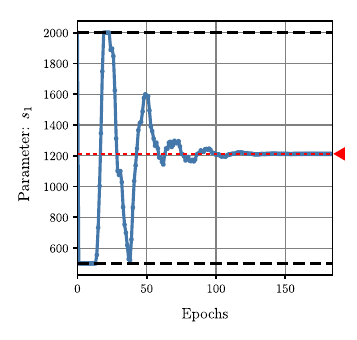}
			\caption{}
			\label{subfig:vm_model_parameter_history_a}
		\end{subfigure}
		\begin{subfigure}[b]{0.38\textwidth}
			\centering
			\includegraphics[width=\textwidth]{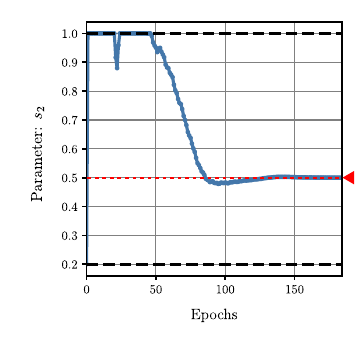}
			\caption{}
			\label{subfig:vm_model_parameter_history_b}
		\end{subfigure}
		\caption{VM model global discovery from VM displacement-force data ($E=110\,$GPa, $\nu=0.33$, $s_{0}=900\sqrt{3}\,$MPa, $s_{1}=700\sqrt{3}\,$MPa, $s_{2}=0.5$) with random model initialization:
			\subref{subfig:vm_tensile_dogbone_training_loss_history} Force equilibrium loss history throughout the discovery process;
			\subref{subfig:vm_model_parameter_history_E} Young modulus $E$; \subref{subfig:vm_model_parameter_history_v} Poisson ratio $\nu$;
			\subref{subfig:vm_model_parameter_history_s0} Yield parameter $s_{0}$;
			\subref{subfig:vm_model_parameter_history_a} Yield parameter $s_{1}$;  \subref{subfig:vm_model_parameter_history_b} Yield parameter $s_{2}$. Black dashed lines correspond to the optimization upper and lower bounds. Red dashed lines correspond to the parameters `ground-truth'.}
		\label{fig:vm_tensile_dogbone_parameters}
	\end{figure}

	The global discovery of all VM model parameters is shown in Figure~\ref{fig:vm_tensile_dogbone_parameters}. The same broad, exploratory ranges considered in the local discovery setting are assigned to all parameters. The isotropic elastic constants ($E$, $\nu$) are only included in this first example for demonstrative purposes, as they can be easily found experimentally and hence be excluded from the optimization problem. Figures~\ref{subfig:vm_model_parameter_history_E}-\ref{subfig:vm_model_parameter_history_b} show that the `ground-truth' values of all parameters are successfully found after ca. 200 epochs (see Figure~\ref{subfig:vm_tensile_dogbone_training_loss_history}). The elastic constants and the nonlinear hardening law, often required by elasto-plastic models, are thus found solely from displacement-force data.
	
\begin{figure}[h!]
		\centering
		\begin{subfigure}[b]{0.45\textwidth}
			\centering
			\includegraphics[width=\textwidth]{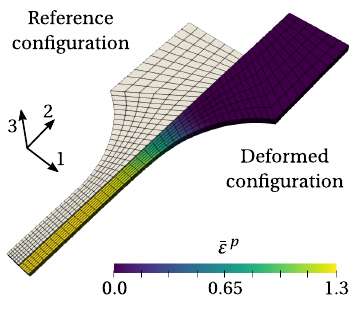}
			\caption{}
			\label{subfig:dp_tensile_dogbone_acc_p_strain}
		\end{subfigure}
		\begin{subfigure}[b]{0.42\textwidth}
			\centering
			\includegraphics[width=\textwidth]{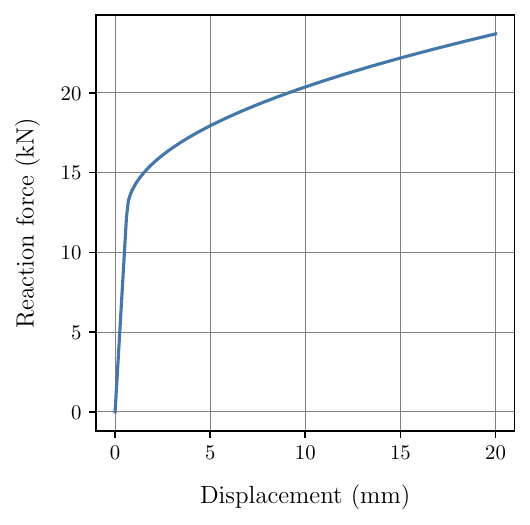}
			\caption{}
			\label{subfig:dp_tensile_dogbone_force_displacement}
		\end{subfigure}\hfill
		\begin{subfigure}[b]{0.45\textwidth}
			\centering
			\includegraphics[width=\textwidth]{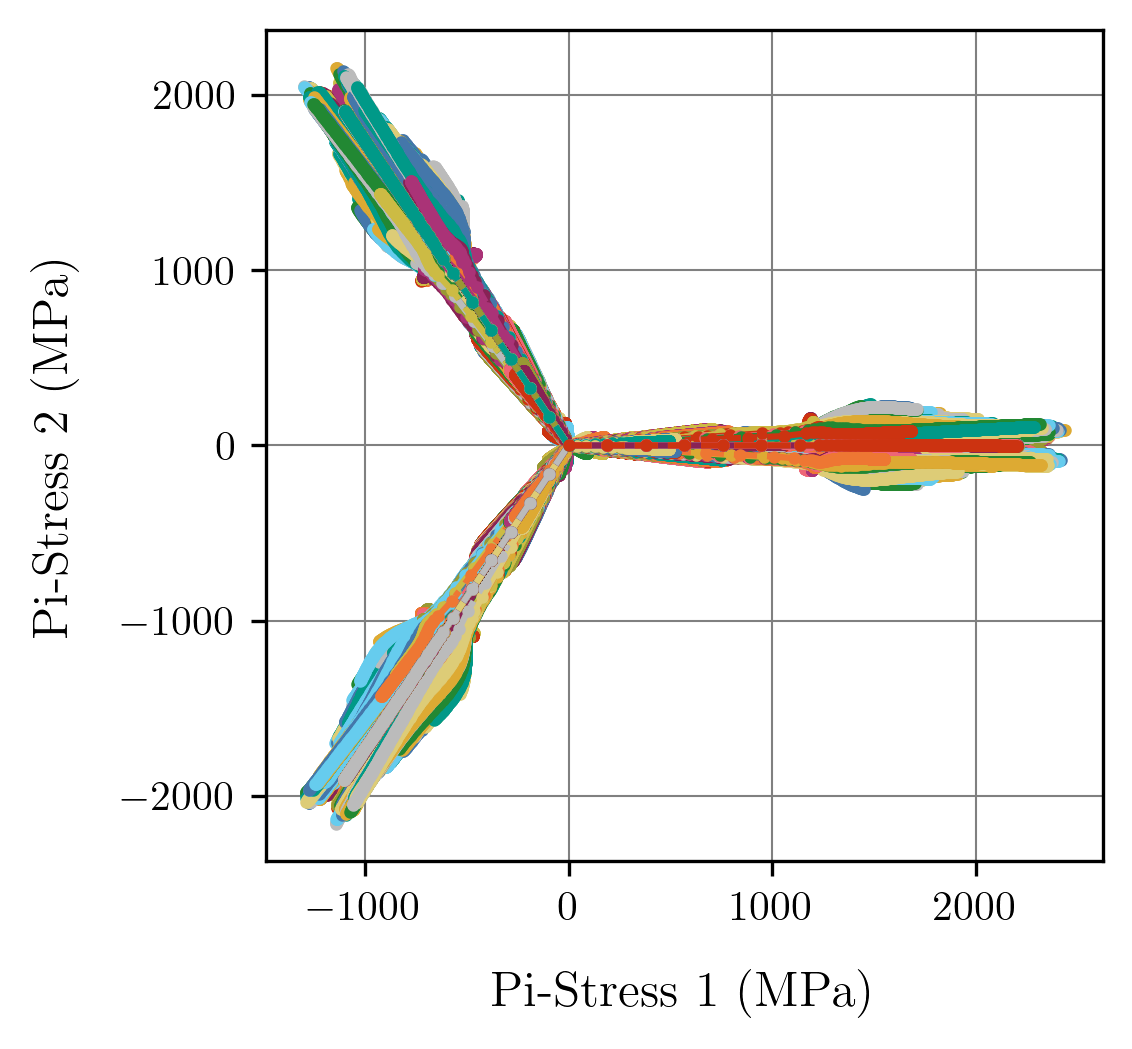}
			\caption{}
			\label{subfig:dp_tensile_dogbone_stress_path_pi_stress_1v2}
		\end{subfigure}
		\begin{subfigure}[b]{0.45\textwidth}
			\centering
			\includegraphics[width=\textwidth]{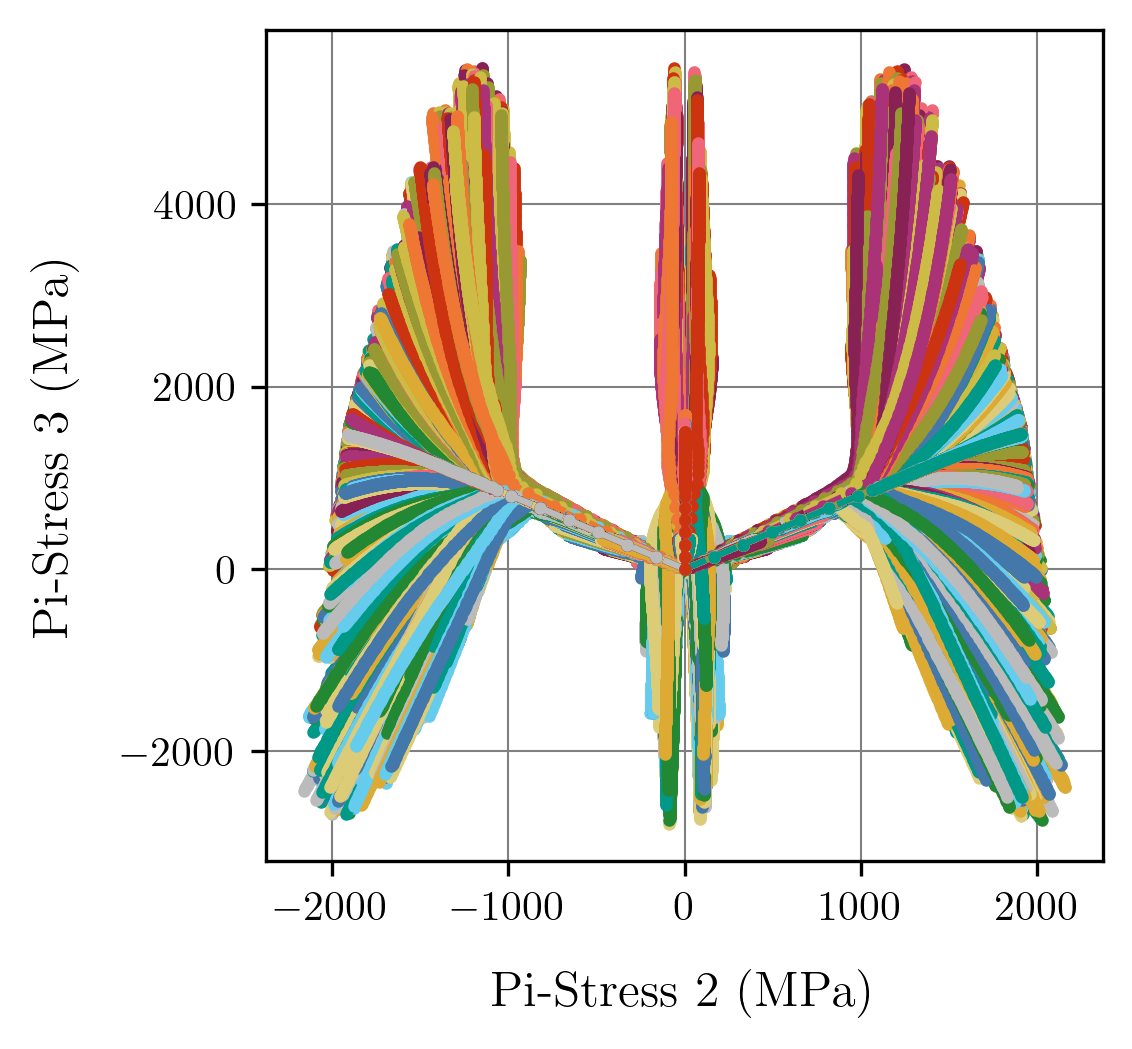}
			\caption{}
			\label{subfig:dp_tensile_dogbone_stress_path_pi_stress_2v3}
		\end{subfigure}\hfill
		\caption{Tensile dogbone D-P data ($E=110$GPa, $\nu=0.33$, $\phi=5^{\circ}\approx 0.0873 \, \mathrm{rad}$, $s_{0}=900/\xi\,$MPa, $s_{1}=700/\xi\,$MPa, $s_{2}=0.5$): \subref{subfig:dp_tensile_dogbone_acc_p_strain} Accumulated plastic strain field resulting from the total prescribed displacement; \subref{subfig:dp_tensile_dogbone_force_displacement} Uniaxial force-displacement response;
			\subref{subfig:dp_tensile_dogbone_stress_path_pi_stress_1v2} $\pi_{1}-\pi_{2}$ stress projection view of the induced local strain-stress paths;
			\subref{subfig:dp_tensile_dogbone_stress_path_pi_stress_2v3} $\pi_{2}-\pi_{3}$ stress projection view of the induced local strain-stress paths.}
		\label{fig:dp_tensile_dogbone_training}
	\end{figure}

\begin{figure}[H]
		\centering
		\begin{subfigure}[b]{0.38\textwidth}
			\centering
			\includegraphics[width=\textwidth]{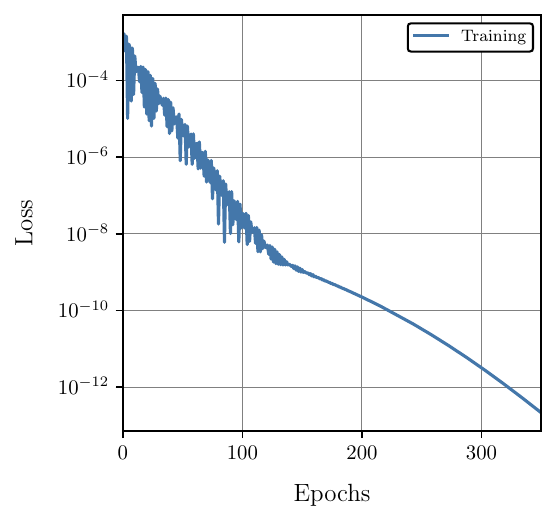}
			\caption{}
			\label{subfig:dp_tensile_dogbone_training_loss_history}
		\end{subfigure}
		\begin{subfigure}[b]{0.38\textwidth}
			\centering
			\includegraphics[width=\textwidth]{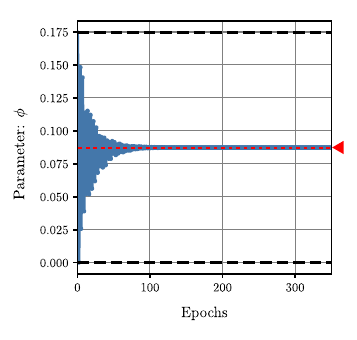}
			\caption{}
			\label{subfig:dp_model_parameter_history_friction_angle}
		\end{subfigure}\hfill
		\begin{subfigure}[b]{0.38\textwidth}
			\centering
			\includegraphics[width=\textwidth]{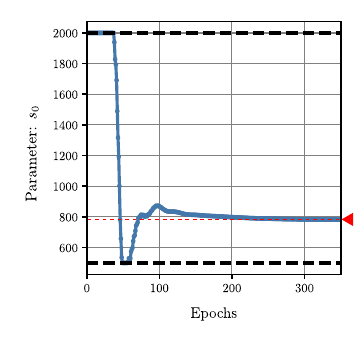}
			\caption{}
			\label{subfig:dp_model_parameter_history_s0}
		\end{subfigure}
		\begin{subfigure}[b]{0.38\textwidth}
			\centering
			\includegraphics[width=\textwidth]{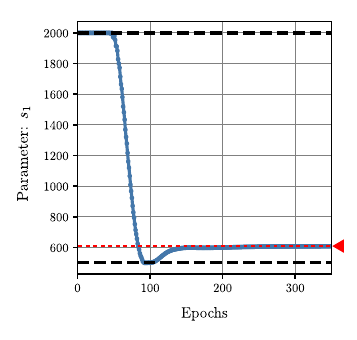}
			\caption{}
			\label{subfig:dp_model_parameter_history_a}
		\end{subfigure}\hfill
		\begin{subfigure}[b]{0.38\textwidth}
			\centering
			\includegraphics[width=\textwidth]{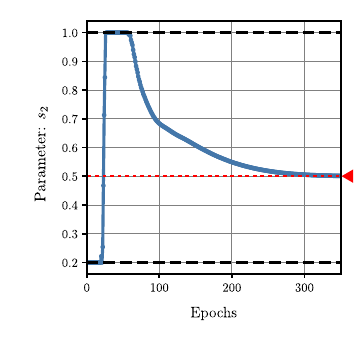}
			\caption{}
			\label{subfig:dp_model_parameter_history_b}
		\end{subfigure}
		\caption{D-P model global discovery from D-P displacement-force data ($E=110$GPa, $\nu=0.33$, $\phi=5^{\circ}\approx 0.0873 \, \mathrm{rad}$, $s_{0}=900/\xi\,$MPa, $s_{1}=700/\xi\,$MPa, $s_{2}=0.5$) with random model initialization:
			\subref{subfig:dp_tensile_dogbone_training_loss_history} Force equilibrium loss history throughout the discovery process;
			\subref{subfig:dp_model_parameter_history_friction_angle} Friction angle $\phi$;
			\subref{subfig:dp_model_parameter_history_s0} Yield parameter $s_{0}$;
			\subref{subfig:dp_model_parameter_history_a} Yield parameter $s_{1}$;
			\subref{subfig:dp_model_parameter_history_b} Yield parameter $s_{2}$. Black dashed lines correspond to the optimization upper and lower bounds. Red dashed lines correspond to the parameters `ground-truth'.}
		\label{fig:dp_tensile_dogbone_parameters}
	\end{figure}

	In the second example, let us consider the same uniaxial tensile test of a dogbone specimen but now assuming a D-P material (see Figures~\ref{subfig:dp_tensile_dogbone_acc_p_strain} and \ref{subfig:dp_tensile_dogbone_force_displacement}). Despite the uniaxial nature of this tensile test, note how valuable yield pressure dependency data is still available from the specimen local strain-stress paths shown in Figure~\ref{subfig:dp_tensile_dogbone_stress_path_pi_stress_2v3}, where the 
	`ground-truth' D-P material model conical yield surface is evidenced (see Figure~\ref{subfig:dp_gt_yield_surface}). It is thus understandable that the D-P model friction angle is successfully found, together with all the hardening law parameters, in the global discovery process shown in Figure~\ref{fig:dp_tensile_dogbone_parameters}.

	\begin{remark}
		It is important to highlight that the yield pressure dependency plays a major role in the discovery of the elasto-plastic model parameters, namely in the global indirect context. From an optimization standpoint, several numerical examples suggest that the remaining parameters begin converging only after the `ground-truth' yield pressure dependency is reasonably approached.
	\end{remark}
	
	The third and last example involves discovering the parameters of the elasto-plastic Lou-Zhang-Yoon (LZY) conventional model (see ~\ref{sec:lzy_model}) from the displacement-force data collected in the uniaxial tensile test of a double notched specimen (see Figure~\ref{fig:double_notched_specimen}). The specimen geometry and loading conditions are detailed in ~\ref{ssec:tensile_double_notched}, while the LZY `ground-truth' parameters can be found in ~\ref{sec:material_parameters}.
	
	\begin{figure}[h!]
		\centering
		\begin{subfigure}[b]{0.45\textwidth}
			\centering
			\includegraphics[width=\textwidth]{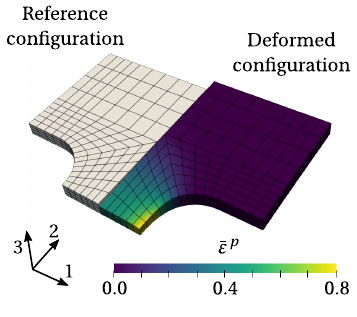}
			\caption{}
			\label{subfig:lou_tensile_double_notched_acc_p_strain}
		\end{subfigure}
		\begin{subfigure}[b]{0.42\textwidth}
			\centering
			\includegraphics[width=\textwidth]{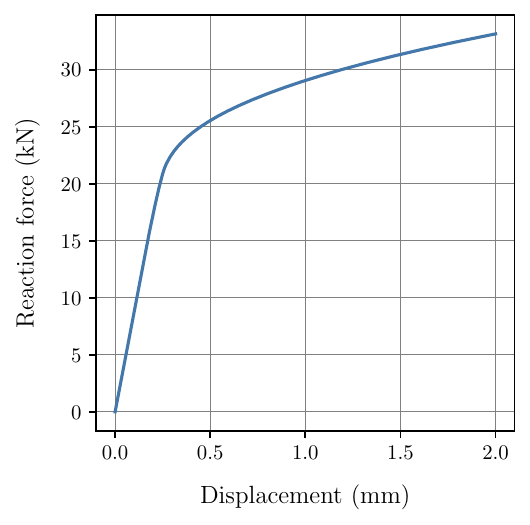}
			\caption{}
			\label{subfig:lou_tensile_double_notched_force_displacement}
		\end{subfigure}\hfill
		\begin{subfigure}[b]{0.45\textwidth}
			\centering
			\includegraphics[width=\textwidth]{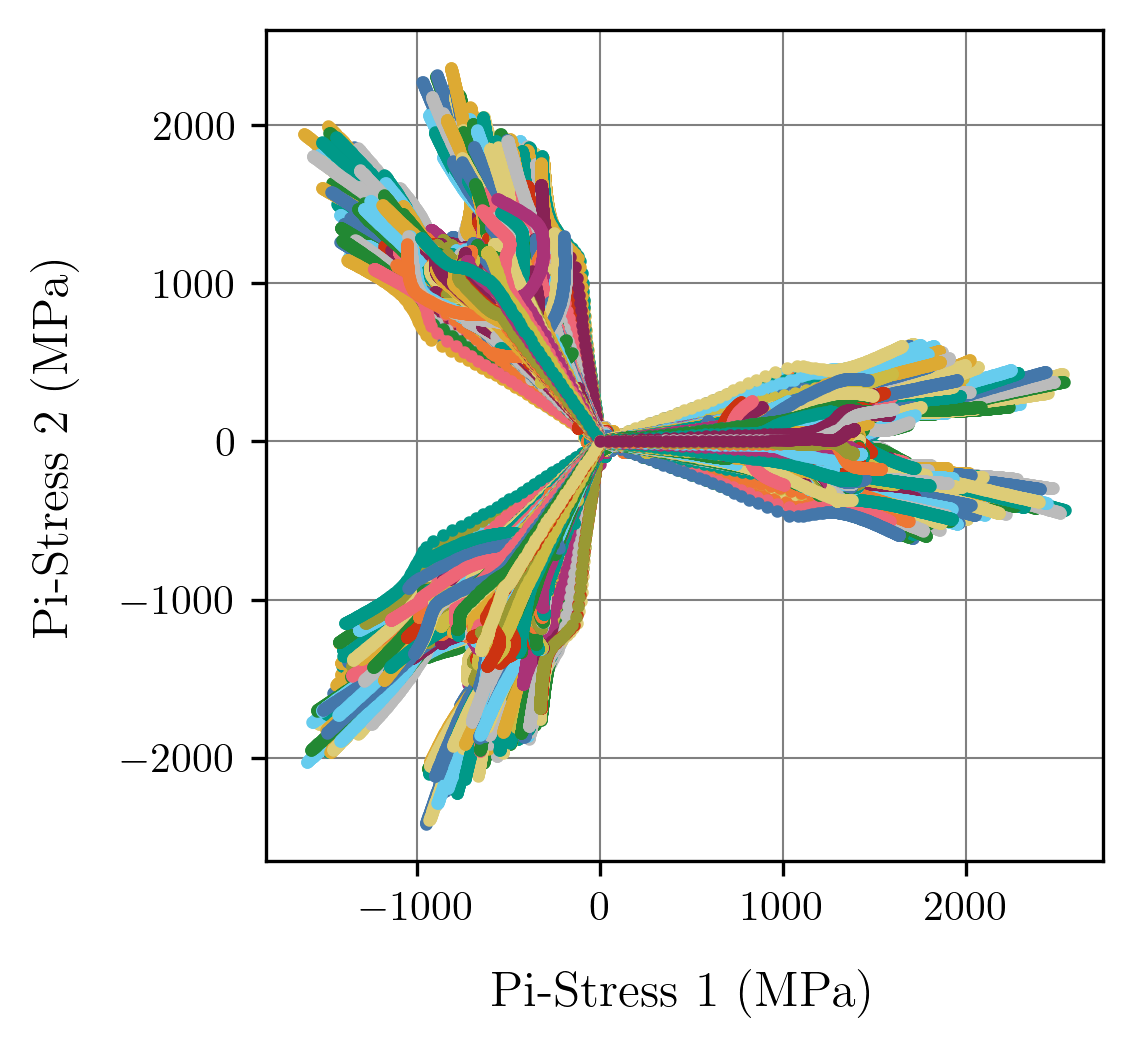}
			\caption{}
			\label{subfig:lou_tensile_double_notched_stress_path_pi_stress_1v2}
		\end{subfigure}
		\begin{subfigure}[b]{0.45\textwidth}
			\centering
			\includegraphics[width=\textwidth]{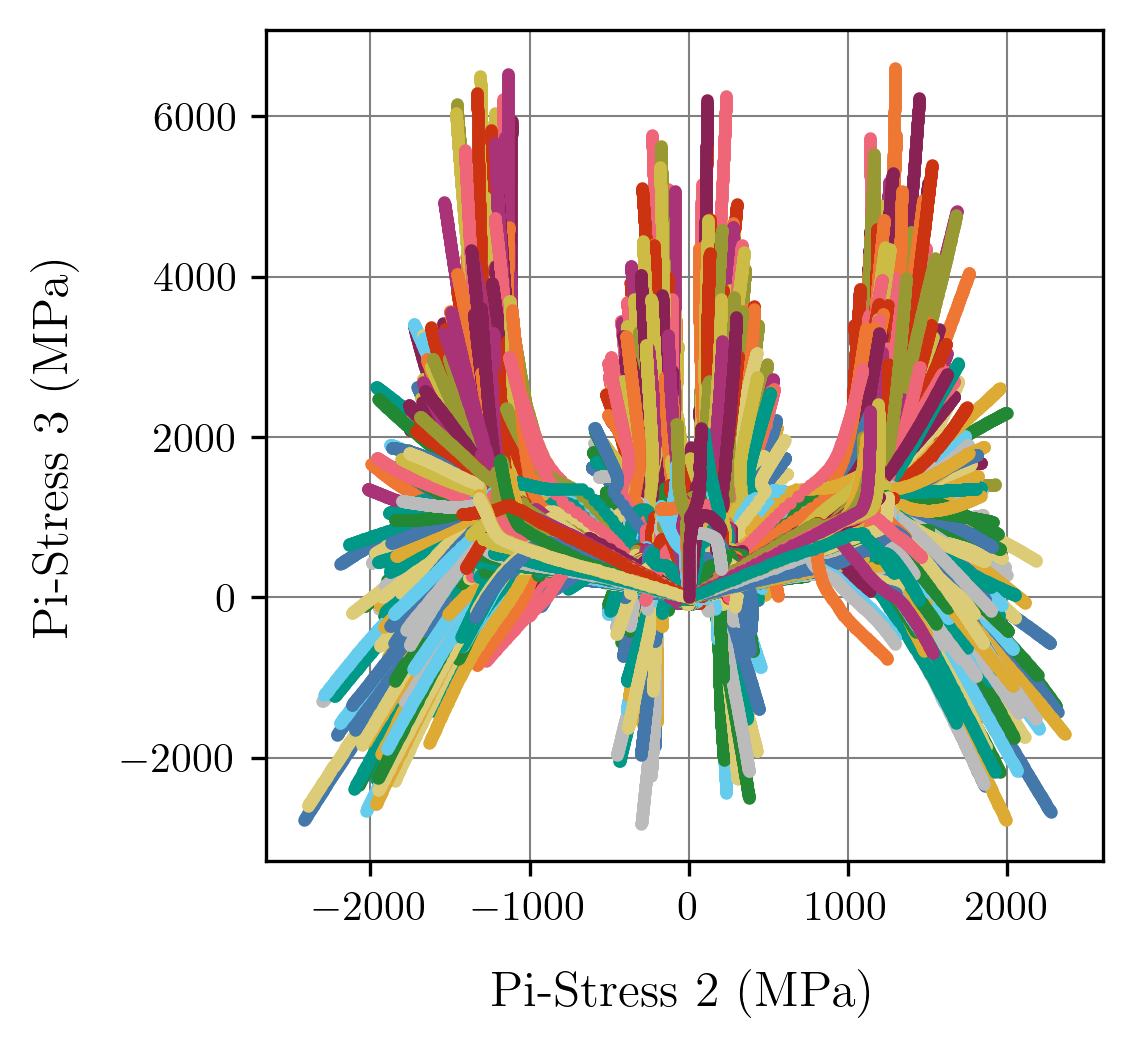}
			\caption{}
			\label{subfig:lou_tensile_double_notched_stress_path_pi_stress_2v3}
		\end{subfigure}\hfill
		\caption{Tensile double notched LZY data ($E=110$GPa, $\nu=0.33$, $a=1.0$, $b=0.05$, $c=1.0$, $d=0.5$, $s_{0}=900\,$MPa, $s_{1}=700\,$MPa, $s_{2}=0.5$): \subref{subfig:lou_tensile_double_notched_acc_p_strain} Accumulated plastic strain field resulting from the total prescribed displacement; \subref{subfig:lou_tensile_double_notched_force_displacement} Uniaxial force-displacement response;
			\subref{subfig:lou_tensile_double_notched_stress_path_pi_stress_1v2} $\pi_{1}-\pi_{2}$ stress projection view of the induced local strain-stress paths;
			\subref{subfig:lou_tensile_double_notched_stress_path_pi_stress_2v3} $\pi_{2}-\pi_{3}$ stress projection view of the induced local strain-stress paths.}
		\label{fig:lou_tensile_double_notched_training}
	\end{figure}
	
\begin{figure}[h!]
		\centering
		\begin{subfigure}[b]{0.45\textwidth}
			\centering
			\includegraphics[width=\textwidth]{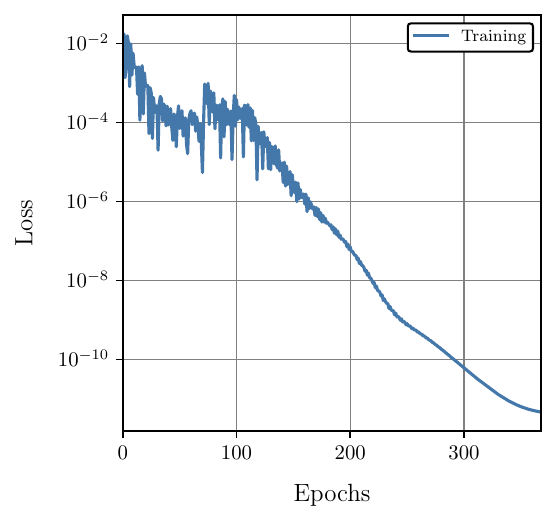}
			\caption{}
			\label{subfig:lou_tensile_double_notched_training_loss_history}
		\end{subfigure}
		\begin{subfigure}[b]{0.445\textwidth}
			\centering
			\includegraphics[width=\textwidth]{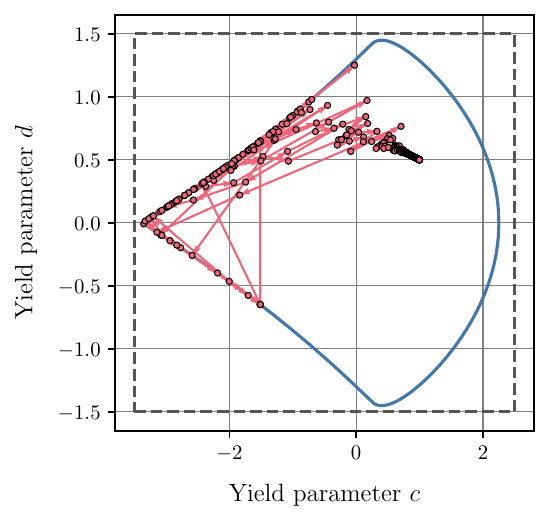}
			\caption{}
			\label{subfig:lou_model_global_convexity_history}
		\end{subfigure}\hfill
		\caption{LZY model global discovery from LZY displacement-force data ($E=110$GPa, $\nu=0.33$, $a=1.0$, $b=0.05$, $c=1.0$, $d=0.5$, $s_{0}=900\,$MPa, $s_{1}=700\,$MPa, $s_{2}=0.5$) with random model initialization:
			\subref{subfig:lou_tensile_double_notched_training_loss_history} Force equilibrium loss history throughout the discovery process;
			\subref{subfig:lou_model_global_convexity_history} Yield surface convexity enforcement throughout the discovery process.}
		\label{fig:lou_tensile_double_notched_loss}
	\end{figure}

\begin{figure}[h!]
		\centering
		\begin{subfigure}[b]{0.38\textwidth}
			\centering
			\includegraphics[width=\textwidth]{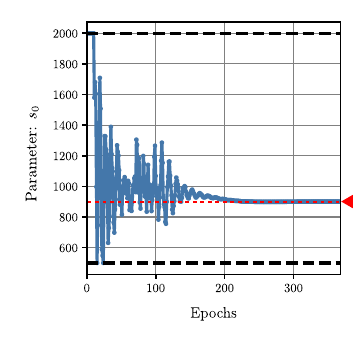}
			\caption{}
			\label{subfig:lou_model_parameter_history_s0}
		\end{subfigure}
		\begin{subfigure}[b]{0.38\textwidth}
			\centering
			\includegraphics[width=\textwidth]{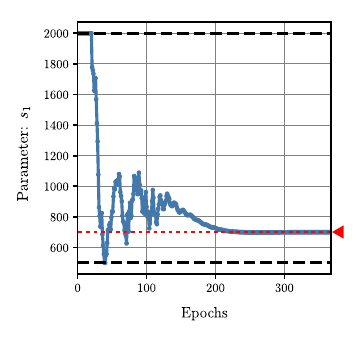}
			\caption{}
			\label{subfig:lou_model_parameter_history_a}
		\end{subfigure}\hfill
		\begin{subfigure}[b]{0.38\textwidth}
			\centering
			\includegraphics[width=\textwidth]{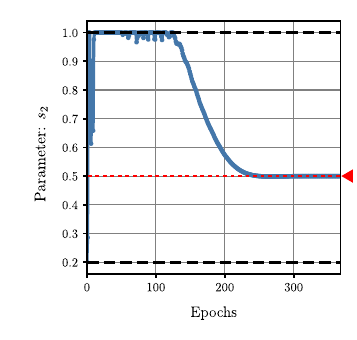}
			\caption{}
			\label{subfig:lou_model_parameter_history_b}
		\end{subfigure}
		\begin{subfigure}[b]{0.38\textwidth}
			\centering
			\includegraphics[width=\textwidth]{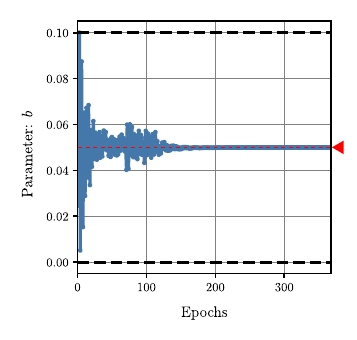}
			\caption{}
			\label{subfig:lou_model_parameter_history_yield_b_s0}
		\end{subfigure}\hfill
		\begin{subfigure}[b]{0.38\textwidth}
			\centering
			\includegraphics[width=\textwidth]{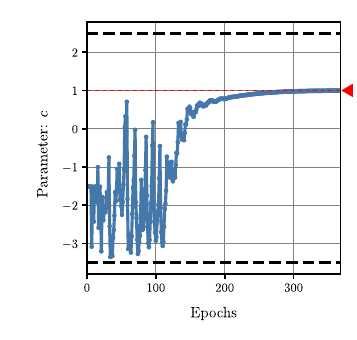}
			\caption{}
			\label{subfig:lou_model_parameter_history_yield_c_s0}
		\end{subfigure}
		\begin{subfigure}[b]{0.38\textwidth}
			\centering
			\includegraphics[width=\textwidth]{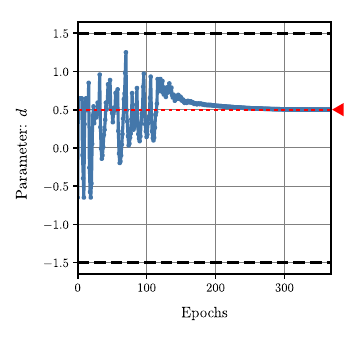}
			\caption{}
			\label{subfig:lou_model_parameter_history_yield_d_s0}
		\end{subfigure}
		\caption{LZY model global discovery from LZY displacement-force data ($E=110$GPa, $\nu=0.33$, $a=1.0$, $b=0.05$, $c=1.0$, $d=0.5$, $s_{0}=900\,$MPa, $s_{1}=700\,$MPa, $s_{2}=0.5$) with random model initialization:
			\subref{subfig:lou_model_parameter_history_s0} Yield parameter $s_{0}$;
			\subref{subfig:lou_model_parameter_history_a} Yield parameter $s_{1}$;
			\subref{subfig:lou_model_parameter_history_b} Yield parameter $s_{2}$;
			\subref{subfig:lou_model_parameter_history_yield_b_s0} Yield parameter $b$;
			\subref{subfig:lou_model_parameter_history_yield_c_s0} Yield parameter $c$; \subref{subfig:lou_model_parameter_history_yield_d_s0} Yield parameter $d$. Black dashed lines correspond to the optimization upper and lower bounds. Red dashed lines correspond to the parameters `ground-truth'.}
		\label{fig:lou_tensile_double_notched_parameters}
	\end{figure}

	Compared to the previous dogbone geometry, this specimen promotes a greater diversity of the local strain-stress paths due to the double notch (see Figures~\ref{subfig:lou_tensile_double_notched_stress_path_pi_stress_1v2} and \ref{subfig:lou_tensile_double_notched_stress_path_pi_stress_2v3}). This diversity is essential, as it enables the discovery of not only the yield pressure dependency and nonlinear hardening law parameters but also the LZY yield surface curvature and strength differential effect (see Figures~\ref{subfig:lou_yield_b_parametric}-\ref{subfig:lou_yield_d_parametric}). In what concerns the global discovery of the LZY parameters, illustrated in Figures~\ref{fig:lou_tensile_double_notched_loss} and \ref{fig:lou_tensile_double_notched_parameters}, the same wide, exploratory ranges used in the local discovery setting are applied to all parameters.
	As shown in Figures~\ref{subfig:lou_model_parameter_history_s0}-\ref{subfig:lou_model_parameter_history_yield_d_s0}, the `ground-truth' values of all parameters are successfully found after ca. 350 epochs, being the yield surface convexity enforced throughout the discovery process (see Figure~\ref{subfig:lou_model_global_convexity_history}).

	\begin{remark}
		The global discovery process for all previous examples is performed three times with random parameter initialization. Despite traversing distinct paths in the loss landscape and requiring a varying number of epochs, all parameters converged to their respective `ground-truth' values.
	\end{remark}
	
	Before wrapping up this section, the effectiveness of ADiMU's global indirect model discovery when the data includes noise is addressed. As previously mentioned, noise stems essentially from the experimental displacement field measurements performed with Digital Image Correlation (DIC) or Digital Volume Correlation (DVC). Although effective denoising methods exist for post-processing experimental data, artificial noise is introduced into the specimen's ground-truth displacement history in what follows. In these preliminary analyses, we assume a Uniform distribution to model the synthetic noise. Based on the in-plane measurement resolution of modern data acquisition setups, the noise amplitude is defined as $u_{\mathrm{amp}}=2.7\times 10^{-4}\,$mm and the corresponding zero-centered Uniform distribution as\footnote{The synthetic noise follows a uniform distribution defined directly in the displacement space.}
	\begin{equation}
		\tilde{u} \sim \mathcal{U}(-0.5 \, u_{\mathrm{amp}}, \, 0.5 \, u_{\mathrm{amp}}) \, \mathrm{(mm)} \, .
		\label{eq:uniform_noise_dic}
	\end{equation}
	The synthetic noise, which is homoscedastic in this context, is then sampled independently for each mesh node, displacement component, and time step. Lastly, the noisy displacement field is generated by superimposing the sampled artificial noise onto the noiseless displacement field.
	
\begin{figure}[h!]
		\centering
		\begin{subfigure}[b]{0.9\textwidth}
			\centering
			\includegraphics[width=\textwidth]{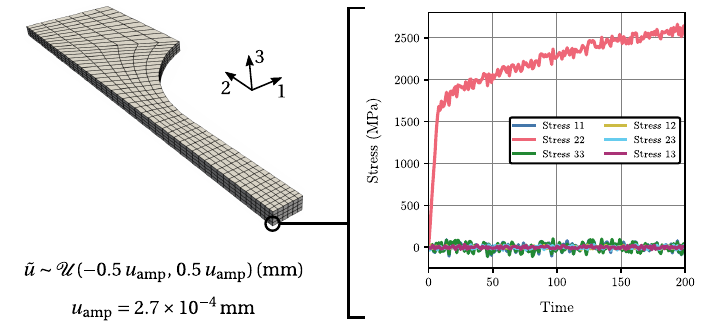}
			\caption{}
			\label{subfig:vm_tensile_dogbone_noisy_displacements}
		\end{subfigure}\vspace*{10pt} \hfill
		\begin{subfigure}[b]{0.9\textwidth}
			\centering
			\includegraphics[width=\textwidth]{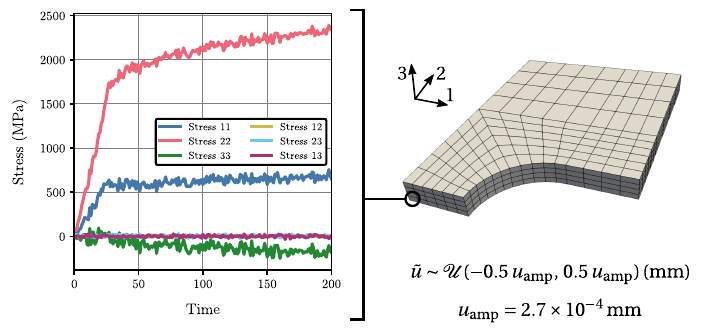}
			\caption{}
			\label{subfig:lou_tensile_double_notched_noisy_displacements}
		\end{subfigure}
		\caption{The effect of superimposing artificial noise onto the displacement field in the local stress field: \subref{subfig:vm_tensile_dogbone_noisy_displacements} Tensile dobgone specimen with VM material model ($E=110\,$GPa, $\nu=0.33$, $s_{0}=900\sqrt{3}\,$MPa, $s_{1}=700\sqrt{3}\,$MPa, $s_{2}=0.5$); \subref{subfig:vm_tensile_dogbone_noisy_displacements} Tensile double notched specimen with LZY model ($E=110$GPa, $\nu=0.33$, $a=1.0$, $b=0.05$, $c=1.0$, $d=0.5$, $s_{0}=900\,$MPa, $s_{1}=700\,$MPa, $s_{2}=0.5$). Synthetic noise, $\tilde{u}$, is modeled by Uniform distribution with amplitude $u_{\mathrm{amp}}$.}
		\label{fig:tensile_noisy_displacements}
	\end{figure}

\begin{figure}[h!]
		\centering
		\begin{subfigure}[b]{0.488\textwidth}
			\centering
			\includegraphics[width=\textwidth]{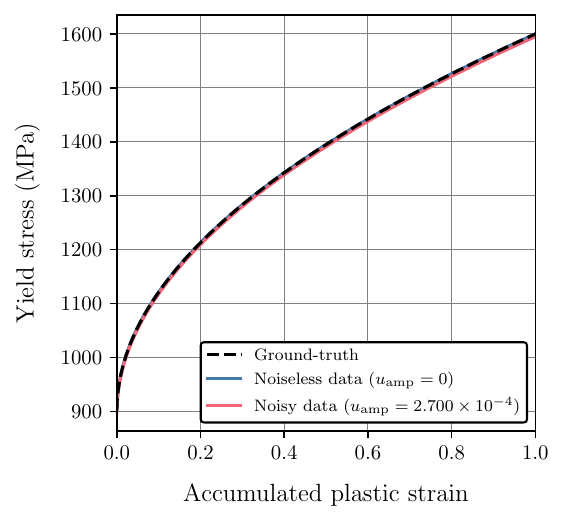}
			\caption{}
			\label{subfig:lou_tensile_double_notched_noisy_hard_laws}
		\end{subfigure} \hfill
		\begin{subfigure}[b]{0.499\textwidth}
			\centering
			\includegraphics[width=\textwidth]{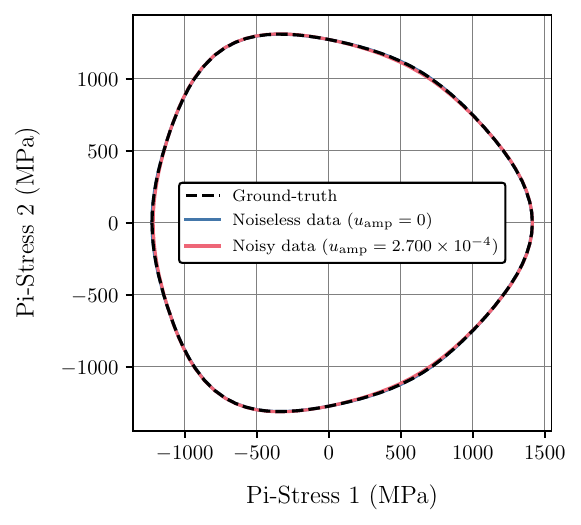}
			\caption{}
			\label{subfig:lou_tensile_double_notched_noisy_yield_pi}
		\end{subfigure}
		\caption{Comparison between the converged parameters found in the LZY model global discovery from LZY noiseless and noisy displacement-force data ($E=110$GPa, $\nu=0.33$, $a=1.0$, $b=0.05$, $c=1.0$, $d=0.5$, $s_{0}=900\,$MPa, $s_{1}=700\,$MPa, $s_{2}=0.5$): \subref{subfig:lou_tensile_double_notched_noisy_hard_laws} Nadai-Ludwik isotropic strain hardening law; \subref{subfig:lou_tensile_double_notched_noisy_yield_pi} Yield surface.}
		\label{fig:tensile_noisy_converged_parameters}
	\end{figure}

	The impact of adding synthetic noise to the displacement field on the local stress field of both dogbone and double notched tensile specimens is shown in Figure~\ref{fig:tensile_noisy_displacements}. The same global model discovery process previously outlined is then applied to the VM tensile dogbone\footnote{Only the hardening law parameters ($s_0$, $s_1$, $s_2$) are considered in the optimization process.} and LZY tensile double notched specimens with noisy data. While the `ground-truth' values of the VM model parameters are successfully found, there are minor discrepancies in the converged values of the LZY model parameters. These are shown in Table~\ref{tab:lzy_noisy_convergence_parameters} for different noise levels, ranging from a noiseless scenario to the assumed noise amplitude ($u_{\mathrm{amp}}=2.7\times 10^{-4}\,$mm). Such discrepancies, however, do not prevent an excellent match of the strain hardening law and LZY yield surface as shown in Figure~\ref{fig:tensile_noisy_converged_parameters}.
	
	\begin{remark}
		For the noisy displacement data defined by $\mathit{u_{\mathrm{amp}}=2.7\times 10^{-4}}\,$mm, the force equilibrium loss achieved with the converged parameters is indeed lower than the one achieved with the `ground-truth' parameters. In addition, the converged parameters are consistently found for multiple realizations with a random parameter initialization. These observations are sensible from a force equilibrium point of view: if the noise leads to a different displacement field for which the corresponding internal forces must balance the same external forces,
		the only solution consists in finding a different set of model parameters that best explain the `new' equilibrium problem. Nevertheless, given that the model does no longer perfectly explain the displacement-force data, the force equilibrium loss is always finite and higher than the one achieved with the `ground-truth' parameters in the noiseless scenario (perfect equilibrium), as expected.
	\end{remark}
	
	\begin{table}[H]
		\caption{LZY model global discovery from LZY noiseless and noisy displacement-force data ($E=110$GPa, $\nu=0.33$, $a=1.0$, $b=0.05$, $c=1.0$, $d=0.5$, $s_{0}=900\,$MPa, $s_{1}=700\,$MPa, $s_{2}=0.5$) with random model initialization. The reference noise amplitude defines the highest noise level, with lower levels derived by dividing it by 2, 4, and 8.}
		\label{tab:lzy_noisy_convergence_parameters}
		\centering
		\setlength{\tabcolsep}{0.5cm}
		\renewcommand{\arraystretch}{1.5}
		\begin{tabular}{ccccccc}
			\toprule
			\multirow{3}{*}{\centering Noise amplitude ($u_{\mathrm{amp}}$)} & \multicolumn{6}{c}{\centering \textbf{Converged parameters}} \\ \cmidrule(l){2-7}
			& \multicolumn{3}{c}{\centering \textbf{Yield surface}} & \multicolumn{3}{c}{\centering \textbf{Isotropic hardening}} \\ \cmidrule(lr){2-4} \cmidrule(l){5-7}
			& $b$ & $c$ & $d$ & $s_{0}$ (MPa) & $s_{1}$ (MPa) & $s_{2}$  \\ \midrule[0.08em]
			$0$ &  0.05 & 1.00 & 0.50 & 900 & 700 & 0.500 \\ \midrule
			$3.375 \times 10^{-5}$ & 0.05 & 1.00 & 0.50 & 899 & 700 & 0.499 \\ \midrule
			$6.750 \times 10^{-5}$ & 0.05 & 1.00 & 0.50 & 900 & 699 & 0.497 \\ \midrule
			$1.350 \times 10^{-4}$ & 0.05 & 0.98 & 0.51 & 897 & 700 & 0.496 \\ \midrule
			$2.700 \times 10^{-4}$ & 0.05 & 0.87 & 0.53 & 896 & 699 & 0.498 \\
			\bottomrule
		\end{tabular}
	\end{table}
	
	\begin{remark}
		Incorporating elastic parameters into the global discovery process significantly impairs performance when handling noisy displacement data. In certain cases, it can even prevent the convergence of all parameters entirely, similar to high noise levels. This is hypothesized to stem from the high-frequency sequence of loading/unloading states induced by noisy displacements, where the unloading is solely dependent on the elastic parameters. Such observation further supports the exclusion of the elastic constants from the optimization problem whenever possible.
	\end{remark}

	\subsection{Neural network and hybrid models}
	
	This last section is focused on the global discovery of neural network and hybrid models from displacement-force data. Unlike the previous examples, where the discovery process heavily relies on a physics-based conventional model, the discovery of a GRU material model is entirely data-driven. Unsurprisingly, the local strain-stress paths generated in the uniaxial tensile tests of both dogbone and double-notched specimens (see Figures~\ref{subfig:dp_tensile_dogbone_stress_path_pi_stress_1v2} and \ref{subfig:lou_tensile_double_notched_stress_path_pi_stress_1v2}) lack the richness needed for an accurate, data-driven learning of the underlying model behavior.
	
	\begin{figure}[h!]
		\centering
		\begin{subfigure}[b]{0.45\textwidth}
			\centering
			\includegraphics[width=\textwidth]{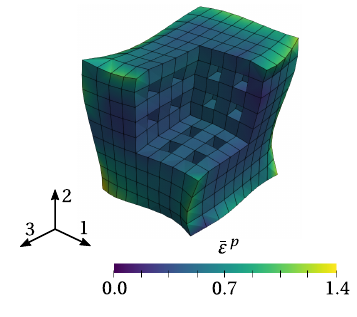}
			\caption{}
			\label{subfig:vm_random_material_patch_acc_p_strain}
		\end{subfigure}
		\begin{subfigure}[b]{0.45\textwidth}
			\centering
			\includegraphics[width=\textwidth]{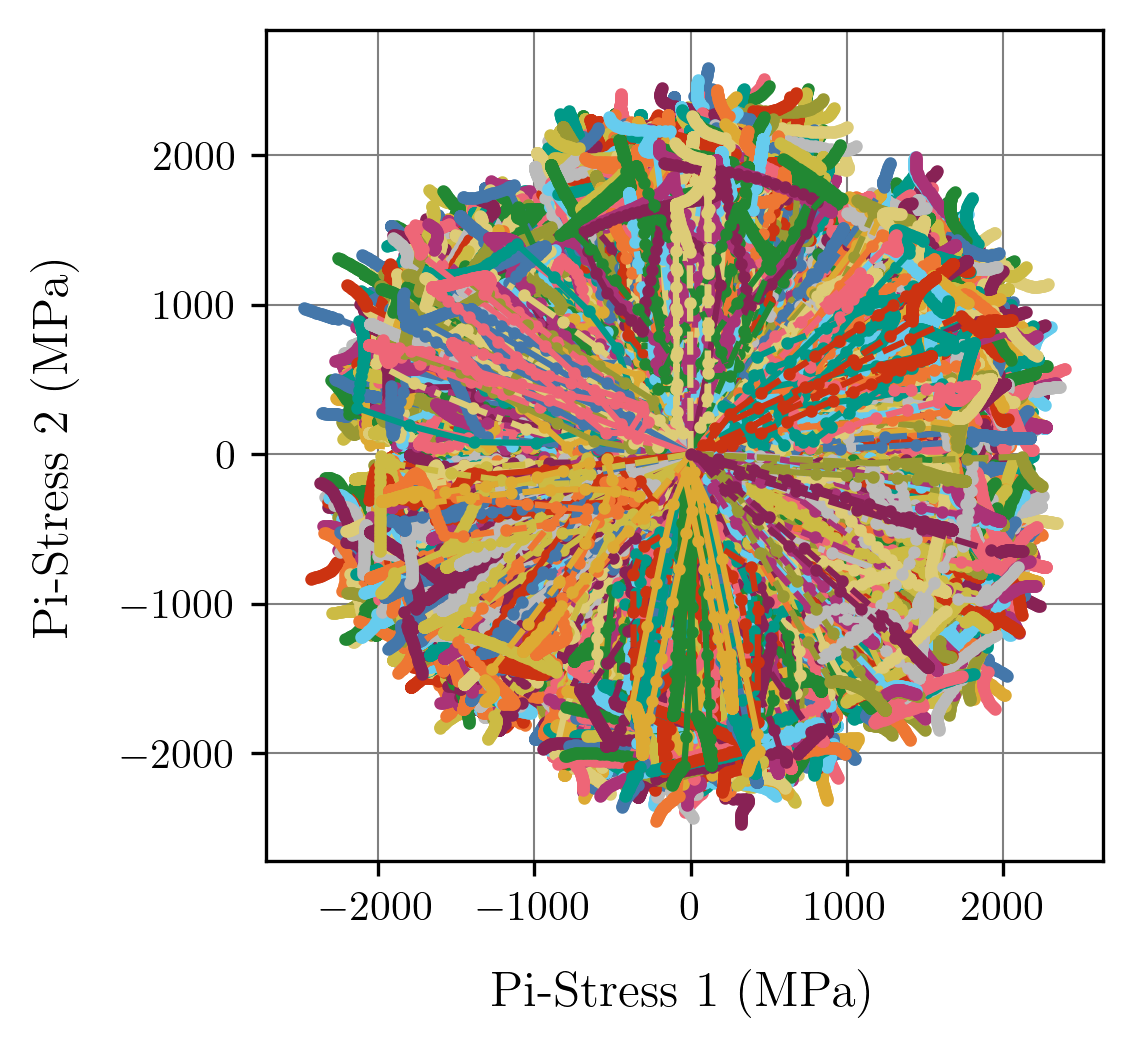}
			\caption{}
			\label{subfig:vm_random_stress_path_pi_stress_1v2}
		\end{subfigure}
		\caption{Randomly deformed material patch specimen VM data ($E=110\,$GPa, $\nu=0.33$, $s_{0}=900\sqrt{3}\,$MPa, $s_{1}=700\sqrt{3}\,$MPa, $s_{2}=0.5$): \subref{subfig:vm_random_material_patch_acc_p_strain} Accumulated plastic strain field resulting from the prescribed displacement ($\lambda(t) = -0.92$);
			\subref{subfig:vm_random_stress_path_pi_stress_1v2} $\pi_{1}-\pi_{2}$ stress projection view of the induced local strain-stress paths.}
		\label{fig:vm_random_material_patch_training}
	\end{figure}
	
	For this reason, we demonstrate this challenging scenario by resorting to a randomly deformed material patch specimen generated with the Stochastic Patch Deformation Generator (SPDG) that we created (see ~\ref{sec:spdg}). This idealized specimen, whose details can be found in ~\ref{ssec:random_patch}, promotes a broad coverage of the strain-stress space (see Figure~\ref{fig:vm_random_material_patch_training}) by coupling (i) prescribed random boundary displacements, (ii) an internal regular grid of cubic voids, and (iii) a randomly polynomial non-monotonic proportional loading scheme. Importantly, the specimen shown in Figure~\ref{fig:vm_random_material_patch_training}) undergoes only one loading-unloading cycle (just like the uniaxial test), but where the boundary monotonically increases according to a displacement obeying a polynomial surface. Without any loss of generality, we assume the `ground-truth' von Mises (VM) model to compute the displacement-force data (see ~\ref{sec:material_parameters}).
	
	\begin{remark}
		The deformation bounds of the randomly deformed material patch specimen are set to ensure that the resulting local strain-stress paths span stress magnitude ranges comparable to those of the random polynomial strain-stress paths. As a consequence, the following analyses assess the performance of the discovered models when extrapolating in terms of strain-stress path diversity.
	\end{remark}

	\begin{table}[h!]
		\caption{Comparison of the average prediction Normalized Root Mean Squared Error (NRMSE) of different models discovered from VM data ($E=110\,$GPa, $\nu=0.33$, $s_{0}=900\sqrt{3}\,$MPa, $s_{1}=700\sqrt{3}\,$MPa, $s_{2}=0.5$) and tested in a (unseen) local synthetic data set of 512 random polynomial strain-stress paths. The local discovery process involves strain-stress data, while the global discovery process only accounts for displacement-force data. The local `random polynomial' data set denotes 2560 random polynomial strain-stress paths, while the local `material patch' data set denotes 4256 strain-stress local paths (80\%/20\% training/validation split) of the randomly deformed material patch specimen that arise from a single boundary loading condition (see Figure \ref{fig:vm_random_material_patch_training}). The hybrid model has a candidate-corrector type architecture, where the candidate is a VM model with wrong parameters ($E=90\,$GPa, $\nu=0.29$, $\sigma_{y,\,0}=1200\,$MPa, $H=600\,$MPa) and the corrector is a GRU model.}
		\label{tab:global_material_patch_nrmse}
		\centering
		\setlength{\tabcolsep}{0.30cm}
		\renewcommand{\arraystretch}{1.5}
		\begin{tabular}{ccccccccc}
			\toprule
			\multirow{2}{*}{\parbox[c]{1.5cm}{\centering \bfseries \shortstack{Discovery \\[0.3em] process}}} & \multirow{2}{*}{\parbox[c]{2cm}{\centering \bfseries \shortstack{Training \\[0.3em] data set}}} & \multirow{2}{*}{\parbox[c]{2cm}{\centering \bfseries Model}} & \multicolumn{6}{c}{\centering \bfseries NRMSE (\%)}  \\ \cmidrule(l){4-9}
			& & & $\sigma_{11}$ & $\sigma_{22}$ & $\sigma_{33}$ & $\sigma_{12}$ & $\sigma_{23}$ & $\sigma_{13}$ \\ \midrule[0.08em]
			\multirow{4}{*}{Local} & \multirow{2}{*}{\parbox{2cm}{\centering Random polynomial}} & Hybrid & 0.23 & 0.22 & 0.23 & 0.82 & 0.82 & 0.79 \\ \cmidrule(l){3-9}
			& & GRU & \cellcolor{lightgray} 0.63 & \cellcolor{lightgray} 0.61 & \cellcolor{lightgray} 0.64 & \cellcolor{lightgray}1.51 & \cellcolor{lightgray} 1.58 & \cellcolor{lightgray} 1.47 \\ \cmidrule(l){2-9}
			& \multirow{2}{*}{\parbox{2cm}{\centering Material patch}} & Hybrid & 1.29 & 1.19 & 1.28 & 5.43 & 5.43 & 5.33 \\ \cmidrule(l){3-9}
			& & GRU & \cellcolor{lightgray} 4.95 & \cellcolor{lightgray} 5.03 & \cellcolor{lightgray} 5.26 & \cellcolor{lightgray} 16.45 & \cellcolor{lightgray} 15.97 & \cellcolor{lightgray} 16.26 \\ \midrule
			\multirow{2}{*}{Global} & \multirow{2}{*}{\parbox{2cm}{\centering Material patch}} & Hybrid & 3.82 & 3.64 & 3.94 & 15.96 & 16.27 & 16.44  \\ \cmidrule(l){3-9}
			& & GRU & \cellcolor{lightgray} 8.17 & \cellcolor{lightgray} 7.80 & \cellcolor{lightgray} 8.15 & \cellcolor{lightgray} 30.79 & \cellcolor{lightgray} 32.21 & \cellcolor{lightgray} 30.21 \\ \bottomrule
		\end{tabular}
	\end{table}
	
	\begin{figure}[H]
		\centering
		\includegraphics[width=0.98\textwidth]{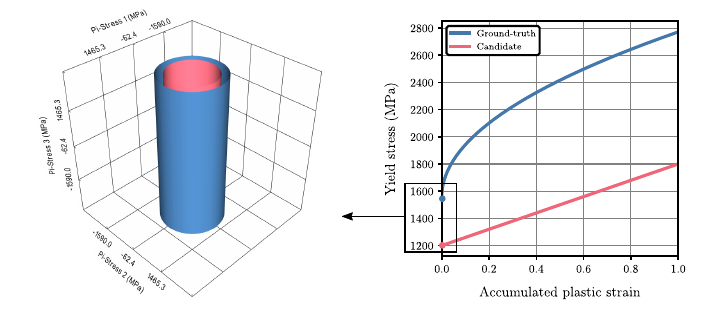}
		\caption{Comparison between the `ground-truth' VM model ($E=110\,$GPa, $\nu=0.33$, $s_{0}=900\sqrt{3}\,$MPa, $s_{1}=700\sqrt{3}\,$MPa, $s_{2}=0.5$) and the candidate VM model ($E=90\,$GPa, $\nu=0.29$, $\sigma_{y,\,0}=1200\,$MPa, $H=600\,$MPa). The candidate VM model has a linear hardening law with initial stress $\sigma_{y,\,0}$ and hardening slope $H$.}
		\label{fig:hybrid_model_vm_candidate}
	\end{figure}
	
\begin{figure}[h!]
		\centering
		\begin{subfigure}[b]{0.49\textwidth}
			\centering
			\includegraphics[width=\textwidth]{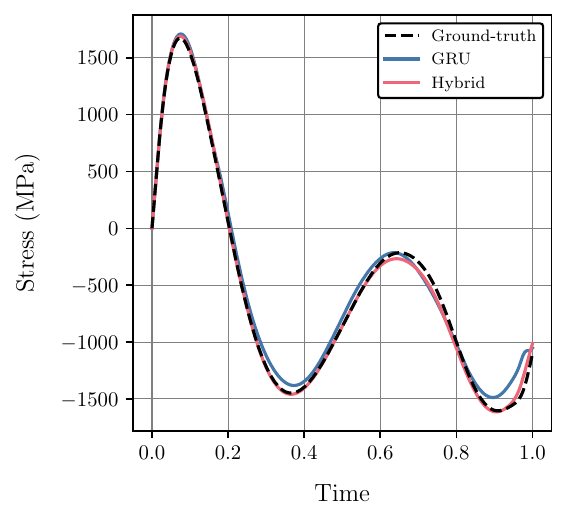}
			\caption{}
			\label{subfig:global_gru_hybrid_stress_paths_stress_11_sample}
		\end{subfigure}
		\begin{subfigure}[b]{0.49\textwidth}
			\centering
			\includegraphics[width=\textwidth]{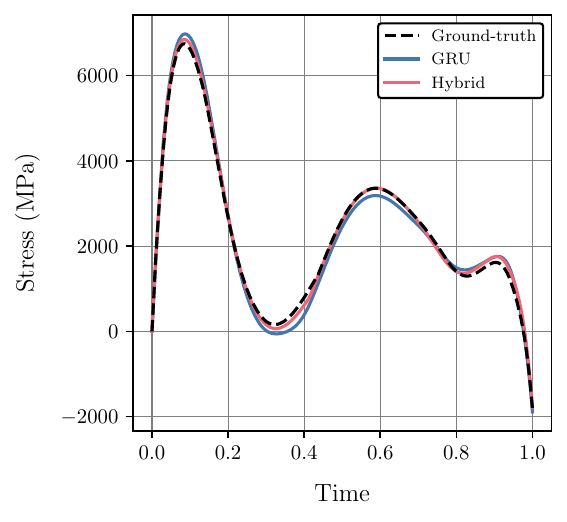}
			\caption{}
			\label{subfig:global_gru_hybrid_stress_paths_stress_22_sample}
		\end{subfigure}\hfill
		\begin{subfigure}[b]{0.48\textwidth}
			\centering
			\includegraphics[width=\textwidth]{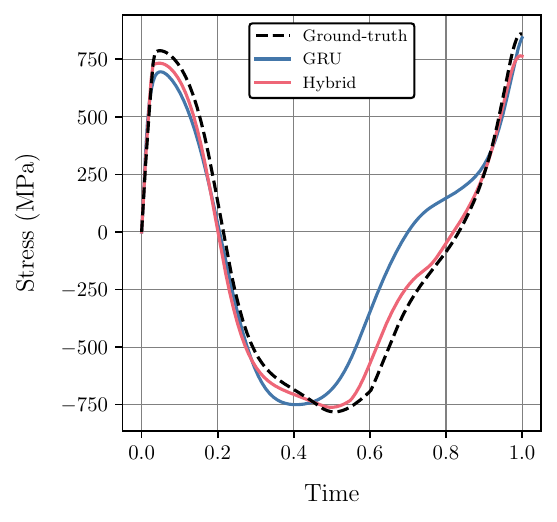}
			\caption{}
			\label{subfig:global_gru_hybrid_stress_paths_stress_12_sample}
		\end{subfigure}
		\begin{subfigure}[b]{0.49\textwidth}
			\centering
			\includegraphics[width=\textwidth]{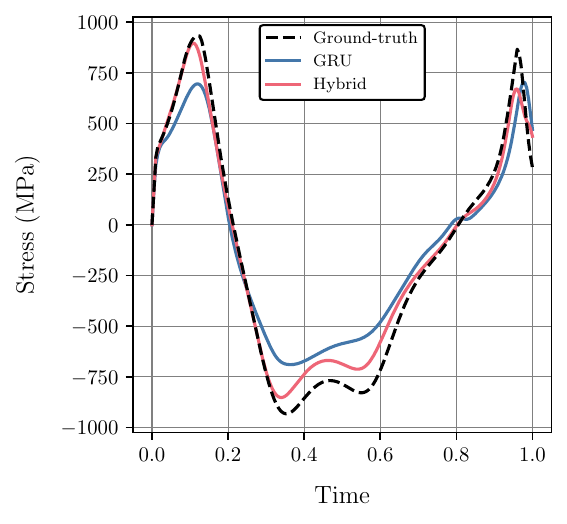}
			\caption{}
			\label{subfig:global_gru_hybrid_stress_paths_stress_23_sample}
		\end{subfigure}\hfill
		\caption{Performance of the GRU and hybrid material model discovered from VM displacement-force data ($E=110\,$GPa, $\nu=0.33$, $s_{0}=900\sqrt{3}\,$MPa, $s_{1}=700\sqrt{3}\,$MPa, $s_{2}=0.5$), with random model initialization, on (unseen) testing random polynomial strain-stress paths and different stress components: \subref{subfig:global_gru_hybrid_stress_paths_stress_11_sample} Normal stress 11; \subref{subfig:global_gru_hybrid_stress_paths_stress_22_sample} Normal stress 22;
			\subref{subfig:global_gru_hybrid_stress_paths_stress_12_sample} Shear stress 12;
			\subref{subfig:global_gru_hybrid_stress_paths_stress_23_sample} Shear stress 23.}
		\label{fig:global_gru_hybrid_stress_paths}
	\end{figure}

	In what follows, we compare the performance of two different models discovered from displacement-force data. The first model is a multi-layer gated recurrent unit (GRU) material model, akin to the one discussed in Section~\ref{ssec:local_ml_models}. The second is a hybrid model featuring the simple candidate-corrector architecture illustrated in Figure~\ref{fig:hybrid_architecture_cc}. The candidate is a VM model with erroneous elastic parameters ($E=90\,$GPa, $\nu=0.29$) and a highly inaccurate linear hardening law (see Figure~\ref{fig:hybrid_model_vm_candidate}). This choice aims to highlight the potential benefits of incorporating a physics-based model in the global discovery setting, even when such model is innacurate. The performance of both models is evaluated on an unseen testing data set comprising 512 random polynomial strain-stress paths, with the average prediction NRMSE for each stress component summarized in Table~\ref{tab:global_material_patch_nrmse}. Examples of testing samples matching the average NRMSE are also provided in Figure~\ref{fig:global_gru_hybrid_stress_paths} for different components.
	
\begin{figure}[hbt]
		\centering
		\begin{subfigure}[b]{0.453\textwidth}
			\centering
			\includegraphics[width=\textwidth]{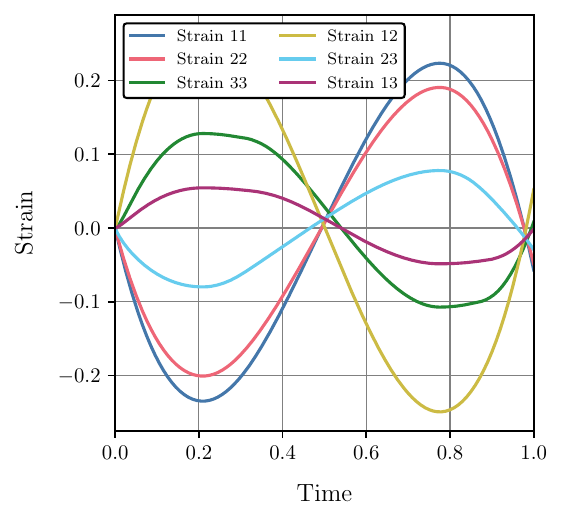}
			\caption{}
			\label{subfig:vm_random_patch_path}
		\end{subfigure}
		\begin{subfigure}[b]{0.47\textwidth}
			\centering
			\includegraphics[width=\textwidth]{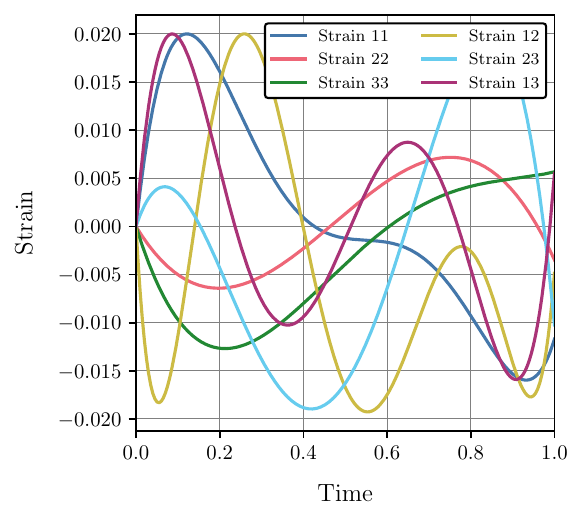}
			\caption{}
			\label{subfig:vm_random_polynomial_path}
		\end{subfigure}
		\caption{Comparison between \subref{subfig:vm_random_patch_path} an induced local strain path in the randomly deformed material patch specimen and \subref{subfig:vm_random_polynomial_path} a random polynomial strain path.}
		\label{fig:vm_comparison_random_specimen_path}
	\end{figure}

	It is remarkable that a fully data-driven GRU material model, discovered solely from displacement-force data, achieves prediction errors of roughly 8\% for normal and 30\% for shear stress components in random polynomial strain-stress paths. Besides the absence of a physics-based foundation, the results transpire two main reasons to explain these accuracy shortcomings. First, discovering a material model indirectly from displacement-force data is inherently complex. Note that the performance improves significantly when the GRU material model is discovered directly (around 5\% and 16\% error for normal and shear stress components, respectively) from the specimen local strain-stress paths (see Table~\ref{tab:global_material_patch_nrmse}). Second, despite the random nature of the specimen deformation, the diversity is still limited when compared to the random polynomial strain-stress paths. On the one hand, the specimen strain-stress paths are approximately synchronized with the total load factor, while the different components are completely asynchronous in a random polynomial strain-stress path (see Figure~\ref{fig:vm_comparison_random_specimen_path}). On the other hand, it is challenging to promote substantial multi-axial shear deformations in the specimen. Note that the accuracy of the GRU material model improves substantially when discovered from a data set of random polynomial strain-stress paths, achieving errors of approximately 0.6\% for normal and 1.5\% for shear stress components. While the same observations hold for the hybrid material model, the physics-based candidate model significantly enhances its prediction accuracy compared to the GRU material model. In the global discovery setting, the prediction error drops to approximately 4\% for normal and 16\% for shear stress components.

	Last but not least, a fundamental aspect concerning the diversity of the specimen local strain-stress paths is in order. As seen throughout Section~\ref{sec:local_model_discovery}, the performance of both neural network and hybrid models improves with the training data set size. Provided the model architecture is sufficiently flexible, such accuracy improvements are expected as long as additional data enhances the coverage of the strain-stress space. This is essentially the case with random polynomial strain-stress paths, where each new path is independently generated. However, the local strain-stress paths induced by a given specimen are strongly influenced by its geometry and loading conditions, not to mention that the number of paths is closely tied to the specimen spatial discretization. Consequently, it comes as no surprise that these paths often exhibit redundancy, with similar strain-stress coverage from spatially proximate points and/or those undergoing comparable deformation paths.
	
	\begin{figure}[h!]
		\centering
		\begin{subfigure}[b]{0.482\textwidth}
			\centering
			\includegraphics[width=\textwidth]{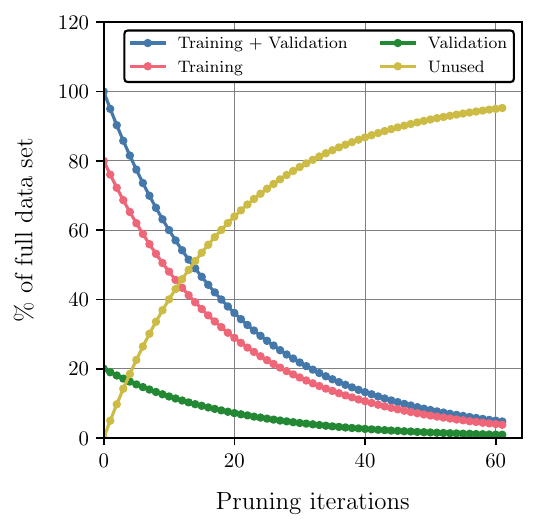}
			\caption{}
			\label{subfig:pruning_iterations_dataset_sizes}
		\end{subfigure} 
		\begin{subfigure}[b]{0.49\textwidth}
			\centering
			\includegraphics[width=\textwidth]{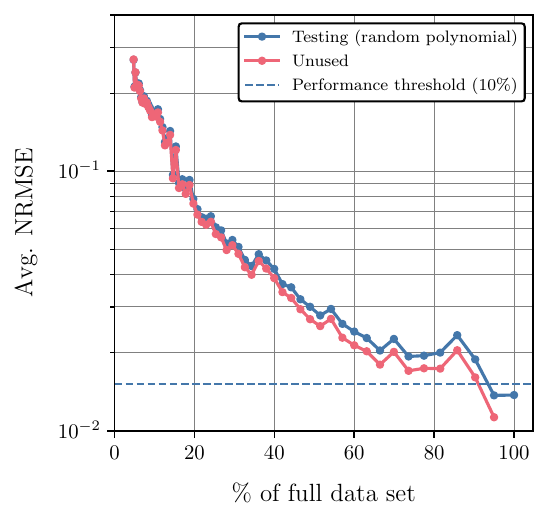}
			\caption{}
			\label{subfig:random_polynomial_pruning_testing_convergence}
		\end{subfigure} \vspace*{5pt} \hfill \\
		\begin{subfigure}[b]{0.49\textwidth}
			\centering
			\includegraphics[width=\textwidth]{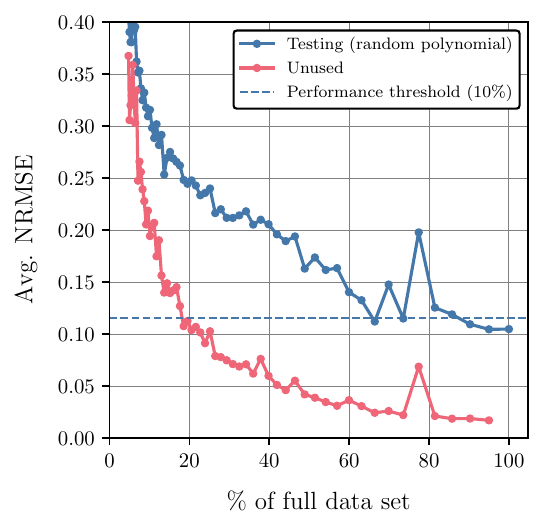}
			\caption{}
			\label{subfig:random_patch_pruning_testing_convergence}
		\end{subfigure}
		\begin{subfigure}[b]{0.49\textwidth}
			\centering
			\includegraphics[width=\textwidth]{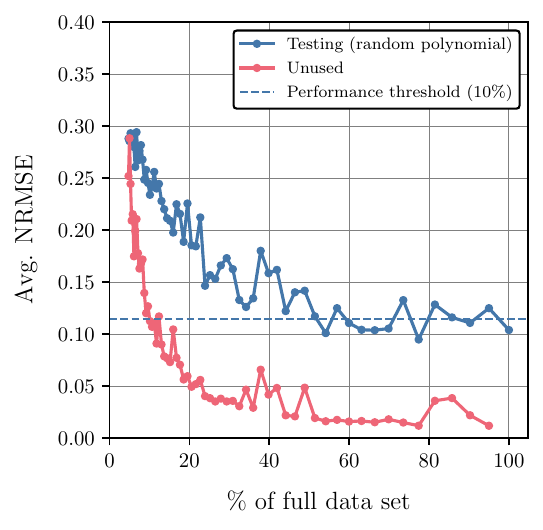}
			\caption{}
			\label{subfig:random_patch_0d25_pruning_testing_convergence}
		\end{subfigure}\hfill
		\caption{Average prediction Normalized Root Mean Squared Error (NRMSE) (averaged over all stress components) of the GRU material model discovered from VM strain-stress data ($E=110\,$GPa, $\nu=0.33$, $s_{0}=900\sqrt{3}\,$MPa, $s_{1}=700\sqrt{3}\,$MPa, $s_{2}=0.5$) and tested in a (unseen) local synthetic data set of 512 random polynomial strain-stress paths. The training data stems from pruning different data sets with the same pruning algorithm scheduler: \subref{subfig:pruning_iterations_dataset_sizes} Data set relative sizes throughout pruning history;
			\subref{subfig:random_polynomial_pruning_testing_convergence} 2560 random polynomial strain-stress paths; \subref{subfig:random_patch_pruning_testing_convergence} 5320 local strain-stress paths of randomly deformed material patch specimen: \subref{subfig:random_patch_0d25_pruning_testing_convergence} 5320 local strain-stress paths of randomly deformed material patch specimen with 75\% less deformation magnitude ($\lambda^{*}(t) = 0.25 \lambda(t)$).}
		\label{fig:pruning_datasets_gru}
	\end{figure}
	
	These arguments are quantified in Figure~\ref{fig:pruning_datasets_gru}, where a data set pruning algorithm (see ~\ref{sec:pruning_algorithm}) based on the GRU material model is explored to evaluate the redundancy of local strain-stress paths data sets. Redundancy is evaluated below using a performance threshold defined as a 10\% testing accuracy degradation (based on the average NRMSE over all stress components) compared to the model trained on the full data set. Data that can be removed without surpassing this accuracy degradation threshold is deemed redundant. The same pruning algorithm scheduler, shown in Figure~\ref{subfig:pruning_iterations_dataset_sizes}, is applied to three different data sets, namely: (i) a data set of 2560 random polynomial strain-stress paths (see Figure~\ref{subfig:random_polynomial_pruning_testing_convergence}); (ii) the  randomly deformed material patch specimen local data set of 5320 strain-stress paths (see Figure~\ref{subfig:random_patch_pruning_testing_convergence}); and (iii) the local data set of the randomly deformed material patch specimen loaded with 75\% less deformation magnitude (see Figure~\ref{subfig:random_patch_0d25_pruning_testing_convergence}). To account for randomness in model initialization and data set splitting, three pruning realizations are conducted for each data set. The pruning results yield a redundancy range of [0-10]\%, [15-40]\% and [50-60]\% for the aforementioned data sets, respectively. These expected results are further supported by the high, stable testing accuracy on the unused data set throughout the redundancy range.

	\begin{remark}
		Assuming a 40\% redundancy upper bound for the local data set of the randomly deformed material patch specimen, and applying an 80\%/20\% training/validation split, the 5320 local strain-stress paths (extracted from the integration points) reduce to approximately 2553 effective (non-redundant) paths. This aligns with the 2560 random polynomial strain-stress paths in the comparison shown in Table~\ref{tab:global_material_patch_nrmse}, further emphasizing the limited diversity of the specimen data set compared to the random polynomial data set and the resulting performance degradation.
	\end{remark}

	Despite the idealized nature of the specimen explored here, the previous findings underscore the critical role of specimen design in the global discovery of neural network or hybrid material models. Beyond experimental feasibility in manufacturing and testing, maximizing the diversity of the induced local strain-stress paths is essential for discovering the underlying material behavior and achieving accurate, reliable predictions. This challenge is outside the scope of this contribution and is addressed in the future remarks.
	
	\section{Conclusions and future remarks \label{sec:conclusions}}
	This paper introduces Automatically Differentiable Model Updating (ADiMU), the first fully automatically differentiable framework that finds any history-dependent material model from full-field displacement and global force data (global, indirect discovery) or from strain-stress data (local, direct discovery). Remarkably, ADiMU requires no fine-tuning of hyperparameters or additional quantities beyond those ineherent to the user-selected material model architecture and optimizer. We demonstrate ADiMU's performance in numerous examples involving history-dependency, encompassing conventional (physics-based), neural network (data-driven), and hybrid material models. Furthermore, ADiMU is released as an open-source computational tool, integrated into a carefully designed and documented software named HookeAI. This facilitates the integration, evaluation and application of future material model architectures by openly supporting the research community.
	
	Akin to conventional material models in finite element simulation software, ADiMU allows a straightforward selection or integration of any kind of parameterized material model, ranging from a few to millions of parameters. Beyond its practical implications, this establishes a significant milestone in the literature: the novel material model architectures that have been extensively developed by the community can now be seamlessly tested and benchmarked in both local and global model discovery contexts. The numerous examples presented in this paper illustrate that the diversity, not merely the quantity, of directly or indirectly available strain-stress data is fundamental for discovering accurate material models. This is particularly crucial as material models become more expressive, particularly with neural networks, where the lack of physics-based knowledge demands a richer and more diverse data set to effectively learn the material behavior. Developing hybrid material model architectures is therefore essential, leveraging the strengths of two approaches: the thermodynamically consistent foundations of conventional material modeling and the remarkable capacity of neural network models to capture a wide range of material behaviors.
	
	Among the several challenges and opportunities laying ahead concerning the application of ADiMU, some are discussed in what follows. First and foremost, despite some preliminary analyses with synthetic noise, it is of the utmost importance to test ADiMU's performance with actual experimental data. Moreover, developing a manufacturable and testable specimen that maximizes the diversity of induced local strain-stress paths \citep{ihuaneyi:2024a, tung:2024a} is essential to discover neural network and hybrid material models. Similarly, developing novel experimental techniques to measure additional data, such as local elastic strains, could significantly enhance the model discovery process \citep{wang:2024a}. Second, every material model discovered using ADiMU is fully differentiable, making it well-suited for integration into advanced topology optimization methods \citep{jia:2025a, vijayakumaran:2025a}. Such integration embodies a completely automated pipeline ranging from experimental material data to the multi-objective, multi-material design of structures with topology optimization. Third, ADiMU's hybrid material model architecture provides a versatile foundation for building a wide variety of models. For instance, it can easily combine high-fidelity conventional models from multiphasic materials with neural network models that learn unknown thermomechanical dependencies and/or interphasic interaction mechanisms. Exploring cooperative data-driven modeling \citep{dekhovich:2023a} to develop models that continuously adapt to new material behaviors is also a promising avenue. Finally, the need to quantify the uncertainty of neural network material models, particularly in data-scarce and/or highly noisy data scenarios, underscores the importance of adopting a Bayesian approach to the model discovery process \citep{yi2025cooperative,coscia:2025a}.

	\section{Acknowledgments}
	The authors acknowledge the support from U.S. Defense Advanced Research Projects Agency (DARPA) Award HR0011-24-2-0333. The authors also acknowledge the expertise shared by Antonios Kontsos concerning the noise arising from DIC experimental measurements. Bernardo P. Ferreira is deeply grateful to António Carneiro and Gawel Kus for several insightful discussions during the development of this paper.
	
	\clearpage
	
	\appendix
	
	\section{HookeAI: An open-source ADiMU framework \label{sec:hookeai}}
	HookeAI is an open-source Python software that bridges computational mechanics and artificial intelligence to discover material models from available data. Its main purpose aligns with Automatically Differentiable Model Updating (ADiMU), a framework that finds any history-dependent material model from full-field displacement and global force data (global, indirect discovery) or from strain-stress data (local, direct discovery). Carefully designed and extensively documented, this software also aims to facilitate the integration, evaluation and application of novel material model architectures by openly supporting the research community.
	
	\begin{figure}[hbt]
		\centering
		\includegraphics[width=0.5\textwidth]{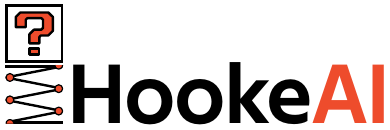}
		\caption{HookeAI - An open-source ADiMU framework.}
		\label{fig:adimu_logo}
	\end{figure}
	
	HookeAI was fully designed and implemented by Bernardo P. Ferreira (Copyright \textcopyright  ~2025).

	\section{Fully implicit Lou-Zhang-Yoon model \label{sec:lzy_model}}
	This section includes the fully implicit Lou-Zhang-Yoon (LZY) elasto-plastic constitutive model proposed in this paper. The model is developed based on the highly versatile yield function proposed by Lou and coworkers \citep{lou:2022a} defined as
	\begin{equation}
		\Phi (\bm{\sigma} \, ; \, a, \, b, \, c, \, d, \, \sigma_{y}) = a \left\{ b I_{1} + \left[ (J_{2}^{3} - c J_{3}^{2})^{\frac{1}{2}} - d J_{3} \right]^{\frac{1}{3}} \right\} - \sigma_{y} \, ,
		\label{eq:lzy_yield_surface}
	\end{equation}
	where $\bm{\sigma}$ is the Cauchy stress tensor, $I_1$ is the first stress invariant, $(J_{2}, \, J_{3})$ are the second and third deviatoric stress invariants, $(a, \, b, \, c, \, d)$ are yield parameters, and $\sigma_{y}$ is the uniaxial yield stress. For instance, it can be easily demonstrated that Equation~\eqref{eq:lzy_yield_surface} recovers well-established smooth yield surfaces for particular sets of parameters $(a, \, b, \, c, \, d)$, such as von Mises \citep{mises:1913}, Drucker-Prager \citep{drucker:1952a} and Cazacu-Barlat \citep{barlat:2004a}. This yield surface has been the object of extensive development and experimental validation by Lou and coworkers \citep{lou:2023a}, and has been employed to predict the yielding of metallic materials under multi-axial loading conditions and accounting for complex behavior (e.g., anisotropic hardening, strength differential effects, and tension-compression asymmetry).
	
	\begin{remark}
		Although Lou and coworkers \citep{lou:2022a} mention an ABAQUS/Explicit model implementation and show sounding numerical results, to the best of the authors' knowledge neither the formulation nor the implementation details of such a model have been published in the aforementioned contributions. Nonetheless, we propose a fully implicit formulation in the present paper, in addition to improving the model robustness with respect to the yield surface apex singularity.
	\end{remark}
	
	The following section establishes the main ingredients of the elasto-plastic model formulation and provides all the required details for a fully implicit computational implementation. The formulation follows closely the formalism of thermodynamics with internal variables and overall computational framework described by Souza Neto and coworkers \citep{desouzaneto:2008a}.
	
	\subsection{Formulation \label{ssec:lzy_formulation}}
\begin{figure}[hbt]
		\centering
		\begin{subfigure}[b]{0.48\textwidth}
			\centering
			\includegraphics[width=\textwidth]{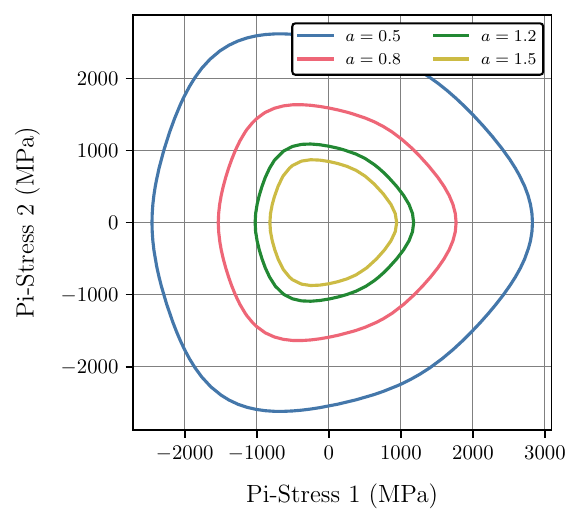}
			\caption{}
			\label{subfig:lou_yield_a_parametric}
		\end{subfigure}
		\begin{subfigure}[b]{0.47\textwidth}
			\centering
			\includegraphics[width=\textwidth]{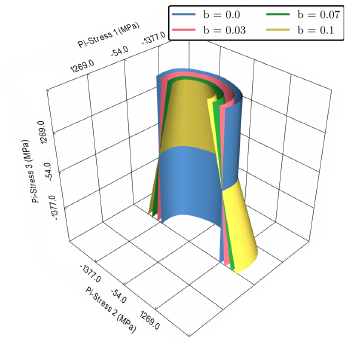}
			\caption{}
			\label{subfig:lou_yield_b_parametric}
		\end{subfigure}\hfill
		\begin{subfigure}[b]{0.47\textwidth}
			\centering
			\includegraphics[width=\textwidth]{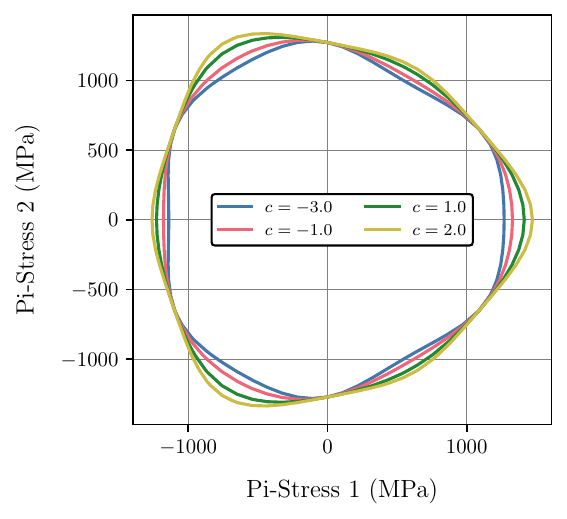}
			\caption{}
			\label{subfig:lou_yield_c_parametric}
		\end{subfigure}
		\begin{subfigure}[b]{0.495\textwidth}
			\centering
			\includegraphics[width=\textwidth]{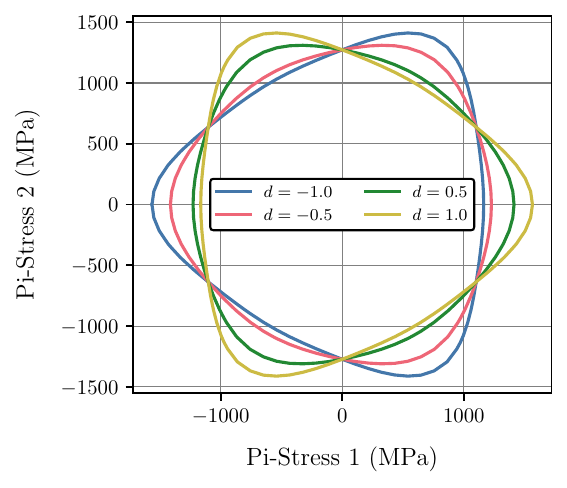}
			\caption{}
			\label{subfig:lou_yield_d_parametric}
		\end{subfigure}\hfill
		\caption{The highly versatile LZY model yield surface as a function of the four yield parameters, $(a, \, b, \, c, \, d)$: \subref{subfig:lou_yield_a_parametric} $(a, \, b, \, c, \, d) = (\Delta, \, 0.05, \, 1.0, \, 0.5)$; \subref{subfig:lou_yield_b_parametric} $(a, \, b, \, c, \, d) = (1.0, \, \Delta, \, 1.0, \, 0.5)$; \subref{subfig:lou_yield_c_parametric} $(a, \, b, \, c, \, d) = (1.0, \, 0.05, \, \Delta, \, 0.5)$; \subref{subfig:lou_yield_d_parametric} $(a, \, b, \, c, \, d) = (1.0, \, 0.05, \, 1.0, \, \Delta)$. The default parameters are $(a, \, b, \, c, \, d) = (1.0, \, 0.05, \, 1.0, \, 0.5)$ and $\Delta$ denotes the only yield parameter being updated.}
		\label{fig:lou_yield_parametric}
	\end{figure}

	The essential elements of the LZY elasto-plastic model can be summarized as follows:
	\begin{itemize}
		\item \textbf{Additive strain decomposition.} The total strain tensor, $\bm{\varepsilon}$, is additively decomposed as
		\begin{equation}
			\bm{\varepsilon} = \bm{\varepsilon}^{e} + \bm{\varepsilon}^{p} \, ,
		\end{equation}
		where $\bm{\varepsilon}^{e}$ and $\bm{\varepsilon}^{p}$ are the elastic and plastic strain tensors, respectively.
		
		\item \textbf{Linear elastic law.} A linear elastic law is assumed as
		\begin{equation}
			\bm{\sigma} = \bm{\mathsf{D}}^{e} : \bm{\varepsilon}^{e} \, ,
		\end{equation}
		where $\bm{\mathsf{D}}^{e}$ is the fourth-order isotropic elasticity tensor.
		\item \textbf{Yield criterion.} The LZY yield function is established as
		\begin{equation}
			\Phi (\bm{\sigma} \, ; \, a, \, b, \, c, \, d, \, \sigma_{y}) = a  \left\{ b I_{1} + \left[ (J_{2}^{3} - c J_{3}^{2})^{\frac{1}{2}} - d J_{3} \right]^{\frac{1}{3}} \right\} - \sigma_{y} (\bar{\varepsilon}^{\, p}) \, ,
			\label{eq:lzy_yield_surface_model}
		\end{equation}
		where $\bar{\varepsilon}^{\, p}$ is the accumulated plastic strain and $\sigma_{y}(\bar{\varepsilon}^{p})$ denotes the isotropic strain hardening curve. As illustrated in Figure~\ref{fig:lou_yield_parametric}, each yield parameter plays a different role on the yield surface: (i) parameter $a$ controls the size, (ii) parameter $b$ controls the pressure dependency, (iii) parameter $c$ controls the curvature, and (iv) parameter $d$ controls the strength differential effect. It is remarked that an independent hardening rule can be also postulated for each yield parameter, i.e., $a(\bar{\varepsilon}^{p}), \, b(\bar{\varepsilon}^{p}), \, c(\bar{\varepsilon}^{p}), \, d(\bar{\varepsilon}^{p})$, allowing a complex evolution of the yield surface with plastic strain. In addition, it is convenient to define the LZY effective stress as
		\begin{equation}
			\bar{q} =  a  \left\{ b I_{1} + \left[ (J_{2}^{3} - c J_{3}^{2})^{\frac{1}{2}} - d J_{3} \right]^{\frac{1}{3}} \right\}  \, ,
			\label{eq:lzy_effective_stress}
		\end{equation}
		and to note that the yield surface apex pressure is given by
		\begin{equation}
			p^{\mathrm{apex}} = \dfrac{1}{3 a b} \, \sigma_{y} (\bar{\varepsilon}^{\, p}) \, ,
			\label{eq:lzy_apex_equation}
		\end{equation}
		which can be easily obtained from Equation~\eqref{eq:lzy_yield_surface_model} by taking into account that neither $c$ or $d$ influence the yield surface apex. With the yield function at hand, the elastic domain and the plastically admissible domains are then, respectively, defined as
		\begin{equation}
			\mathcal{E} = \left\{ \bm{\sigma} \, | \, \Phi (\bm{\sigma}) < 0 \right\} \, , \qquad \bar{\mathcal{E}} = \left\{ \bm{\sigma} \, | \, \Phi (\bm{\sigma}) \leq 0 \right\} \, ,
		\end{equation}
		where $\mathcal{Y} = \left\{ \bm{\sigma} \, | \, \Phi (\bm{\sigma}) = 0 \right\}$ defines the yield surface of the model.
		\item \textbf{Plastic flow rule.} An associative plastic flow rule is postulated by taking the yield function as the flow potential, i.e., $\bm{\Psi} \equiv \bm{\Phi}$, such that
		\begin{equation}
			\dot{\bm{\varepsilon}}^{p} = \dot{\gamma} \, \bm{N} \, , \qquad \bm{N} = \dfrac{\partial \bm{\Phi}(\bm{\sigma})}{\partial \bm{\sigma}} \, ,
		\end{equation}	 
		where $\dot{\gamma}$ is the plastic multiplier and $\bm{N}$ denotes the plastic flow vector. Associativity implies that the plastic strain rate, $\dot{\bm{\varepsilon}}^{p}$, is a tensor normal to the yield surface in the stress space. However, in resemblance with the well-known Drucker-Prager model, the LZY yield surface is smooth everywhere but has a singularity at the apex. This means that the flow vector at the apex singularity is an element of the subdifferential of the yield function, $\bm{N} \in \partial_{\sigma} \bm{\Phi}$, and lies within the complementary cone to the LZY yield surface.
		\item \textbf{Hardening rule.} The hardening law is postulated as
		\begin{equation}
			\dot{\bar{\bm{\varepsilon}}}^{\, p} = \omega \, \dot{\gamma} \, ,
		\end{equation}
		where $\omega$ is a constant parameter\footnote{The parameter $\omega$ is solely introduced here as a convenience to recover the hardening rule of well-known constitutive models, being set as $\omega=1.0$ by default. For instance, while $\omega=1.0$ recovers the associative hardening rule of the von Mises elasto-plastic model, setting $\omega = \xi$, where $\xi$ is yield surface cohesion parameter, recovers the Drucker-Prager associative hardening rule.}. For the particular case of the plastic flow at the yield surface apex, a Drucker-Prager-like hardening rule is adopted instead as
		\begin{equation}
			\dot{\bar{\bm{\varepsilon}}}^{\, p} = \dfrac{\xi}{\eta} \, \dot{\varepsilon}^{\, p}_{\mathrm{vol}} \, ,
			\label{eq:lzy_hardening_apex}
		\end{equation}
		where $\varepsilon^{\, p}_{\mathrm{vol}}$ denotes the volumetric plastic strain, and $\xi$ and $\eta$ can be determined from the equivalence between LZY and Drucker-Prager yield surfaces as
		\begin{equation}
			\xi = \dfrac{2 \sqrt{3}}{3}\sqrt{1 - \dfrac{1}{3} \eta^{2}} \, , \qquad \eta = 3 a b \, .
		\end{equation}
		\item \textbf{Loading/Unloading conditions.} The model is completed by setting the so-called loading/unloading conditions, i.e., the constraints that establish when plastic flow may occur,
		\begin{equation}
			\Phi (\bm{\sigma}) \leq 0 \, , \qquad \dot{\gamma} \geq 0 \, , \qquad \dot{\gamma} \, \Phi (\bm{\sigma}) = 0 \, .
		\end{equation}
	\end{itemize}
	
	\subsection{State update and consistent tangent operator}
	
	By adopting a backward (or fully implicit) Euler scheme to discretize the LZY model constitutive equations, the incremental elasto-plastic constitutive problem in the generic (pseudo-)time interval $[t_{n}, \, t_{n+1}]$ is stated as follows: Given the values $\bm{\varepsilon}^{e}_{n}$ and $\bar{\bm{\varepsilon}}^{\,  p}_{n}$, and given the prescribed incremental strain $\Delta \bm{\varepsilon}$, solve the following system of equations
	\begin{align}
		\bm{\varepsilon}^{e}_{n+1} &= \bm{\varepsilon}^{e}_{n} + \Delta \bm{\varepsilon} - \Delta \gamma \, \bm{N}_{n+1} \, , \\[5pt]
		\bar{\bm{\varepsilon}}^{\, p}_{n+1} &= \bar{\bm{\varepsilon}}^{\, p}_{n} + \omega \Delta \gamma \, , 
	\end{align}
	for the unknowns $\bm{\varepsilon}^{e}_{n+1}$, $\bar{\bm{\varepsilon}}^{\, p}_{n+1}$ and $\Delta \gamma$, subjected to the constraints
	\begin{equation}
		\Delta \gamma \geq 0 \, , \qquad \bm{\Phi}(\bm{\sigma}_{n+1})  \leq 0 \, , \qquad \Delta \gamma \, \bm{\Phi}(\bm{\sigma}_{n+1}) = 0 \, ,
	\end{equation}
	where
	\begin{equation}
		\bm{\sigma}_{n+1} = \bm{\mathsf{D}}^{e} : \bm{\varepsilon}^{e}_{n+1} \, .
	\end{equation}
	
	\clearpage
	
	The previous problem can be conveniently solved by the well-known elastic prediction/return-mapping algorithm. First, the so-called elastic trial state ($\bm{\varepsilon}^{e, \, \mathrm{trial}}_{n+1}$, $\bar{\bm{\varepsilon}}^{\, p, \, \mathrm{trial}}_{n+1}$, $\bm{\sigma}^{\mathrm{trial}}_{n+1}$) is computed assuming that $\Delta \gamma = 0$. If the elastic trial state lies within the elastic domain, $\mathcal{E}$, then it is accepted as the problem solution at $t_{n+1}$. Otherwise, the elastic trial state is not plastically admissible and a return-mapping systems of equations must be solved to determine the state ($\bm{\varepsilon}^{e}_{n+1}$, $\bar{\bm{\varepsilon}}^{\, p}_{n+1}$, $\bm{\sigma}_{n+1}$), depending on whether the plastic flow occurs at the smooth yield surface or at the yield surface apex. To determine which of the return-mapping problems must be solved\footnote{The solution of the Drucker-Prager model follows a different strategy. First, the return-mapping to the smooth yield surface is solved. If the solution is admissible, then it is accepted as the problem solution at $t_{n+1}$. Otherwise, the return-mapping to the yield surface apex is solved to find the solution instead.}, a simple but efficient criterion is proposed based on the trial pressure, $p^{\mathrm{trial}}_{n+1} = \frac{1}{3} \mathrm{tr}(\bm{\sigma}^{\mathrm{trial}}_{n+1})$.
	If the trial pressure is greater than the apex pressure, $p^{\mathrm{trial}}_{n+1} > p^{\mathrm{apex}}_{n+1}$, then the return-mapping to the yield surface apex is solved, otherwise the return-mapping to the smooth yield surface is solved instead\footnote{Note that this simple criterion can be deemed `approximate' as it discards some admissible solutions to the smooth yield surface for which the elastic trial state lies outside the complementary cone. Nevertheless, the approximation holds well for most of the pressure dependency levels often found in practice, but numerical convergence issues in the transition to the apex can be alleviated by introducing a switching tolerance, $\epsilon^{\mathrm{apex}} > 0$, in the criterion, such that $p^{\mathrm{trial}}_{n+1} > (1.0 - \epsilon^{\mathrm{apex}}) \, p^{\mathrm{apex}}_{n+1}$.}.
	
	The return-mapping to the smooth yield surface comprises the solution of the following system of residual equations
	\begin{equation}
		\begin{cases}
			\bm{R}_{1} = \bm{\varepsilon}^{e}_{n+1} - \bm{\varepsilon}^{e, \, \mathrm{trial}}_{n+1} + \Delta \gamma \, \bm{N}_{n+1} = \bm{0}  \\[10pt]
			R_{2} = \bar{\bm{\varepsilon}}^{\, p}_{n+1} - \bar{\bm{\varepsilon}}^{\, p}_{n} - \omega \Delta \gamma = 0 \\[10pt]
			R_{3} = \dfrac{\bar{q}_{n+1} - \sigma_{y} (\bar{\varepsilon}^{\, p}_{n+1})}{\sigma_{y} (\bar{\varepsilon}^{\, p}_{0})} = 0
		\end{cases} , \qquad w.r.t. \quad \bm{\alpha}_{n+1} = \left\{ \bm{\varepsilon}_{n+1}^{e}, \,  \bar{\bm{\varepsilon}}^{\, p}_{n+1}, \, \Delta \gamma \right\} \, ,
		\label{eq:lzy_retmap_smooth}
	\end{equation}
	with initial conditions $\bm{\varepsilon}^{e}_{n+1} = \bm{\varepsilon}^{e, \, \mathrm{trial}}_{n+1}$, $\bar{\bm{\varepsilon}}^{\, p}_{n+1} = \bar{\bm{\varepsilon}}^{\, p}_{n}$ and $\Delta \gamma = 0$. The solution can be efficiently found through the well-known Newton-Raphson method, where each iteration $(k)$ consists in finding the solution of the linearized version of Equation~\eqref{eq:lzy_retmap_smooth} as
	\begin{equation}
		\bm{\upalpha}^{(k)}_{n+1} = \bm{\upalpha}^{(k-1)}_{n} + \delta \bm{\upalpha}^{(k)} \, , \qquad \delta   \bm{\upalpha}^{(k)} = - \mathbf{J}^{-1}\left(\bm{\upalpha}^{(k-1)}_{n+1} \right) \mathbf{R}\left(\bm{\upalpha}^{(k-1)}_{n+1} \right) \, ,
		\label{eq:lzy_nr_method}
	\end{equation}
	where $\bm{\upalpha}^{(k)}_{n+1}$ is the matricial form of the vector of unknowns, $\mathbf{R}$ is the vectorial form of the joint residual function, and $\mathbf{J}$ is the matricial form of the Jacobian matrix
	\begin{equation}
		\bm{J} = \dfrac{\partial \bm{R}}{\partial \bm{\alpha}_{n+1}} = \begin{bmatrix}
			\bm{J}_{11} & \bm{J}_{12} & \bm{J}_{13} \\
			\bm{J}_{21} & J_{22} & J_{23} \\
			\bm{J}_{31} & J_{32} & J_{33}
		\end{bmatrix}
		= \begin{bmatrix}
			\dfrac{\partial \bm{R}_{1}}{\partial  \bm{\varepsilon}_{n+1}^{e}} &  \dfrac{\partial \bm{R}_{1}}{\partial  \bar{\bm{\varepsilon}}^{\, p}_{n+1}} & \dfrac{\partial \bm{R}_{1}}{\partial  \Delta \gamma} \\[20pt]
			\dfrac{\partial R_{2}}{\partial  \bm{\varepsilon}_{n+1}^{e}} &  \dfrac{\partial R_{2}}{\partial  \bar{\bm{\varepsilon}}^{\, p}_{n+1}} & \dfrac{\partial R_{2}}{\partial  \Delta \gamma} \\[20pt]
			\dfrac{\partial R_{3}}{\partial  \bm{\varepsilon}_{n+1}^{e}} &  \dfrac{\partial R_{3}}{\partial  \bar{\bm{\varepsilon}}^{\, p}_{n+1}} & \dfrac{\partial R_{3}}{\partial  \Delta \gamma}
		\end{bmatrix} \, . \label{eq:lzy_retmap_smooth_jac}
	\end{equation}
	
	\clearpage
	
	The return-mapping to the yield surface apex involves the solution of the residual equation
	\begin{equation}
		R_{4} = \sigma_{y} (\bar{\varepsilon}^{\, p}_{n+1}) \beta - p^{\mathrm{trial}}_{n+1} + K \Delta \varepsilon^{\, p}_{\mathrm{vol}} = 0 \, , \quad  \bar{\varepsilon}^{\, p}_{n+1} = \bar{\varepsilon}^{\, p}_{n} + \alpha \Delta \varepsilon^{\, p}_{\mathrm{vol}} \, , \quad w.r.t.  \quad \bm{\alpha}_{n+1} = \left\{ \Delta \varepsilon^{\, p}_{\mathrm{vol}} \right\}  \, ,
	\end{equation}
	where $\alpha = \xi/\eta$ and $\beta = (3ab)^{-1}$. The previous residual equation is obtained by discretizing the hardening rule in Equation~\eqref{eq:lzy_hardening_apex} with the backward Euler scheme,
	\begin{equation}
		\bar{\varepsilon}^{p}_{n+1} = \bar{\varepsilon}^{p}_{n} + \alpha \Delta \varepsilon^{p}_{\mathrm{vol}} \, ,
	\end{equation}
	and replacing the volumetric component of the linear elastic law\footnote{In this case, note that both the elastic strain and the stress are purely volumetric,
		\begin{equation*}
			\bm{\varepsilon}^{e}_{n+1} = \dfrac{1}{3K} \, p_{n+1} \bm{I} \, , \qquad \bm{\sigma}_{n+1} = p_{n+1} \bm{I} \, ,
		\end{equation*}
		where $K$ denotes the elastic bulk modulus and $\bm{I}$ is the second-order identity tensor.},
	\begin{equation}
		p_{n+1} = p^{\mathrm{trial}}_{n+1} - K \Delta \varepsilon^{\, p}_{\mathrm{vol}} \, ,
	\end{equation}
	in the LZY yield surface apex defined by Equation~\eqref{eq:lzy_apex_equation}. Assuming the initial condition $\Delta \varepsilon^{\, p}_{\mathrm{vol}} = 0$, such that $\bar{\bm{\varepsilon}}^{\, p}_{n+1} = \bar{\bm{\varepsilon}}^{\, p}_{n}$, the solution can be once again found with the Newton-Raphson method (see Equation~\eqref{eq:lzy_nr_method}), being the matricial form of the Jacobian matrix simply given by
	\begin{equation}
		\bm{J} = \dfrac{\partial \bm{R}}{\partial \bm{\alpha}_{n+1}} = \begin{bmatrix} \dfrac{\partial R_{4}}{\partial \Delta \varepsilon^{\, p}_{\mathrm{vol}}} \end{bmatrix} \, .
		\label{eq:lzy_retmap_apex_jac}
	\end{equation}
	
	Having established the LZY model state update, it remains to determine the consistent tangent operator, $\bm{\mathsf{D}}_{n+1}$, required to complete the fully implicit implementation. Three different consistent tangent operators are required according with the followed state update path:
	\begin{itemize}
		\item \textbf{Elastic update.} If the incremental step is elastic, where the elastic trial state is accepted as the solution, the consistent tangent operator is the fourth-order isotropic elasticity tensor,
		\begin{equation}	
			\bm{\mathsf{D}}_{n+1} = \bm{\mathsf{D}}^{e} \, .
		\end{equation}
		\item \textbf{Plastic update to smooth yield surface.} If the incremental step is elasto-plastic and results from the return-mapping to the smooth yield surface, the consistent tangent operator is given by
		\begin{equation}
			\bm{\mathsf{D}}_{n+1} = \dfrac{\partial \bm{\sigma}_{n+1}}{\partial \bm{\varepsilon}_{n+1}} = \bm{\mathsf{D}}^{e} : \dfrac{\partial \bm{\varepsilon}_{n+1}^{e}}{\partial \bm{\varepsilon}_{n+1}} \, , \qquad \dfrac{\partial \bm{\varepsilon}_{n+1}^{e}}{\partial \bm{\varepsilon}_{n+1}} = (\bm{J}^{-1})_{11} : \bm{\mathsf{I}}_{s} \, ,
		\end{equation}
		where $\bm{J}^{-1}$ denotes the inverse of the Jacobian matrix (see Equation~\eqref{eq:lzy_retmap_smooth_jac}) and $\bm{\mathsf{I}}_{s}$ denotes the fourth-order symmetric identity operator.
		\item \textbf{Plastic update to yield surface apex.} If the incremental step is elasto-plastic and results from the return-mapping to the yield surface apex, the consistent tangent operator is
		\begin{equation}
			\bm{\mathsf{D}}_{n+1} = \dfrac{\partial \bm{\sigma}_{n+1}}{\partial \bm{\varepsilon}_{n+1}} = K ( 1 - K \bm{J}^{-1}) \, \bm{I} \otimes \bm{I} \, ,
		\end{equation}
		where $\bm{J}^{-1}$ denotes the inverse of the Jacobian matrix (see Equation~\eqref{eq:lzy_retmap_apex_jac}) and $\bm{I}$ denotes the second-order identity tensor.
	\end{itemize}
	
	All derivatives required to (explicitly) compute the flow vector and both state update Jacobian matrices (see Equations~\eqref{eq:lzy_retmap_smooth_jac} and \eqref{eq:lzy_retmap_apex_jac}) are provided for convenience in Section~\ref{ssec:lzy_linearization}.

	\subsection{Convexity return-mapping based on GINCA method \label{ssec:convexity_ginca}}
\begin{figure}[hbt]
		\centering
		\includegraphics[width=0.95\textwidth]{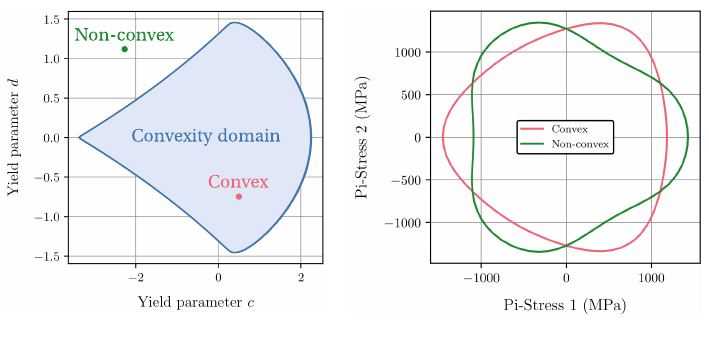}
		\caption{Convexity domain of the LZY model yield surface. If given set $(c, \, d)$ falls outside the convexity domain, then the yield surface is non-convex, otherwise convexity is guaranteed.}
		\label{fig:lzy_convexity_domain}
	\end{figure}

	From a numerical point of view, it is well-known that the convexity of the yield surface should be guaranteed to ensure the uniqueness of the state update solution (i.e., avoiding that the yield surface is `pierced' as the result of a given strain increment). Therefore, it is important to highlight that the convexity of the LZY yield surface depends on the set of yield parameters $(c, \, d)$ (see Figure~\ref{fig:lzy_convexity_domain}) and that Lou and coworkers \citep{lou:2022a} proposed a numerical method called Geometry-Inspired Numerical Convex Analysis (GINCA) to check if a given set $(c, \, d)$ corresponds to a convex LZY yield surface.
	
	However, in the context of the present paper it is not only necessary to know if the LZY yield surface is convex, but also to actually enforce such convexity during the model discovery process\footnote{When independent hardening rules are postulated for the yield parameters $c$ and $d$, i.e., $ c=c(\bar{\varepsilon}^{p}), \, d=d(\bar{\varepsilon}^{p})$, enforcing convexity is also relevant in the state update algorithm.}. This motivated the proposal of the convexity return-mapping described below, based on the GINCA method, where the convexity domain boundary can be roughly interpreted as the yield surface in the elastic prediction/return-mapping algorithm.
	
\begin{figure}[hbt]
		\centering
		\begin{subfigure}[b]{0.490\textwidth}
			\centering
			\includegraphics[width=\textwidth]{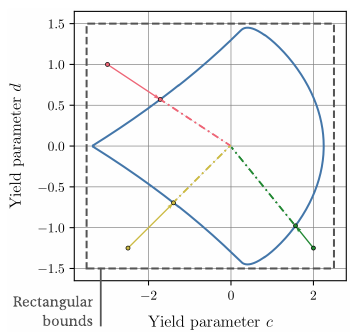}
			\caption{}
			\label{subfig:lzy_convexity_retmap_points}
		\end{subfigure}
		\begin{subfigure}[b]{0.49\textwidth}
			\centering
			\includegraphics[width=\textwidth]{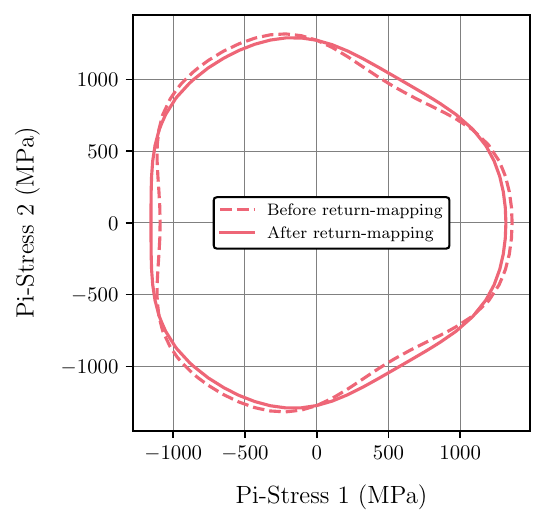}
			\caption{}
			\label{subfig:lzy_convexity_retmap_p1}
		\end{subfigure}\hfill
		\begin{subfigure}[b]{0.490\textwidth}
			\centering
			\includegraphics[width=\textwidth]{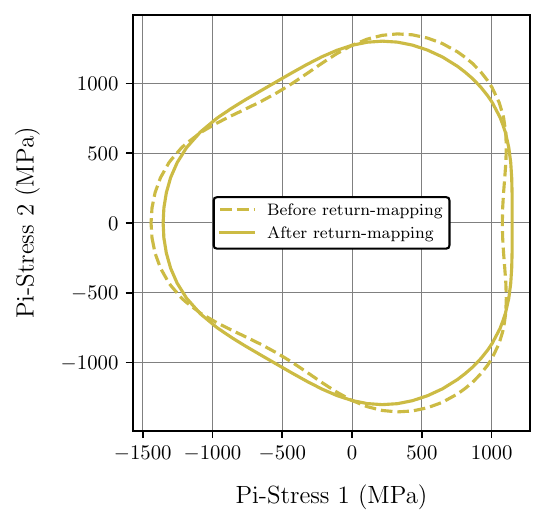}
			\caption{}
			\label{subfig:lzy_convexity_retmap_p2}
		\end{subfigure}
		\begin{subfigure}[b]{0.495\textwidth}
			\centering
			\includegraphics[width=\textwidth]{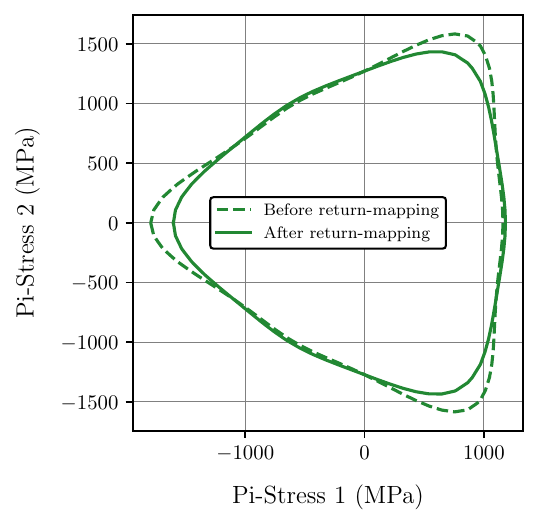}
			\caption{}
			\label{subfig:lzy_convexity_retmap_p3}
		\end{subfigure}\hfill
		\caption{LZY model yield surface convexity return-mapping: \subref{subfig:lzy_convexity_retmap_points} Return-mapping of three different sets $(c^{*}, \, d^{*})$ to the convexity domain; \subref{subfig:lzy_convexity_retmap_p1}-\subref{subfig:lzy_convexity_retmap_p3} Comparison between the yield surface before and after convexity return-mapping.}
		\label{fig:lzy_convexity_retmap}
	\end{figure}

	Assume that the set $(c^{*}, \, d^{*})$ results from a given optimization step. If $(c^{*}, \, d^{*})$ lies inside the convexity domain as determined by the GINCA method, then it is accepted as the admissible solution, $(c, \, d) = (c^{*}, \, d^{*})$. Otherwise, the GINCA method is employed to determine the set $(c^{*}, \, d^{*})$ in the convex domain boundary along the direction given by $\mathrm{tan}(\theta) = d/c$, thus `returning' $(c, \, d)$ to the convex domain. To further improve the efficiency of this method, rectangular bounds containing the whole convexity domain can be first enforced ($c \in [-3.5, \, 2.5]$, $d \in [-1.5, \, 1.5]$), and then the GINCA method is applied with an upper search radius of $R_u = 5.0$. The convexity return-mapping is illustrated in Figure~\ref{fig:lzy_convexity_retmap}.
	
	\clearpage
	
	\subsection{Linearization \label{ssec:lzy_linearization}}
	The derivatives required to implement the fully implicit LZY model are summarized below.
	
	\begin{remark}
		The following linearization assumes the general case where an independent hardening rule is postulated for each yield parameter, i.e., $a=a(\bar{\varepsilon}^{p}), \, b=b(\bar{\varepsilon}^{p}), \, c=c(\bar{\varepsilon}^{p}), \, d=d(\bar{\varepsilon}^{p})$.
	\end{remark}
	
	\subsubsection{Auxiliary variables and derivatives}
	
	For the sake of compactness, the increment indexing is omitted and the following auxiliary variables are defined
	\begin{equation}
		w_{1} = b I_{1} \, , \qquad w_{2} = c J_{3}^{2} \, , \qquad  w_{3} = d J_{3} \, , \qquad w_{4} = J_{2}^{3} - w_{2} \, , \qquad w_{5} = w_{4}^{\frac{1}{2}} - w_{3} \, .
	\end{equation}
	
	The corresponding first- and second-order derivatives w.r.t. stress are
	\begin{equation}
		\begin{gathered}
			\dfrac{\partial w_{1}}{\partial \bm{\sigma}} = b \dfrac{\partial I_{1}}{\partial \bm{\sigma}} \, , \qquad \dfrac{\partial w_{2}}{\partial \bm{\sigma}} = 2 c J_{3} \dfrac{\partial J_{3}}{\partial \bm{\sigma}} \, , \qquad \dfrac{\partial w_{3}}{\partial \bm{\sigma}} = d \dfrac{\partial J_{3}}{\partial \bm{\sigma}} \, , \qquad \dfrac{\partial w_{4}}{\partial \bm{\sigma}} = 3 J_{2}^{2} \, \dfrac{\partial J_{2}}{\partial \bm{\sigma}} - \dfrac{\partial w_{2}}{\partial \bm{\sigma}} \, , \\[10pt]
			\dfrac{\partial w_{5}}{\partial \bm{\sigma}} = \dfrac{1}{2} w_{4}^{-\frac{1}{2}} \, \dfrac{\partial w_{4}}{\partial \bm{\sigma}} - \dfrac{\partial w_{3}}{\partial \bm{\sigma}} \, .
		\end{gathered}
	\end{equation}
	and
	\begin{equation}
		\begin{gathered}
			\dfrac{\partial^{2} w_{2}}{\partial \bm{\sigma}^{2}} = 2c \left( \dfrac{\partial J_{3}}{\partial \bm{\sigma}} \otimes \dfrac{\partial J_{3}}{\partial \bm{\sigma}} + J_{3} \, \dfrac{\partial^{2} J_{3}}{\bm{\sigma}^{2}}\right) \, , \qquad \dfrac{\partial^{2} w_{3}}{\partial \bm{\sigma}^{2}} = d \, \dfrac{\partial^{2} J_{3}}{\partial \bm{\sigma}^{2}} \, , \\[10pt]
			\dfrac{\partial^{2} w_{4}}{\partial \bm{\sigma}^{2}} = 6 J_{2} \, \dfrac{\partial J_{2}}{\partial \bm{\sigma}} \otimes \dfrac{\partial J_{2}}{\partial \bm{\sigma}} + 3 J_{2}^{2} \, \dfrac{\partial^{2} J_{2}}{\partial \bm{\sigma}^{2}} - \dfrac{\partial^{2} w_{2}}{\partial \bm{\sigma}^{2}} \, , \\[10pt]
			\dfrac{\partial^{2} w_{5}}{\partial \bm{\sigma}^{2}} = -\dfrac{1}{4} w_{4}^{-\frac{3}{2}} \, \dfrac{\partial w_{4}}{\partial \bm{\sigma}} \otimes \dfrac{\partial w_{4}}{\partial \bm{\sigma}} + \dfrac{1}{2} w_{4}^{-\frac{1}{2}} \, \dfrac{\partial^{2} w_{4}}{\partial \bm{\sigma}^{2}} - \dfrac{\partial^{2} w_{3}}{\partial \bm{\sigma}^{2}} \, .
		\end{gathered}
	\end{equation}
	
	In addition, the first- and cross second-order derivatives w.r.t. accumulated plastic strain are
	\begin{equation}
		\begin{gathered}
			\dfrac{\partial w_{1}}{\partial \bar{\varepsilon}^{p}} = I_{1} \dfrac{\partial b}{\partial \bar{\varepsilon}^{p}} \, , \qquad
			\dfrac{\partial w_{2}}{\partial \bar{\varepsilon}^{p}} = J_{3}^{2} \, \dfrac{\partial c}{\partial \bar{\varepsilon}^{p}} \, , \qquad \dfrac{\partial w_{3}}{\partial \bar{\varepsilon}^{p}} = J_{3} \, \dfrac{\partial d}{\partial \bar{\varepsilon}^{p}} \, , \qquad \dfrac{\partial w_{4}}{\partial \bar{\varepsilon}^{p}} = - \dfrac{\partial w_{2}}{\partial \bar{\varepsilon}^{p}} \, , \\[10pt] \dfrac{\partial w_{5}}{\partial \bar{\varepsilon}^{p}} = \dfrac{1}{2} w_{4}^{-\frac{1}{2}} \, \dfrac{\partial w_{4}}{\partial \bar{\varepsilon}^{p}} - \dfrac{\partial w_{3}}{\partial \bar{\varepsilon}^{p}} \, ,
		\end{gathered}
	\end{equation}
	and
	\begin{equation}
		\begin{gathered}
			\dfrac{\partial^{2} w_{1}}{\partial \bm{\sigma} \partial \bar{\varepsilon}^{p}} = \dfrac{\partial b}{\partial \bar{\varepsilon}^{p}} \bm{I} \, , \qquad	
			\dfrac{\partial^{2} w_{2}}{\partial \bm{\sigma} \partial \bar{\varepsilon}^{p}} = 2 J_{3} \, \dfrac{\partial c}{\partial \bar{\varepsilon}^{p}} \dfrac{\partial J_{3}}{\partial \bm{\sigma}} \, , \qquad \dfrac{\partial^{2} w_{3}}{\partial \bm{\sigma} \partial \bar{\varepsilon}^{p}} = \dfrac{\partial d}{\partial \bar{\varepsilon}^{p}} \dfrac{\partial J_{3}}{\partial \bm{\sigma}} \, , \\[10pt] \dfrac{\partial^{2} w_{4}}{\partial \bm{\sigma} \partial \bar{\varepsilon}^{p}} = - \dfrac{\partial^{2} w_{2}}{\partial \bm{\sigma} \partial \bar{\varepsilon}^{p}} \, , \qquad
			\dfrac{\partial^{2} w_{5}}{\partial \bm{\sigma} \partial \bar{\varepsilon}^{p}} = -\dfrac{1}{4} w_{4}^{-\frac{3}{2}} \, \dfrac{\partial w_{4}}{\partial \bar{\varepsilon}^{p}} \dfrac{\partial w_{4}}{\partial \bm{\sigma}} + \dfrac{1}{2} w_{4}^{-\frac{1}{2}} \dfrac{\partial^{2} w_{4}}{\partial \bm{\sigma} \partial \bar{\varepsilon}^{p}} - \dfrac{\partial^{2} w_{3}}{\partial \bm{\sigma} \partial \bar{\varepsilon}^{p}}
		\end{gathered}
	\end{equation}
	
	The first- and second-order derivatives of the stress invariants w.r.t. stress are
	\begin{equation}
		\dfrac{\partial I_{1}}{\partial \bm{\sigma}} = \bm{I} \, , \qquad \dfrac{\partial J_{2}}{\partial \bm{\sigma}} = \bm{\sigma}_{d} \, , \qquad \dfrac{\partial J_{3}}{\partial \bm{\sigma}} = \mathrm{det}(\bm{\sigma}_{d}) \, \bm{\sigma}_{d}^{-1} : \bm{\mathsf{I}}_d \, ,
	\end{equation}
	where $\bm{\sigma}_{d}$ denotes the deviatoric stress tensor and 
	$\bm{\mathsf{I}}_d$ denotes the fourth-order deviatoric projection tensor, and
	\begin{equation}
		\dfrac{\partial^{2} J_{2}}{\partial \bm{\sigma}^{2}} = \bm{\mathsf{I}}_d \, , \qquad \dfrac{\partial^{2} J_{3}}{\partial \bm{\sigma}^{2}} = \left( \bm{\sigma}_{d}^{-1} : \bm{\mathsf{I}}_d \right) \otimes \dfrac{\partial J_{3}}{\partial \bm{\sigma}} +
		\mathrm{det}(\bm{\sigma}_{d}) \left(  \bm{\mathsf{I}}_d :\dfrac{\partial \bm{\sigma}_{d}^{-1}}{\partial \bm{\sigma}} : \bm{\mathsf{I}}_d  \right) \, ,
	\end{equation}
	where
	\begin{equation}
		\left( \dfrac{\partial \bm{\sigma}_{d}^{-1}}{\partial \bm{\sigma}}  \right)_{ijkl} = -\dfrac{1}{2} \Big[ (\bm{\sigma}_{d}^{-1})_{ik} (\bm{\sigma}_{d}^{-1})_{lj} + (\bm{\sigma}_{d}^{-1})_{il} (\bm{\sigma}_{d}^{-1})_{kj} \Big] \, .
	\end{equation}
	
	\subsubsection{Return-mapping to smooth yield surface}
	
	The flow vector results from the derivative of the yield function (see Equation~\eqref{eq:lzy_yield_surface_model}) w.r.t. stress,
	\begin{equation}
		\bm{N} = \dfrac{\partial \bm{\Phi}}{\partial \bm{\sigma}} = a \left( \dfrac{\partial w_{1}}{\partial \bm{\sigma}} + \dfrac{1}{3} w_{5}^{-\frac{2}{3}}  \, \dfrac{\partial w_{5}}{\partial \bm{\sigma}} \right) \, ,
	\end{equation}
	and the corresponding first-order derivatives are
	\begin{equation}
		\begin{gathered}
			\dfrac{\partial \bm{N}}{\partial \bm{\sigma}} = \dfrac{1}{3} a \left( -\dfrac{2}{3} w_{5}^{-\frac{5}{3}} \, \dfrac{\partial w_{5}}{\partial \bm{\sigma}} \otimes \dfrac{\partial w_{5}}{\partial \bm{\sigma}} + w_{5}^{-\frac{2}{3}} \, \dfrac{\partial^{2} w_{5}}{\partial \bm{\sigma}^{2}} \right) \, , \\[10pt]
			\dfrac{\partial \bm{N}}{\partial \bm{\varepsilon}^{e}} = \dfrac{\partial \bm{N}}{\partial \bm{\sigma}} : \bm{\mathsf{D}}^{e} \, , \\[10pt]
			\dfrac{\partial \bm{N}}{\partial \bar{\varepsilon}^{p}} = \dfrac{\partial a}{\partial \bar{\varepsilon}^{p}} \left( \dfrac{\partial w_{1}}{\partial \bm{\sigma}} + \dfrac{1}{3} w_{5}^{-\frac{2}{3}} \, \dfrac{\partial w_{5}}{\partial \bm{\sigma}} \right) + a \left( \dfrac{\partial^{2} w_{1}}{\partial \bm{\sigma} \partial \bar{\varepsilon}^{p}} - \dfrac{2}{9} w_{5}^{-\frac{5}{3}} \, \dfrac{\partial w_{5}}{\partial \bar{\varepsilon}^{p}} \dfrac{\partial w_{5}}{\partial \bm{\sigma}} + \dfrac{1}{3} w_{5}^{-\frac{2}{3}} \dfrac{\partial^{2} w_{5}}{\partial \bm{\sigma} \partial \bar{\varepsilon}^{p}}  \right) \, , \\[10pt]
		\end{gathered}
	\end{equation}
	
	The first-order derivatives of the effective stress (see Equation~\eqref{eq:lzy_effective_stress}) are
	\begin{equation}
		\begin{gathered}
			\dfrac{\partial \bar{q}}{\partial \bm{\varepsilon}^{e}} = \bm{N} : \bm{\mathsf{D}}^{e} \, , \\[10pt]
			\dfrac{\partial \bar{q}}{\partial \bar{\varepsilon}^{p}} = \dfrac{\partial a}{\partial \bar{\varepsilon}^{p}} \left( w_{1} + w_{5}^{\frac{1}{3}} \right) + a \left( \dfrac{\partial w_{1}}{\partial \bar{\varepsilon}^{p}} + \dfrac{1}{3} w_{5}^{-\frac{2}{3}} \dfrac{\partial w_{5}}{\partial \bar{\varepsilon}^{p}} \right) \, .
		\end{gathered}
	\end{equation}
	
	Finally, the first-order derivatives of the residuals required to the compute the Jacobian matrix (see Equation~\eqref{eq:lzy_retmap_smooth_jac}) yield 
	\begin{equation}
		\dfrac{\partial \bm{R}_{1}}{\partial  \bm{\varepsilon}^{e}} = \bm{\mathsf{I}}_{s} + \Delta \gamma \, \dfrac{\partial \bm{N}}{\partial \bm{\varepsilon}^{e}} \, , \qquad \dfrac{\partial \bm{R}_{1}}{\partial  \bar{\varepsilon}^{\, p}} = \Delta \gamma \, \dfrac{\partial \bm{N}}{\partial \bar{\varepsilon}^{\, p}} \, , \qquad \dfrac{\partial \bm{R}_{1}}{\partial  \Delta \gamma} = \bm{N} \, ,
	\end{equation}
	\begin{equation}
		\dfrac{\partial R_{2}}{\partial  \bm{\varepsilon}^{e}} = \bm{0} \, , \qquad \dfrac{\partial R_{2}}{\partial  \bar{\varepsilon}^{\, p}} = 1 \, , \qquad \dfrac{\partial R_{2}}{\partial  \Delta \gamma} = -\omega \, ,
	\end{equation}
	\begin{equation}
		\dfrac{\partial R_{3}}{\partial  \bm{\varepsilon}^{e}} = \dfrac{1}{\sigma_{y, \, 0}} \dfrac{\partial \bar{q}}{\partial  \bm{\varepsilon}^{e}} \, , \qquad \dfrac{\partial R_{3}}{\partial  \bar{\varepsilon}^{\, p}} = \dfrac{1}{\sigma_{y, \, 0}} \left( \dfrac{\partial \bar{q}}{\partial  \bar{\varepsilon}^{\, p}} - \dfrac{\sigma_{y}}{\partial \bar{\varepsilon}^{\, p}} \right) \, , \qquad \dfrac{\partial R_{3}}{\partial \Delta \gamma} = 0 \, .
	\end{equation}

	\subsubsection{Return-mapping to the yield surface apex}
	The first-order derivative of the residual required to the compute the Jacobian matrix (see Equation~\eqref{eq:lzy_retmap_apex_jac}) yields
	\begin{equation}
		\dfrac{\partial R_{4}}{\partial \Delta \bm{\varepsilon}^{p}_{\mathrm{vol}}} = \alpha \beta \, \dfrac{\partial \sigma_{y}}{\partial \bar{\varepsilon}^{\, p}} + K \, .
	\end{equation}
	
	\clearpage
	
	\section{Generation of randomly deformed material patches \label{sec:spdg}}
	
\begin{figure}[hbt]
		\centering
		\begin{subfigure}[b]{0.3\textwidth}
			\centering
			\includegraphics[width=\textwidth]{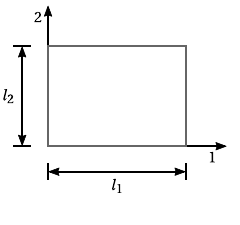}
			\caption{}
			\label{subfig:material_patch_2d_geometry}
		\end{subfigure}
		\begin{subfigure}[b]{0.3\textwidth}
			\centering
			\includegraphics[width=\textwidth]{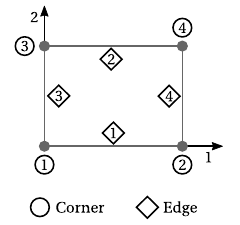}
			\caption{}
			\label{subfig:material_patch_2d_labels}
		\end{subfigure}
		\begin{subfigure}[b]{0.3\textwidth}
			\centering
			\includegraphics[width=\textwidth]{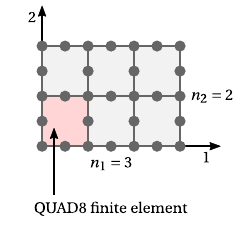}
			\caption{}
			\label{subfig:material_patch_2d_mesh}
		\end{subfigure}\vspace*{5pt} \hfill
		\begin{subfigure}[b]{0.3\textwidth}
			\centering
			\includegraphics[width=\textwidth]{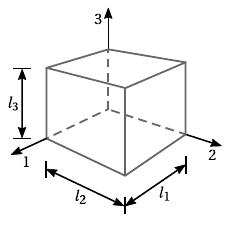}
			\caption{}
			\label{subfig:material_patch_3d_geometry}
		\end{subfigure}
		\begin{subfigure}[b]{0.3\textwidth}
			\centering
			\includegraphics[width=\textwidth]{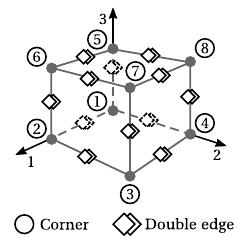}
			\caption{}
			\label{subfig:material_patch_3d_labels}
		\end{subfigure}
		\begin{subfigure}[b]{0.3\textwidth}
			\centering
			\includegraphics[width=\textwidth]{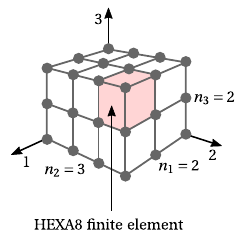}
			\caption{}
			\label{subfig:material_patch_3d_mesh}
		\end{subfigure}\hfill
		\caption{Definition of material patch:
			\subref{subfig:material_patch_2d_geometry} 2D geometry; \subref{subfig:material_patch_2d_labels} 2D corners and edges;
			\subref{subfig:material_patch_2d_mesh} 2D finite element mesh (e.g., 8-node quadratic element);
			\subref{subfig:material_patch_3d_geometry} 3D geometry; \subref{subfig:material_patch_3d_labels} 3D corners and edges;
			\subref{subfig:material_patch_3d_mesh} 3D finite element mesh (e.g., 8-node linear element).}
		\label{fig:material_patch}
	\end{figure}

	This section proposes a new method coined Stochastic Patch Deformation Generator (SPDG) to generate randomly deformed material patches in the context of computational solid mechanics. A material patch is here defined as a quadrilateral (2D) or hexahedral (3D) domain of a given material and with dimensions $l_{i}$, $i=1,\dots, n_{\mathrm{dim}}$ (see Figure~\ref{fig:material_patch}). Moreover, it is assumed that this material patch is discretized in a regular mesh of $n_{e}$ quadrilateral (2D) or hexahedral (3D) finite elements, where $n_{e}=\prod_{i=1}^{n_{\mathrm{dim}}}n_{i}$ and $n_{i}$ denotes the number of elements along dimension $i$.
	
\begin{figure}[hbt]
		\centering
		\begin{subfigure}[b]{0.45\textwidth}
			\centering
			\includegraphics[width=\textwidth]{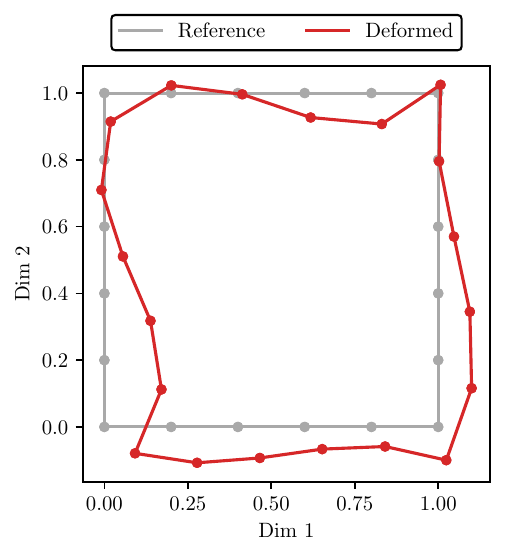}
			\caption{}
			\label{subfig:material_patch_2d_example}
		\end{subfigure}
		\begin{subfigure}[b]{0.48\textwidth}
			\centering
			\includegraphics[width=\textwidth]{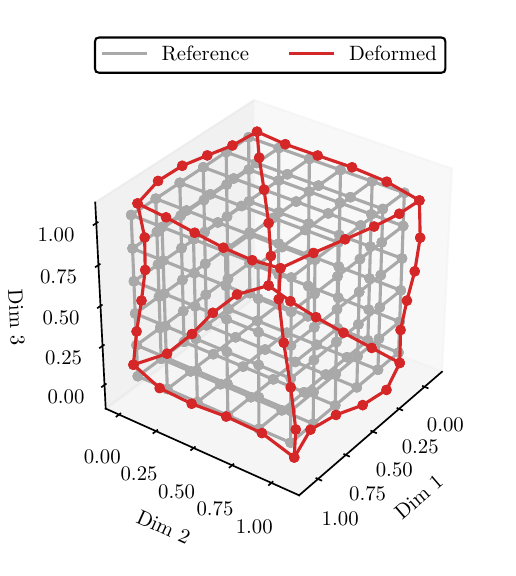}
			\caption{}
			\label{subfig:material_patch_3d_example}
		\end{subfigure}\hfill
		\caption{Example of randomly deformed material patches generated with Stochastic Patch Deformation Generator (SPDG):
			\subref{subfig:material_patch_2d_example} 2D; \subref{subfig:material_patch_3d_example} 3D.}
		\label{fig:material_patch_examples}
	\end{figure}

	Given a set of geometrical constraints, the SPDG method generates a material patch deformed configuration by applying random displacements to its corners and deforming the edges by randomly sampling a polynomial of any order (see Figure~\ref{fig:material_patch_examples}). After selecting a given material model, it transpires that the output  of SPDG can be readily translated into a finite element simulation of a material patch with Dirichlet boundary conditions. Therefore, two immediate applications can be identified:
	\begin{itemize}
		\item \textbf{Material data sets for surrogate modeling.} Large material response data sets can be built by generating samples of randomly deformed material patches with SPDG and solving the corresponding finite element simulations. The finite element mesh, and the nodal and global values stemming from the numerical solution, as well as material model specific quantities, can then be used to compute the desired features to build the surrogate model with a given architecture (e.g., Feed-forward Neural Network, Recurrent Neural Network, Graph Neural Network);
		\item \textbf{Material model robustness testing.} The robustness of material models is of utmost importance when solving real engineering problems where the material usually undergoes complex multi-axial strain paths. However, new material models are often demonstrated in simple or particular sets of strain paths only, thus preventing an accurate assessment of their performance in a real application scenario. SPDG provides an easy-to-use method to test any material model robustness by generating numerical structural simulations with complex strain loading paths.
	\end{itemize}
	
	In the present paper, SPDG is employed to generate a material specimen promoting a broad coverage of the strain-stress space with a single `idealized' loading test. The fundamental concepts and overall strategy of the SPDG method are described in the following section.
	
	\begin{remark}
		The core concepts of the SPDG method apply to both 2D and 3D scenarios. To streamline the method overview in the following section, we focus the discussion on the 2D case, addressing specific 3D considerations when needed.
	\end{remark}
	
	\subsection{Stochastic Patch Deformation Generator (SPDG)}
	
	Given the material patch dimensions, a set of geometrical constraints, and a finite element mesh discretization, the SPDG method generation process involves three main sequential steps: (1) the random sampling of the corners displacements, (2) the random sampling of the edges polynomial deformation, and (3) the random sampling of rigid body translation and/or rotation. It should be remarked that all previous steps are optional, i.e., if the step is not performed then the corresponding displacements are null. Each of the aforementioned steps is described in what follows.
	
	\begin{itemize}
		\item \textbf{Step 1.} The random sampling of the corners displacements is illustrated in Figure~\ref{fig:spdg_step1}. The material patch average deformation bounds must be defined for each corner and along each dimension $i$ as ($\bar{\varepsilon}_{i, \, \mathrm{min}} \, , \bar{\varepsilon}_{i, \, \mathrm{max}}$), $i=1,\dots, n_{\mathrm{dim}}$, where $\bar{\varepsilon}_{i}$ denotes the average deformation with respect to $l_{i}$. Positive/Negative values are associated with a tensile/compressive deformation of the material patch along the corresponding dimension. After $\bar{\varepsilon}_{i}$ is randomly sampled within the previous bounds, the corner node displacement along dimension $i$ is then computed as\footnote{The factor 0.5 accounts for the fact that both corners aligned along a given dimension would have half of the total displacement required to prescribe the material patch average (symmetric) deformation.} $u_{i} = \pm \, 0.5 \, \bar{\varepsilon}_{i} \, l_{i}$, where the sign is determined according to the tensile/compressive nature of the average deformation. Once this process is performed for all corners and dimensions, the corners nodes displacements associated with the material patch deformed configuration are completely defined.

		\begin{figure}[hbt]
			\centering
			\includegraphics[width=0.9\textwidth]{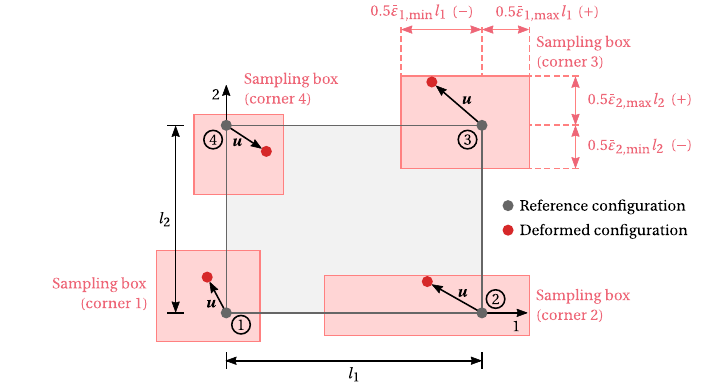}
			\caption{SPDG sampling of material patch corners displacements. Each corner displacement is determined by randomly sampling its position within the corresponding average deformation sampling box.}
			\label{fig:spdg_step1}
		\end{figure}

		\item \textbf{Step 2.} The random sampling of the edges deformation polynomials is illustrated in Figure~\ref{fig:spdg_step2}. It is assumed that the deformation of each edge is a $p$-th order polynomial that intersects both limiting corner nodes. Under this assumption, the polynomial deformation order bounds for each edge must be then defined as $(p_{\mathrm{min}}, \, p_{\mathrm{max}})$, from which a given polynomial order $p$ can be randomly sampled. In turn, the deformation bounds for each edge must be defined as $(\bar{\varepsilon}_{\mathrm{min}}, \, \bar{\varepsilon}_{\mathrm{max}})$. The deformation is always assumed to be orthogonal to the edge deformed orientation and its magnitude is defined with respect to the corresponding patch side length (e.g., the deformation of edge 1 is defined along dimension 2 and with respect to $l_{2}$). Moreover, positive/negative values are associated with a tensile/compressive deformation of the material patch and the edges deformation is measured from the plane defined by the two limiting corner nodes, along the orthogonal dimension.
		
		\clearpage
		
		At this point, the definition of the deformation polynomial falls in one of three cases\footnote{Recall that a general $p$-th order polynomial is completely defined by $p + 1$ control points.}:
		\begin{itemize}
			\item $p=0$. In the particular case of a $0$-th order polynomial, the edge deformation is completely defined by a suitable treatment of the two limiting corner nodes displacements. Given that the edge must remain parallel to its reference configuration, both limiting corner nodes must have the same displacement enforced along the orthogonal dimension. This displacement can be found by averaging both corners displacements determined in step 1;
			\item $p=1$. In the particular case of a $1$-st order polynomial, the edge deformation is completely defined by the two limiting corner nodes displacements determined in step 1;
			\item $p>1$. For the general case of a $(p > 1)$-th order polynomial, the polynomial is randomly computed by sampling the position of the $p - 1$ internal control points. The rectangular sampling domain bounds are defined by (1) the two limiting corner nodes, along the edge deformed orientation, and (2) the range $(0.5 \, \bar{\varepsilon}_{\mathrm{min}} \, l_{\mathrm{orth}}, \, 0.5 \, \bar{\varepsilon}_{\mathrm{max}} \, l_{\mathrm{orth}})$ along the edge orthogonal dimension, where $l_{\mathrm{orth}}$ denotes the material patch size along the edge orthogonal dimension. To avoid extremely sharp deformed configurations, the $p - 1$ internal control points are evenly-spaced along the edge reference orientation. In turn, their coordinate along the edge orthogonal dimension is sampled from a random uniform distribution.
		\end{itemize}
		
		\begin{figure}[hbt]
			\centering
			\includegraphics[width=0.9\textwidth]{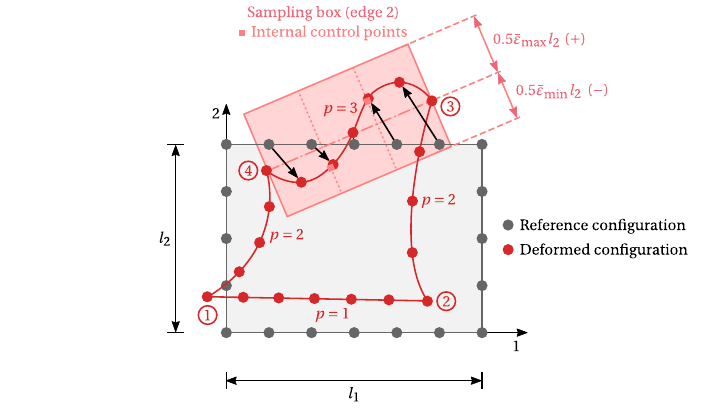}
			\caption{SPDG sampling of material patch edges polynomial deformation. Each edge displacemens are determined by randomly sampling a polynomial deformation order and the required number of internal control points within the corresponding sampling box.  The material patch edge nodes are then projected to the deformation polynomial by enforcing a uniform arc-length spacing.}
			\label{fig:spdg_step2}
		\end{figure}
		
		Having the deformation polynomial completely defined, the material patch edge nodes are then projected  by enforcing a uniform arc-length spacing along the deformed edge and the corresponding displacements are determined. Once this process is performed for all edges, the edges nodes displacement and position in the deformed configuration are completely defined. 
		
		\begin{remark}
			In the 3D case, each material patch face, composed of four edges, is first processed independently following the aforementioned 2D strategy. Each edge is therefore duplicated in the reference configuration, i.e., one edge belonging to each adjacent face (see Figure~\ref{subfig:material_patch_3d_labels}). In addition, the deformation plane of each edge and the corresponding tensile/compressive direction are determined from the face corners deformed coordinates (see Figure~\ref{subfig:spdg_edge_def_plane_3d}). Once all material patch faces are processed, the displacements of each material patch (unique) edge are found by averaging the displacements of the two corresponding edges belonging to adjacent faces (see Figure~\ref{subfig:spdg_double_edge_avg}).

			\begin{figure}[hbt]
				\centering
				\begin{subfigure}[b]{0.495\textwidth}
					\centering
					\includegraphics[width=\textwidth]{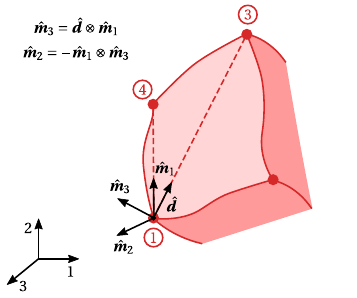}
					\caption{}
					\label{subfig:spdg_edge_def_plane_3d}
				\end{subfigure}
				\begin{subfigure}[b]{0.495\textwidth}
					\centering
					\includegraphics[width=\textwidth]{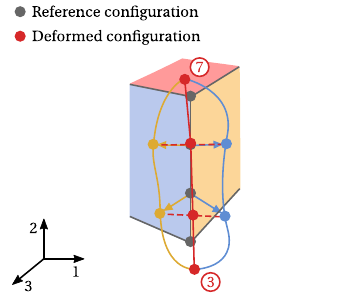}
					\caption{}
					\label{subfig:spdg_double_edge_avg}
				\end{subfigure}\hfill 
				\caption{Specific SPDG procedures required to handle the generation of 3D deformed material patches:
					\subref{subfig:spdg_edge_def_plane_3d} Deformation plane ($\hat{\bm{m}}_{1}\hat{\bm{m}}_{2}$) of given edge (1-4) from the corresponding face corners deformed coordinates, namely the edge reference direction ($\hat{\bm{m}}_{1}$) and face diagonal direction ($\hat{\bm{d}}$); \subref{subfig:spdg_double_edge_avg} Material patch (unique) deformed edge (3-7) determined by averaging the double edge (yellow and blue) nodes displacements.}
				\label{fig:spdg_3d_procedures}
			\end{figure}
			
		\end{remark}
		
		\begin{remark}
			In the 3D case, step 2 may be enriched with the (optional) interpolation of each face node displacements from the corresponding face edges nodes. In particular, the node displacements are found by averaging the interpolated displacements along each dimension. If such interpolation procedure is skipped, then no displacements are prescribed on face nodes.
		\end{remark}
		
		\item \textbf{Step 3.} The random sampling of a rigid body translation and/or rotation is illustrated in Figure~\ref{fig:spdg_step3}. The translation bounds must be defined along each dimension $i$ as $(\bar{u}_{i, \, \mathrm{min}}, \, \bar{u}_{i, \, \mathrm{max}})$, $i=1,\dots, n_{\mathrm{dim}}$, whereas the rotation bounds are defined by means of the Euler angles as $(\bar{\theta}_{i, \, \mathrm{min}}, \, \bar{\theta}_{i, \, \mathrm{max}})$, $i=1, 2, 3$, following the Bunge convention\footnote{The Bunge convention defines Euler angles as a sequence of three rotations about the z-axis, x-axis, and z-axis, here denoted as $(\bar{\theta}_{1}, \, \bar{\theta}_{2}, \, \bar{\theta}_{3})$, respectively.}. The rigid body rotation is assumed to be around the material patch centroid in the deformed configuration.
		
		\begin{figure}[hbt]
			\centering
			\includegraphics[width=0.9\textwidth]{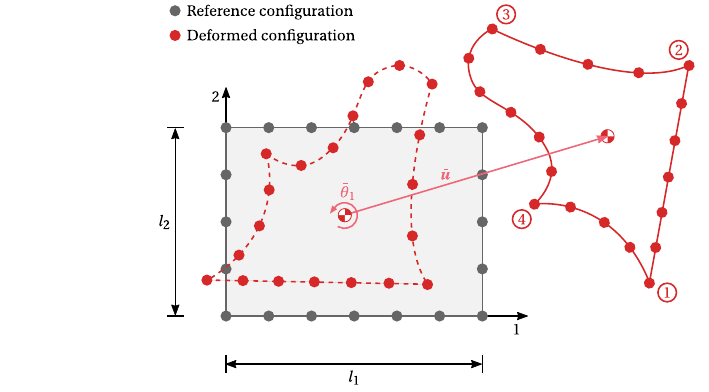}
			\caption{SPDG sampling of a rigid body translation and rotation of the deformed material patch.}
			\label{fig:spdg_step3}
		\end{figure}
		
	\end{itemize}
	
	Finally, when the SPDG-generated material patch is incorporated into a finite element simulation, the Dirichlet boundary conditions prescribed on given boundary nodes can be removed to accommodate any specific multi-axial loading scenario (e.g., free surface).
	
	\begin{remark}
		Given combinations of the sampled geometrical descriptors may naturally yield a non-admissible deformed configuration of the material patch (e.g., intersecting boundary edges), namely when wide bounds are provided. Therefore, a suitable procedure to evaluate the overall deformation admissibility is required to detect non-admissible material patches. Persistent non-admissible deformed configurations suggest the adoption of more conservative bounds.
	\end{remark}

	\section{Material models and parameters \label{sec:material_parameters}}
	This section summarizes the parameters of the several material models employed in the paper to generate the synthetic `ground-truth' data. Three different fully implicit material models are considered: (1) the Von Mises (VM) model with isotropic strain hardening\footnote{The formulation and computational implementation of both Von Mises and Drucker-Prager models employed in this paper are thoroughly discussed by Souza Neto and coworkers \citep{desouzaneto:2008a}.}, (2) the Drucker-Prager (D-P) model with isotropic strain hardening, and (3) the Lou-Zhang-Yoon (LZY) model with isotropic strain hardening described in \ref{sec:lzy_model}.
	
\begin{figure}[hbt]
		\centering
		\begin{subfigure}[b]{0.45\textwidth}
			\centering
			\includegraphics[width=\textwidth]{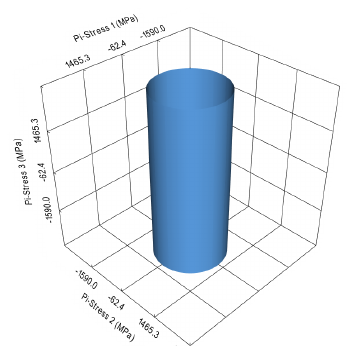}
			\caption{}
			\label{subfig:vm_gt_yield_surface}
		\end{subfigure}
		\begin{subfigure}[b]{0.45\textwidth}
			\centering
			\includegraphics[width=\textwidth]{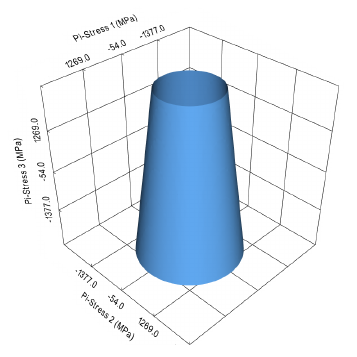}
			\caption{}
			\label{subfig:dp_gt_yield_surface}
		\end{subfigure}\hfill
		\caption{Yield surface: \subref{subfig:vm_gt_yield_surface} Von Mises (VM) `ground-truth' model; \subref{subfig:dp_gt_yield_surface} Drucker-Prager (D-P) `ground-truth' model.}
		\label{fig:vm_dp_yield_surfaces}
	\end{figure}

\begin{figure}[hbt]
		\centering
		\begin{subfigure}[b]{0.45\textwidth}
			\centering
			\includegraphics[width=\textwidth]{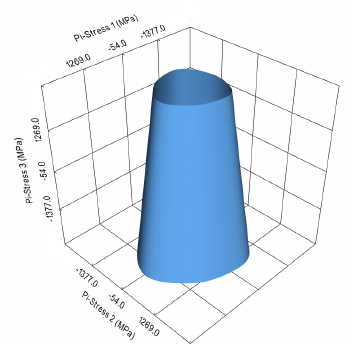}
			\caption{}
			\label{subfig:lou_gt_yield_surface}
		\end{subfigure}
		\begin{subfigure}[b]{0.45\textwidth}
			\centering
			\includegraphics[width=\textwidth]{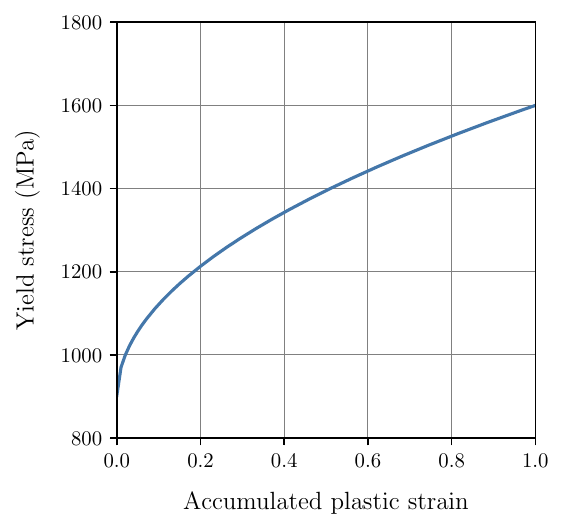}
			\caption{}
			\label{subfig:lou_gt_hardening_law}
		\end{subfigure}\hfill
		\caption{Lou-Zhang-Yoon (LZY) model assumed as the `ground-truth': \subref{subfig:lou_gt_yield_surface} Yield surface; \subref{subfig:lou_gt_hardening_law} Nadai-Ludwik isotropic strain hardening law.}
		\label{fig:lou_local_yield_and_hardening}
	\end{figure}

	The parameters of VM and D-P models are shown in Tables~\ref{tab:vm_gt_material_properties} and \ref{tab:dp_gt_material_properties}, respectively, and the resulting yield surfaces are illustrated in Figure~\ref{fig:vm_dp_yield_surfaces}. The LZY model parameters are shown in Table~\ref{tab:lou_gt_material_properties} and the resulting model yield surface and isotropic hardening law are illustrated in Figure~\ref{fig:lou_local_yield_and_hardening}. The LZY model parameters are chosen such that the `ground-truth' is as complex as possible, namely by including a significant pressure dependency from parameter $b$ (roughly corresponding to a friction angle of $5^{\circ}$) and both curvature and strength differential effect stemming from non-null parameters $c$ and $d$.
	
	\begin{table}[hbt]
		\caption{Material properties of the Von Mises (VM) model assumed as the `ground-truth'.}
		\label{tab:vm_gt_material_properties}
		\centering
		\setlength{\tabcolsep}{0.4cm}
		\renewcommand{\arraystretch}{1.5}
		\begin{tabular}{ccccc}
			\toprule
			\multicolumn{2}{c}{\centering \textbf{Elastic}} & \multicolumn{3}{c}{\centering \textbf{Isotropic hardening}} \\ \cmidrule{1-2} \cmidrule(l){3-5}
			$E$ (GPa) & $\nu$ & $s_{0}$ (MPa) & $s_{1}$ (MPa) & $s_{2}$ \\ \midrule[0.08em]
			110 & 0.33 &  $900\sqrt{3}$ & $700\sqrt{3}$ & 0.5 \\ \bottomrule
			\multicolumn{5}{l}{\small \textbf{Isotropic hardening law:} $\sigma_{y}(\bar{\varepsilon}^{p}) = s_{0} + s_{1} (\bar{\varepsilon}^{p})^{s_{2}}$ (MPa)}
		\end{tabular}
	\end{table}
	
	\begin{table}[hbt]
		\caption{Material properties of the Drucker-Prager (D-P) model assumed as the `ground-truth', where $\phi$ denotes the friction angle, $c$ denotes the cohesion hardening law, $\xi$ is the yield cohesion parameter, $\eta$ is the yield pressure parameter, and $\psi$ is the dilatational angle.}
		\label{tab:dp_gt_material_properties}
		\centering
		\setlength{\tabcolsep}{0.4cm}
		\renewcommand{\arraystretch}{1.5}
		\begin{tabular}{cccccc}
			\toprule
			\multicolumn{2}{c}{\centering \textbf{Elastic}} & \textbf{Yield surface} & \multicolumn{3}{c}{\centering \textbf{Isotropic hardening}} \\ \cmidrule{1-2} \cmidrule(lr){3-3} \cmidrule(l){4-6}
			$E$ (GPa) & $\nu$ & $\phi$ (deg) & $s_{0}$ (MPa) & $s_{1}$ (MPa) & $s_{2}$ \\ \midrule[0.08em]
			110 & 0.33 & 5 &  $900/\xi$ & $700/\xi$ & $0.5$ \\ \bottomrule
			\multicolumn{6}{l}{\small \textbf{Cohesion hardening law:} $c(\bar{\varepsilon}^{p}) = s_{0} + s_{1} (\bar{\varepsilon}^{p})^{s_{2}}$ (MPa)} \\
			\multicolumn{6}{l}{\small \textbf{Additional parameters:} $\xi = (2/\sqrt{3}) \cos(\phi)$, $\eta = (3/\sqrt{3}) \sin(\phi)$, $\psi=\phi$}
		\end{tabular}
	\end{table}
	
	\begin{table}[hbt]
		\caption{Material properties of the Lou-Zhang-Yoon (LZY) model assumed as the `ground-truth'.}
		\label{tab:lou_gt_material_properties}
		\centering
		\setlength{\tabcolsep}{0.4cm}
		\renewcommand{\arraystretch}{1.5}
		\begin{tabular}{ccccccccc}
			\toprule
			\multicolumn{2}{c}{\centering \textbf{Elastic}} & \multicolumn{4}{c}{\centering \textbf{Yield surface}} & \multicolumn{3}{c}{\centering \textbf{Isotropic hardening}} \\ \cmidrule{1-2} \cmidrule(lr){3-6} \cmidrule(l){7-9}
			$E$ (GPa) & $\nu$ & $a$ & $b$ & $c$ & $d$ & $s_{0}$ (MPa) & $s_{1}$ (MPa) & $s_{2}$  \\ \midrule[0.08em]
			110 & 0.33 & 1.0 & 0.05 & 1.0 & 0.5 & 900 & 700 & $0.5$ \\ \bottomrule
			\multicolumn{9}{l}{\small \textbf{Isotropic hardening law:} $\sigma_{y}(\bar{\varepsilon}^{p}) = s_{0} + s_{1} (\bar{\varepsilon}^{p})^{s_{2}}$ (MPa)}
		\end{tabular}
	\end{table}

	\clearpage
	
	\section{Local synthetic data sets \label{sec:local_datasets}}
	This section provides a comprehensive overview of the local synthetic data sets used throughout the paper. All data sets are composed of local synthetic random polynomial strain-stress paths. Taking into account the objectives of this paper, this choice was primarily driven by two key factors: on the one hand, from a model discovery point of view, data diversity is achieved by sampling the whole strain-stress space randomly; on the other hand, from a testing perspective, the model performance is consistently reported on complex strain-stress loading paths rather than a given set of simplified loading conditions.
	
	\subsection{Noiseless random polynomial strain-stress paths}
\begin{figure}[hbt]
		\centering
		\includegraphics[width=0.9\textwidth]{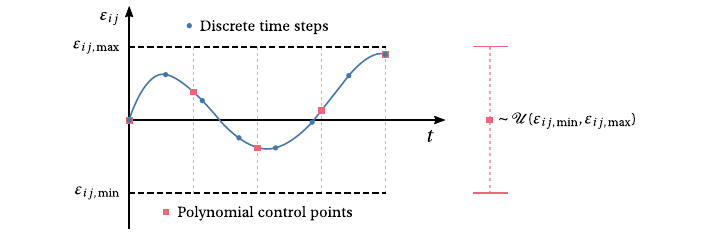}
		\caption{Generation of a random polynomial path of a given strain component $\varepsilon_{ij}$, $\{i, \, j\} \in 1, \, \dots, \, n_\mathrm{dim}$.}
		\label{fig:random_polynomial_path}
	\end{figure}

	The process to generate a random polynomial strain path is akin to the one employed by Mozaffar and coworkers \citep{mozaffar:2019} and is schematically illustrated in Figure~\ref{fig:random_polynomial_path}. The first step consists in setting the total number of time steps, $n_{t}$,
	the polynomial order, $p$, and the sampling bounds for each strain component, $(\varepsilon_{ij, \, \mathrm{min}}, \, \varepsilon_{ij, \, \mathrm{max}})$, $\{i, \, j\} = 1, \, \dots, \, n_\mathrm{dim}$. To ensure more diversity, the polynomial order is here randomly sampled from $(p_{\mathrm{min}}, \, p_{\mathrm{max}})$ for each strain component independently. The second step consists in randomly sampling $p$ control points\footnote{Although a general $p$-th order polynomial is determined from $p+1$ control points, the initial point is always assumed to be zero, corresponding to an initial undeformed configuration.} that are equally spaced on the time axis and randomly sampled from a uniform distribution in the strain axis, i.e., $\varepsilon_{ij} \sim \mathcal{U} (\varepsilon_{ij, \, \mathrm{min}}, \varepsilon_{ij, \, \mathrm{max}})$. In the third step, the $p+1$ control points are used to (deterministically) fit the $p$-order polynomial that describes the strain component $\varepsilon_{ij}$ loading path. At last, the polynomial strain loading path is discretized in $n_{t}$ time steps. After the process is repeated for all strain components, the strain loading path is completely defined.
	
	To complete the local data set sample, the strain loading path is fed to a given material model to compute the stress response path, $(\sigma_{ij})_{t}$, $\{i, \, j\} = 1, \, \dots, \, n_\mathrm{dim}$, and $t = 0, \, \dots, \, n_{t} - 1$, by solving the corresponding state update algorithm\footnote{Despite not being employed in the model discovery process carried out in the present paper, any other relevant material model outputs, such as internal state variables, may also be stored as data set features.}. The whole procedure can be repeated $N$ times to generate a local synthetic data set of $N$ random polynomial strain-stress paths.
	
\begin{figure}[hbt]
		\centering
		\begin{subfigure}[b]{0.495\textwidth}
			\centering
			\includegraphics[width=\textwidth]{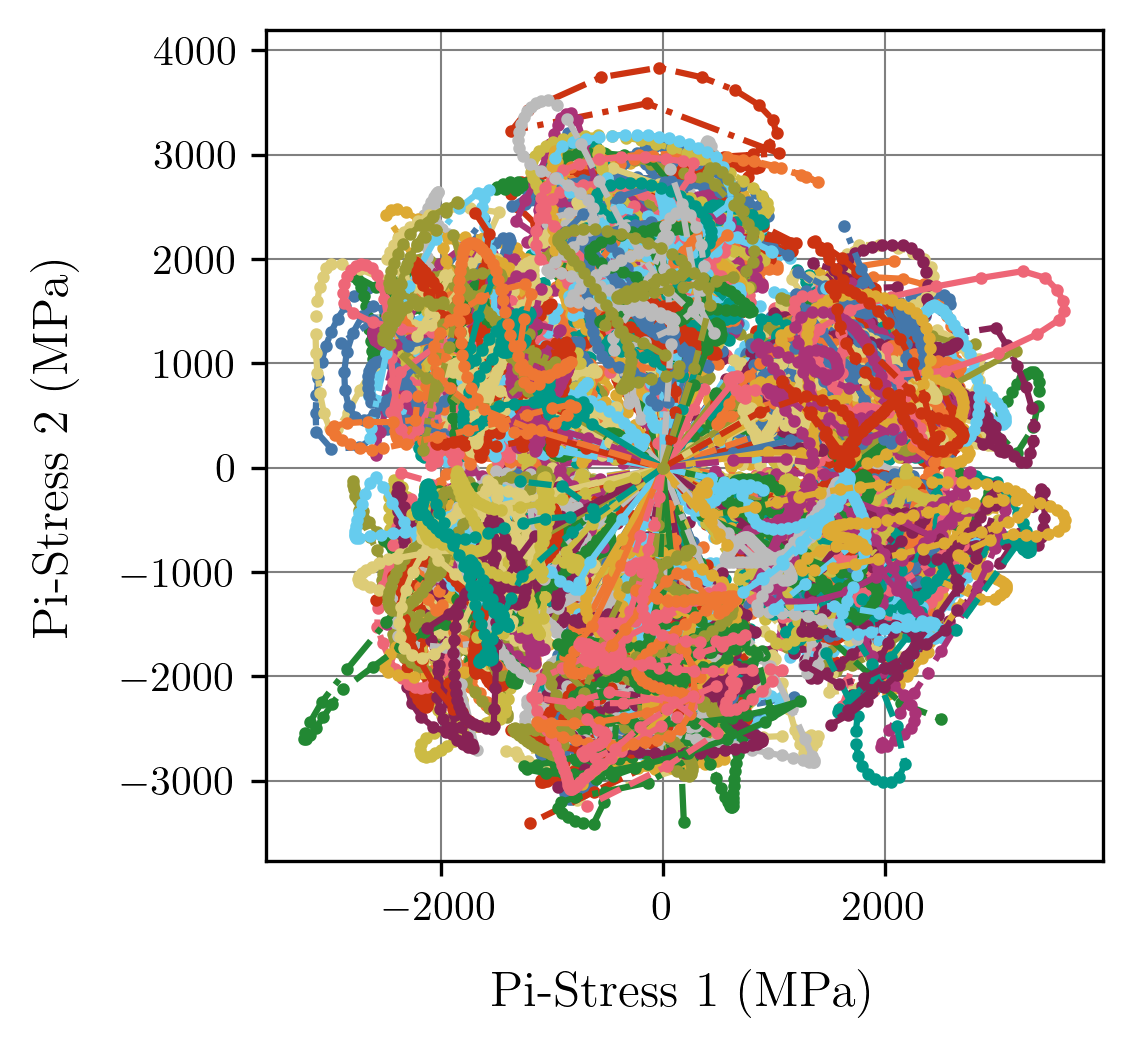}
			\caption{}
			\label{subfig:lou_dataset_pi1_pi2}
		\end{subfigure}
		\begin{subfigure}[b]{0.495\textwidth}
			\centering
			\includegraphics[width=\textwidth]{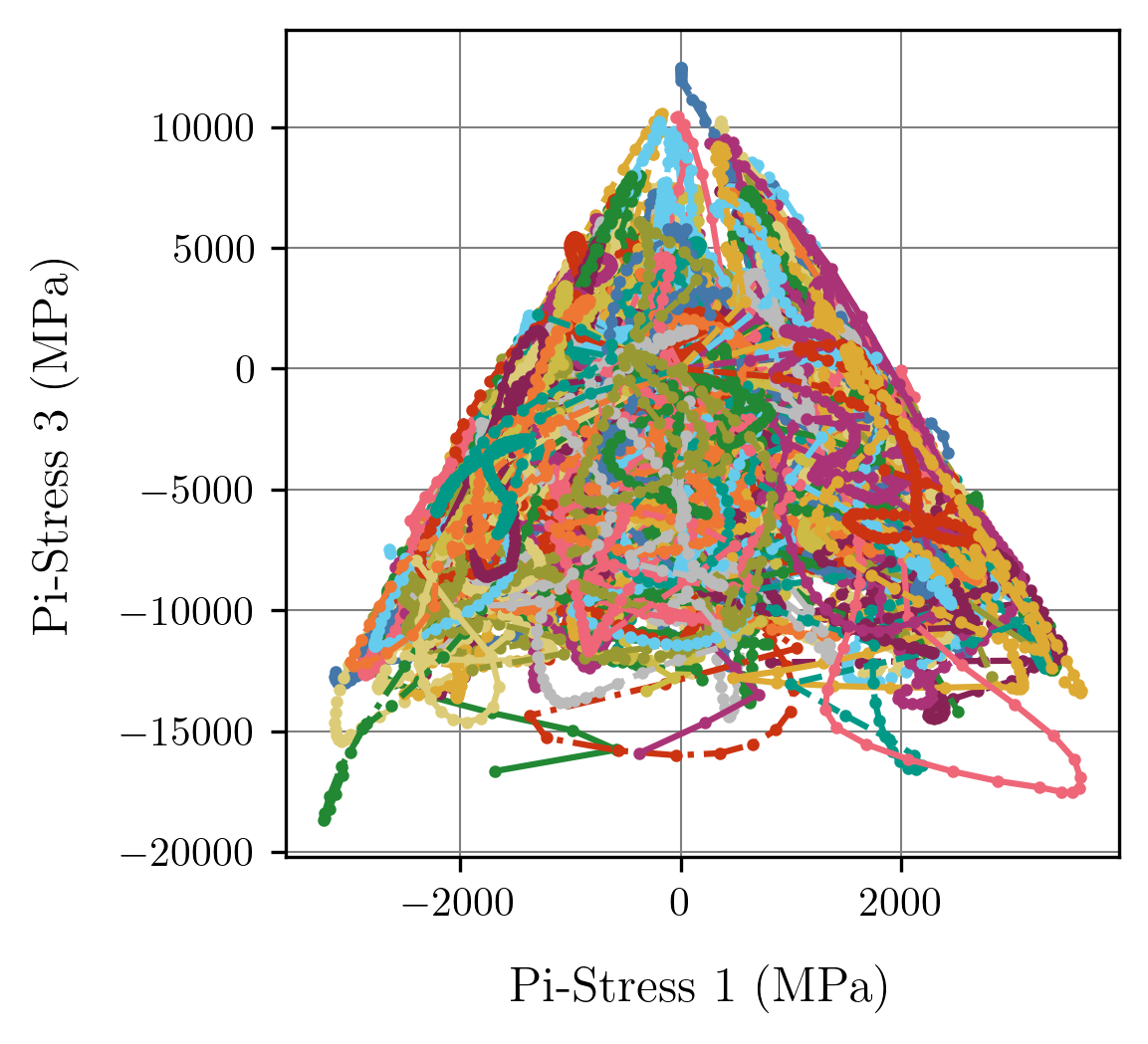}
			\caption{}
			\label{subfig:lou_dataset_pi1_pi3}
		\end{subfigure}\hfill
		\begin{subfigure}[b]{0.495\textwidth}
			\centering
			\includegraphics[width=\textwidth]{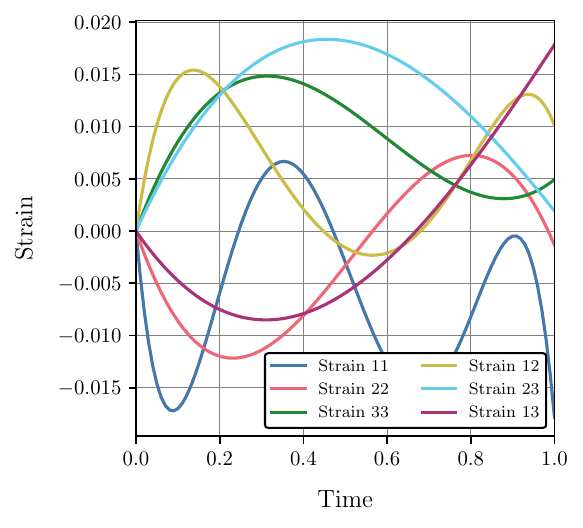}
			\caption{}
			\label{subfig:lou_dataset_strain_path}
		\end{subfigure}
		\begin{subfigure}[b]{0.495\textwidth}
			\centering
			\includegraphics[width=\textwidth]{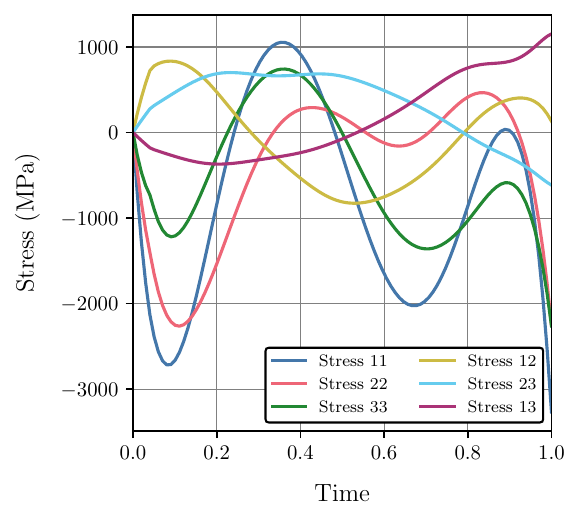}
			\caption{}
			\label{subfig:lou_dataset_stress_path}
		\end{subfigure}\hfill
		\caption{Example of local synthetic data set of 320 random polynomial strain-stress paths generated with LZY model ($E=110$GPa, $\nu=0.33$, $a=1.0$, $b=0.05$, $c=1.0$, $d=0.5$, $s_{0}=900\,$MPa, $s_{1}=700\,$MPa, $s_{2}=0.5$):
			\subref{subfig:lou_dataset_pi1_pi2} $\pi_{1}-\pi_{2}$ stress projection view; \subref{subfig:lou_dataset_pi1_pi3} $\pi_{1}-\pi_{3}$ stress projection view; \subref{subfig:lou_dataset_strain_path} Random sample strain path; \subref{subfig:lou_dataset_stress_path} Random sample stress path.}
		\label{fig:lou_local_dataset_320}
	\end{figure}

	In the present paper, the strain range defined by $(\varepsilon_{ij, \, \mathrm{min}}, \, \varepsilon_{ij, \, \mathrm{max}}) = (-0.02, \, 0.02)$, $\{i, \, j\} = 1, \, \dots, \, n_\mathrm{dim}$, is adopted for all local synthetic data sets, which is compatible with the infinitesimal strains assumption. Moreover, taking advantage of the implicit nature of the different conventional models, the strain-stress paths are discretized in $n_{t}=100$ time steps. An example of a local synthetic data set of 320 random polynomial strain-stress paths generated with the Lou-Zhang-Yoon (LZY) model  is shown in Figure~\ref{fig:lou_local_dataset_320}.
	
	\begin{remark}
		Given the complex nature of the randomly generated polynomial strain paths, involving rich combinations of strain components as well as intricate loading/unloading patterns, it is of utmost importance that the material model state update algorithm is robust. Moreover, the computational cost of such a material model is an important factor to take into account when aiming to generate large data sets.
	\end{remark}
	
	\subsection{Noisy random polynomial strain paths}
	To evaluate the neural network model's performance when discovered from experimental data, synthetic noise is introduced into the local data sets to emulate varying degrees of experimental noise. In the experimental setting, it is assumed that the noise stems essentially from the displacement field measurement (e.g., with Digital Image Correlation (DIC) or Digital Volume Correlation (DVC)), from which the strain field can be computed by means of a suitable numerical treatment. For this reason, synthetic noise is here injected in the local random polynomial strain paths previously described.
	
	Three different noise distribution types (Gaussian, Uniform, and Spiked Gaussian) and four different noise levels (low, medium, high, and very high) are considered as follows:
	\begin{itemize}
		\item \textbf{Gaussian noise}. Gaussian distribution with zero mean and given standard deviation,
		\begin{equation}
			\tilde{\varepsilon} \sim \mathcal{N}(0, \, \epsilon_{\mathrm{std}}) \, ,
			\label{eq:gaussian_noise}
		\end{equation}
		where different noise levels are achieved by varying noise standard deviation;
		\item \textbf{Uniform noise}. Uniform distribution zero-centered and with given amplitude,
		\begin{equation}
			\tilde{\varepsilon} \sim \mathcal{U}(-0.5 \, \epsilon_{\mathrm{amp}}, \, 0.5 \, \epsilon_{\mathrm{amp}}) \, ,
			\label{eq:uniform_noise}
		\end{equation}
		where different noise levels are achieved by varying noise amplitude;
		\item \textbf{Spiked Gaussian noise}. Spiked Gaussian distribution resulting from superimposing a scaled Bernoulli distribution, with a given probability, onto a Gaussian distribution, with zero mean and given standard deviation, 
		\begin{equation}
			\tilde{\varepsilon} \sim \mathcal{N}(0, \, \epsilon_{\mathrm{std}}) + \epsilon_{\mathrm{spk}} \, \mathcal{B}(p) \, ,
			\label{eq:spiked_gaussian_noise}
		\end{equation}
		where different noise levels are achieved by varying the Gaussian noise standard deviation, keeping the noise spike distribution fixed.
	\end{itemize}
	
	In addition, both homoscedastic and heteroscedastic noise variability cases are taken into account. While the noise does not depend on the specific measurement in the homoscedastic case, an additional strain-dependent noise distribution component is considered in the heteroscedastic case. Assume that, at each time step $t$, a scaling parameter defined as $\alpha (t) = \norm{\bm{\varepsilon}_{t}}/\norm{\bm{\varepsilon}}_{\mathrm{max}}$, where the normalizing factor is taken as $\norm{\bm{\varepsilon}}_{\mathrm{max}} = 0.02$, in line with the strain bounds. The heteroscedastic counterparts of the aforementioned homoscedastic noise distributions are described as follows: (i) in both Gaussian and Spiked Gaussian cases, a distribution $\mathcal{N}(0, \, \alpha (t) \, \epsilon_{\mathrm{std}})$ is added to Equations~\eqref{eq:gaussian_noise} and \eqref{eq:spiked_gaussian_noise}, respectively, while (ii) in the Uniform case, a distribution $\mathcal{U}(-0.5 \, \alpha (t) \, \epsilon_{\mathrm{amp}}, \, 0.5 \, \alpha (t) \, \epsilon_{\mathrm{amp}})$ is added to Equation~\eqref{eq:uniform_noise}.

	\clearpage
	
	\begin{remark}
		The noise is sampled independently for each strain path, strain component, and time step. In addition, the noise is conveniently defined and superimposed onto the noiseless strain path in the normalized space, where a min-max scaling based on the strain bounds, $(\varepsilon_{ij, \, \mathrm{min}}, \, \varepsilon_{ij, \, \mathrm{max}})$, $\{i, \, j\} = 1, \, \dots, \, n_\mathrm{dim}$, is employed. The resulting noisy strain path is then denormalized to the original strain space.
	\end{remark}
	
\begin{figure}[hbt]
		\centering
		\begin{subfigure}[b]{0.495\textwidth}
			\centering
			\includegraphics[width=\textwidth]{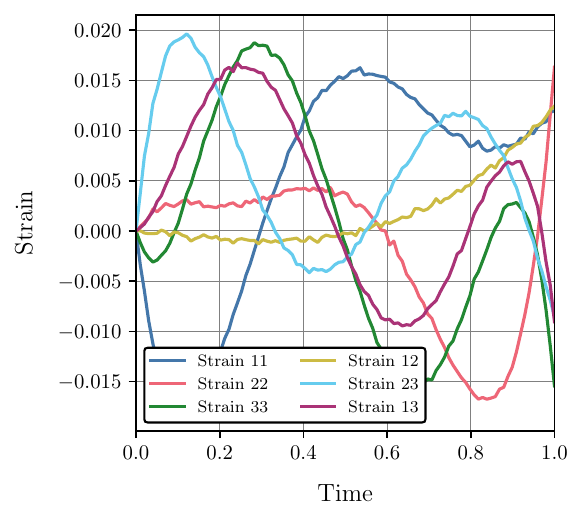}
			\caption{}
			\label{subfig:noisy_strain_path_1}
		\end{subfigure}
		\begin{subfigure}[b]{0.495\textwidth}
			\centering
			\includegraphics[width=\textwidth]{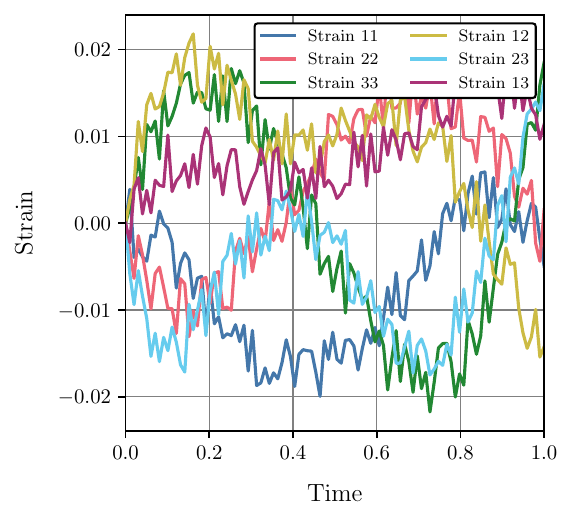}
			\caption{}
			\label{subfig:noisy_strain_path_2}
		\end{subfigure}\hfill
		\begin{subfigure}[b]{0.495\textwidth}
			\centering
			\includegraphics[width=\textwidth]{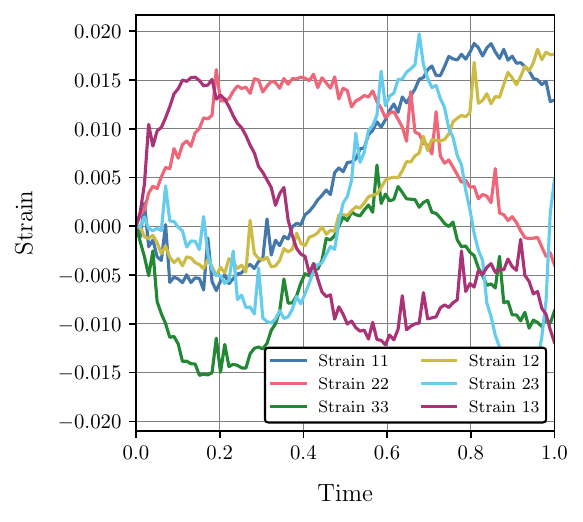}
			\caption{}
			\label{subfig:noisy_strain_path_3}
		\end{subfigure}
		\begin{subfigure}[b]{0.495\textwidth}
			\centering
			\includegraphics[width=\textwidth]{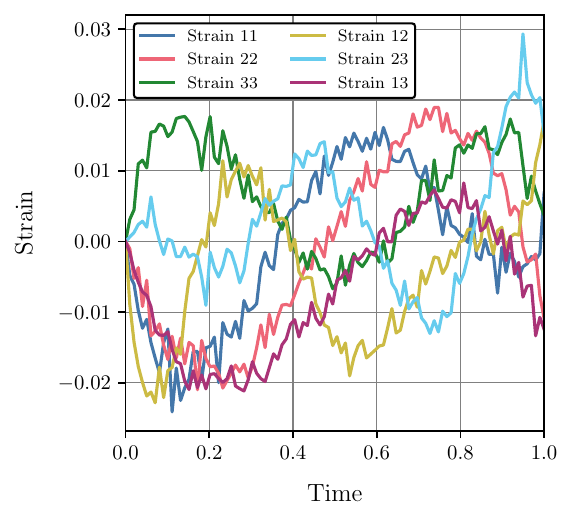}
			\caption{}
			\label{subfig:noisy_strain_path_4}
		\end{subfigure}\hfill
		\caption{Examples of noisy local synthetic random polynomial strain paths:
			\subref{subfig:noisy_strain_path_1} Low level Gaussian homoscedastic noise; \subref{subfig:noisy_strain_path_2} Very high level Uniform homoscedastic noise; \subref{subfig:noisy_strain_path_3} Medium level Spiked Gaussian homoscedastic noise; \subref{subfig:noisy_strain_path_4} High level Uniform heteroscedastic noise.}
		\label{fig:noisy_strain_paths}
	\end{figure}

	The different noise levels considered in this paper are summarized in Table~\ref{tab:noise_dist_params} for each noise distribution type. Examples of noisy random polynomial strain paths are shown in Figure~\ref{fig:noisy_strain_paths}.
	
	\begin{table}[hbt]
		\caption{Parameters of the different noise distribution types and corresponding noise levels employed to generate noisy random polynomial strain paths. The parameter $\epsilon_{\mathrm{std}}$ corresponds to the standard deviation of the Gaussian distribution, while $\epsilon_{\mathrm{amp}}$ corresponds to the amplitude of the zero-centered Uniform distribution. The Uniform distribution amplitude spans two standard deviations of the same noise level Gaussian distribution. A spike of magnitude $\epsilon_{\mathrm{spk}}=0.2$ and probability $p=0.05$ is set for the Spiked Gaussian distribution.}
		\label{tab:noise_dist_params}
		\centering
		\setlength{\tabcolsep}{0.4cm}
		\renewcommand{\arraystretch}{1.0}
		\begin{tabular}{cccccc}
			\toprule
			\multirow{2}{4.0cm}{\centering \textbf{Noise distribution}} & \multirow{2}{2.5cm}{\centering \textbf{Parameter}} & \multicolumn{4}{c}{\textbf{Noise level}} \\ \cmidrule{3-6}
			& & Low & Medium & High & Very high \\ \midrule[0.08em]
			Gaussian & $\epsilon_{\mathrm{std}}$ & 0.01 & 0.025 & 0.05 & 0.1 \\ \midrule
			Uniform & $\epsilon_{\mathrm{amp}}$ & 0.04 & 0.1 & 0.2 & 0.4 \\ \midrule
			Spiked Gaussian & $\epsilon_{\mathrm{std}}$ & 0.01 & 0.025 & 0.05 & 0.1 \\ \bottomrule 
		\end{tabular}
	\end{table}
	
	\subsection{Training, validation and testing}
	Given that several convergence analyses with respect to the training data set size are carried out throughout this paper, local data sets of different sizes are required. In what concerns neural network and hybrid models in a local model discovery setting, the data set sizes considered in this paper are shown in Table~\ref{tab:local_gru_dataset_sizes}. The training data set size is a geometric sequence with a common ratio of 2, while the validation data set size is established as 20\% of the training size. It should be remarked that the validation data set is only used in the early stopping criterion adopted in the training procedure. In turn, the testing data set size is 20\% of the largest training size and is consistently shared across all data sets in the convergence analysis.
	
	\begin{table}[hbt]
		\caption{Size of local synthetic data sets of random polynomial strain-stress paths used in the local discovery and testing of neural network and hybrid material models. The testing data set is consistently shared across all data sets to evaluate the material model performance.}
		\label{tab:local_gru_dataset_sizes}
		\centering
		\setlength{\tabcolsep}{0.3cm}
		\renewcommand{\arraystretch}{1.5}
		\begin{tabular}{cccccccccc}
			\toprule
			\textbf{Data set \#} & \textbf{1} & \textbf{2} & \textbf{ 3} & \textbf{ 4} & \textbf{ 5} & \textbf{ 6} & \textbf{ 7} & \textbf{ 8} & \textbf{ 9} \\ \midrule
			Training & 10 & 20 & 40 & 80 & 160 & 320 & 640 & 1280 & 2560 \\
			Validation & 2 & 4 & 8 & 16 & 32 & 64 & 128 & 256 & 512 \\
			Testing & \multicolumn{9}{c}{512} \\ \bottomrule \\[-5pt]	
		\end{tabular}
	\end{table}

	In the local model discovery of conventional models, involving a small number of parameters, only a single data set is selected for convergence demonstrative purposes. Given that the goal is to recover the known `ground-truth' parameters, i.e., there is no overfitting in this context, validation is not required during the discovery procedure.
	
	\clearpage

	\section{Additional results of local model discovery with noisy data \label{sec:additional_noise_local}}
	
	This section includes additional results concerning the model discovery from noisy data and the corresponding performance (see Section~\ref{ssec:local_ml_models}).
	
	\begin{figure}[hbt]
		\centering
		\begin{subfigure}[b]{0.43\textwidth}
			\centering
			\includegraphics[width=\textwidth]{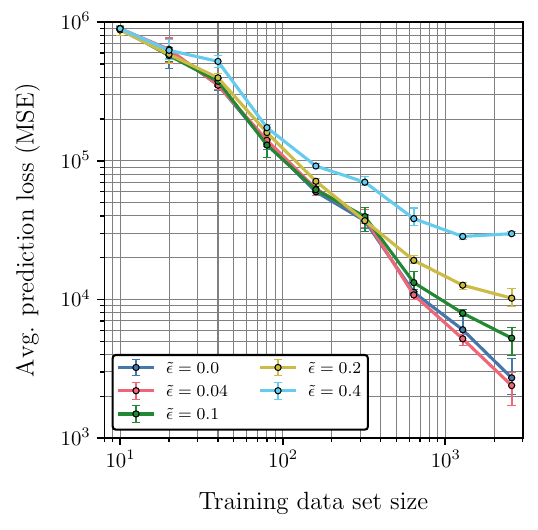}
			\caption{}
			\label{subfig:gru_conv_analysis_avg_mse_homuni_noise}
		\end{subfigure}
		\begin{subfigure}[b]{0.43\textwidth}
			\centering
			\includegraphics[width=\textwidth]{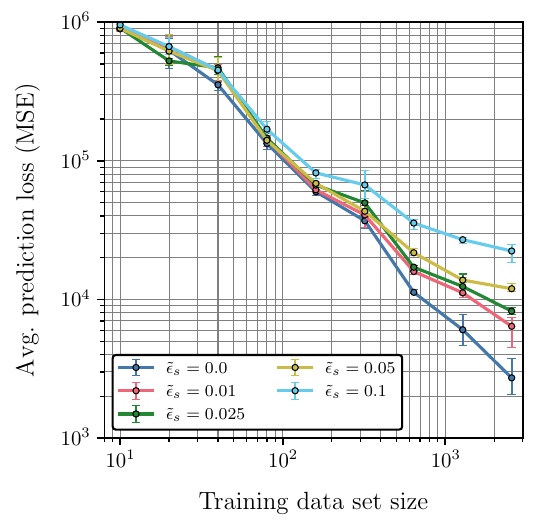}
			\caption{}
			\label{subfig:gru_conv_analysis_avg_mse_homsgau_noise}
		\end{subfigure}\hfill
		\begin{subfigure}[b]{0.43\textwidth}
			\centering
			\includegraphics[width=\textwidth]{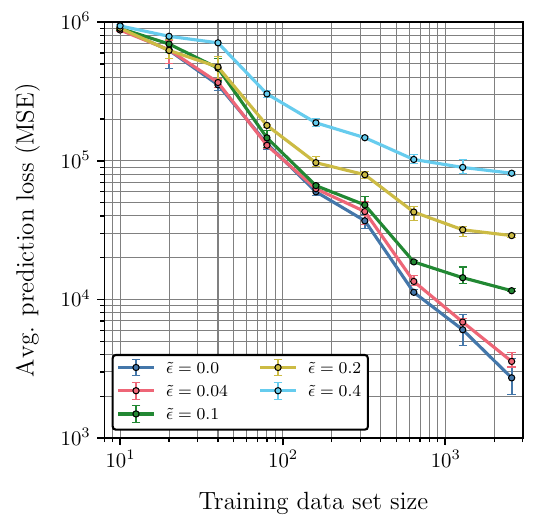}
			\caption{}
			\label{subfig:gru_conv_analysis_avg_mse_hetuni_noise}
		\end{subfigure}
		\begin{subfigure}[b]{0.43\textwidth}
			\centering
			\includegraphics[width=\textwidth]{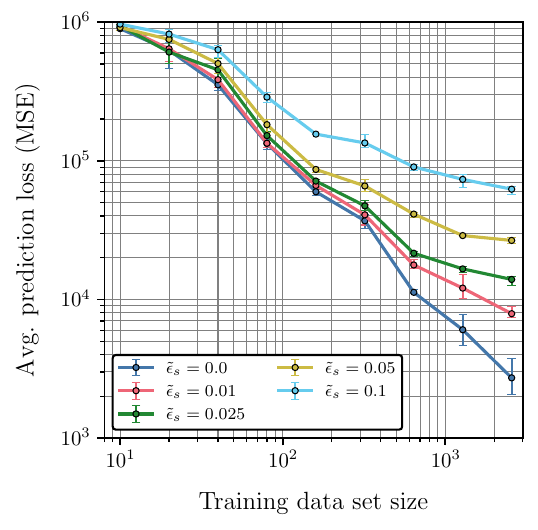}
			\caption{}
			\label{subfig:gru_conv_analysis_avg_mse_hetsgau_noise}
		\end{subfigure}\hfill
		\caption{Performance of the GRU material model discovered from LZY noiseless and noisy strain-stress data ($E=110$GPa, $\nu=0.33$, $a=1.0$, $b=0.05$, $c=1.0$, $d=0.5$, $s_{0}=900\,$MPa, $s_{1}=700\,$MPa, $s_{2}=0.5$) and tested in a (unseen) local synthetic data set of 512 random polynomial strain-stress paths. Average prediction loss (MSE) with respect to the training data set size and uncertainty quantification (3 random model initializations) for different noise distributions:
			\subref{subfig:gru_conv_analysis_avg_mse_homuni_noise} Uniform homoscedastic noise; \subref{subfig:gru_conv_analysis_avg_mse_homsgau_noise} Spiked Gaussian homoscedastic noise; \subref{subfig:gru_conv_analysis_avg_mse_hetuni_noise} Uniform heteroscedastic noise; \subref{subfig:gru_conv_analysis_avg_mse_hetsgau_noise} Spiked Gaussian heteroscedastic noise.}
		\label{fig:gru_conv_analysis_avg_mse_hom_het_noise}
	\end{figure}
	
	\clearpage
	
	\section{Global synthetic data sets}
	This section provides a comprehensive description of the synthetic data sets used in this paper to perform a global indirect model discovery. Each global synthetic data set corresponds to a unique specimen geometry, subjected to specific loading conditions, and is computed assuming a particular material model with a given set of properties. Once the specimen finite element simulation is solved, the finite element mesh nodes displacements and reaction forces are collected to build the corresponding data set.
	
	\begin{remark}
		The finite element simulations required to generate each specimen `ground-truth' data were carried out with LINKS (Large Strain Implicit Non-linear Analysis of Solids Linking Scales), an implicit multi-scale finite element Fortran private code developed by the CM2S research group at the Faculty of Engineering of the University of Porto, Portugal. This choice stems solely from the fact that the first author is a former member of the CM2S research group and developer of LINKS. Any other open-source or commercial finite element software can be alternatively used to perform these simulations.
	\end{remark}
	
	\subsection{Tensile dogbone specimen \label{ssec:tensile_dogbone}}
\begin{figure}[hbt]
		\centering
		\begin{subfigure}[b]{0.45\textwidth}
			\centering
			\includegraphics[width=\textwidth]{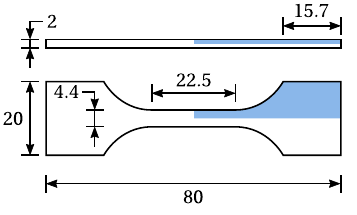}
			\caption{}
			\label{subfig:dogbone_geometry}
		\end{subfigure}
		\begin{subfigure}[b]{0.45\textwidth}
			\centering
			\includegraphics[width=\textwidth]{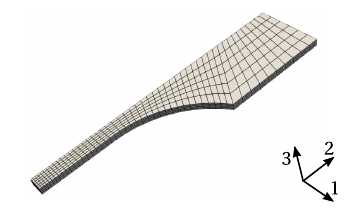}
			\caption{}
			\label{subfig:dogbone_mesh}
		\end{subfigure}\vspace*{10pt}  \hfill
		\begin{subfigure}[b]{0.45\textwidth}
			\centering
			\includegraphics[width=\textwidth]{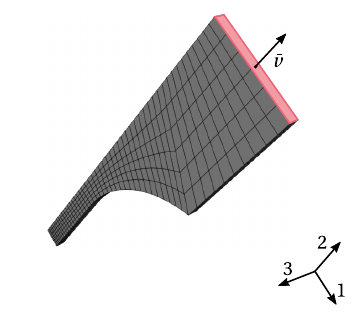}
			\caption{}
			\label{subfig:dogbone_bc}
		\end{subfigure}
		\begin{subfigure}[b]{0.45\textwidth}
			\centering
			\includegraphics[width=\textwidth]{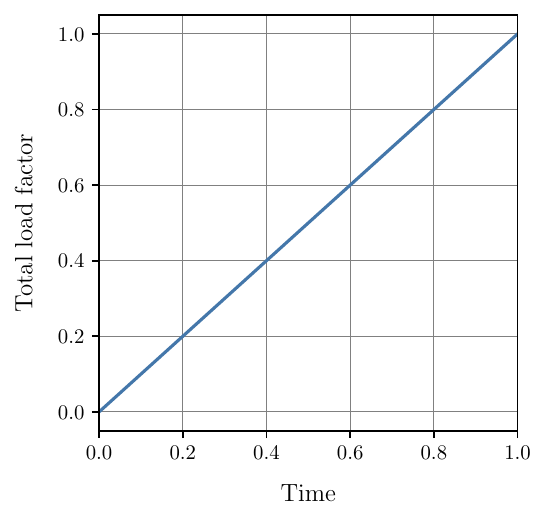}
			\caption{}
			\label{subfig:dogbone_loading_factor}
		\end{subfigure}\hfill
		\caption{Tensile dogbone specimen: \subref{subfig:dogbone_geometry} Geometry; \subref{subfig:dogbone_mesh} Finite element mesh (1200 8-node hexahedral finite elements); \subref{subfig:dogbone_bc} Prescribed uniaxial displacement, $\bar{v}$ ; \subref{subfig:dogbone_loading_factor} Total load factor history. Dimensions are displayed in millimiters.}
		\label{fig:dogbone_specimen}
	\end{figure}

	A tensile dogbone specimen is employed in this paper to illustrate a global indirect discovery of both von Mises and Drucker-Prager models. The specimen geometry is shown in Figure~\ref{subfig:dogbone_geometry} together with the corresponding 8-node hexahedral finite element mesh\footnote{In a purely computational setting, a 20-node hexahedral finite element mesh (or any other finite element type) could be employed instead. However, given that the displacement field experimental measurements are assumed to stem from Digital Image Correlation (DIC) or Digital Volume Correlation (DVC), a 8-node hexahedral finite element discretization is chosen to allow a suitable encoding of the experimental measurements into the finite element mesh.} discretization that takes into account symmetry conditions (see Figure~\ref{subfig:dogbone_mesh}).
	
	To emulate the experimental scenario under a uniaxial tensile testing machine, the loading is applied by means of a prescribed uniaxial displacement, $\bar{v}$, at the specimen top (see Figure~\ref{subfig:dogbone_bc}) and appropriate boundary conditions are applied in the symmetry planes. In addition, a proportional loading scheme is employed such that $v_{t}=\lambda(t) \, \bar{v}$, $t=0, \, \dots, \, n_{t} - 1$, where $\lambda$ denotes the total load factor (see Figure~\ref{subfig:dogbone_loading_factor}). In this context, a total displacement of $\bar{v}=10 \,$mm is prescribed in a total of 200 increments of equal magnitude.
	
	\subsection{Tensile double notched specimen \label{ssec:tensile_double_notched}}
	
\begin{figure}[hbt]
		\centering
		\begin{subfigure}[b]{0.45\textwidth}
			\centering
			\includegraphics[width=\textwidth]{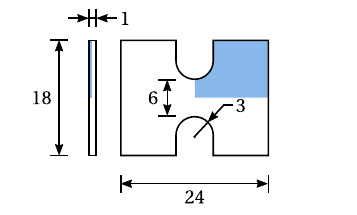}
			\caption{}
			\label{subfig:double_notched_geometry}
		\end{subfigure}
		\begin{subfigure}[b]{0.45\textwidth}
			\centering
			\includegraphics[width=\textwidth]{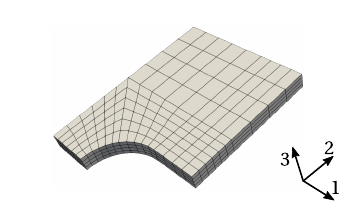}
			\caption{}
			\label{subfig:double_notched_mesh}
		\end{subfigure}\vspace*{10pt}  \hfill
		\begin{subfigure}[b]{0.45\textwidth}
			\centering
			\includegraphics[width=\textwidth]{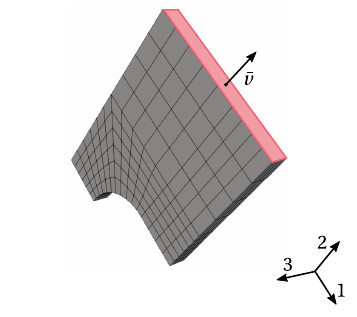}
			\caption{}
			\label{subfig:double_notched_bc}
		\end{subfigure}
		\begin{subfigure}[b]{0.45\textwidth}
			\centering
			\includegraphics[width=\textwidth]{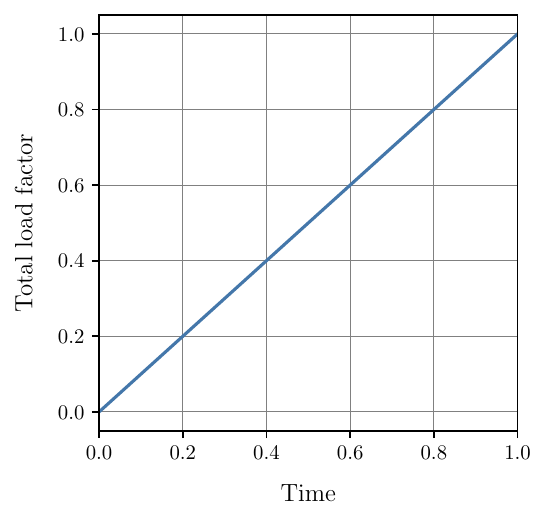}
			\caption{}
			\label{subfig:double_notched_loading_factor}
		\end{subfigure}\hfill
		\caption{Tensile double notched specimen: \subref{subfig:double_notched_geometry} Geometry; \subref{subfig:double_notched_mesh} Finite element mesh (504 8-node hexahedral finite elements); \subref{subfig:double_notched_bc} Prescribed uniaxial displacement, $\bar{v}$ ; \subref{subfig:double_notched_loading_factor} Total load factor history. Dimensions are displayed in millimiters.}
		\label{fig:double_notched_specimen}
	\end{figure}

	A tensile double notched specimen is considered in this paper to perform the global indirect discovery of the Lou-Zhang-Yoon model. The specimen geometry is shown in Figure~\ref{subfig:double_notched_geometry} together with the corresponding 8-node hexahedral finite element mesh discretization that takes into account symmetry conditions (see Figure~\ref{subfig:double_notched_mesh}).
	
	The loading is applied by means of a prescribed uniaxial displacement, $\bar{v}$, at the specimen top (see Figure~\ref{subfig:double_notched_bc}), emulating once again the experimental scenario under a uniaxial tensile testing machine. Appropriate boundary conditions are applied in the symmetry planes and the same proportional loading scheme is adopted (see Figure~\ref{subfig:double_notched_loading_factor}), $v_{t}=\lambda(t) \, \bar{v}$, $t=0, \, \dots, \, n_{t} - 1$. In this context, a total displacement of $\bar{v}=1 \,$mm is prescribed in a total of 200 increments of equal magnitude.
	
	\clearpage
	
	\subsection{Randomly deformed material patch specimen \label{ssec:random_patch}}
\begin{figure}[hbt]
		\centering
		\begin{subfigure}[b]{0.45\textwidth}
			\centering
			\includegraphics[width=\textwidth]{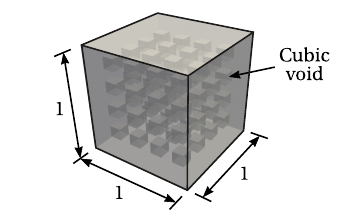}
			\caption{}
			\label{subfig:random_specimen_geometry}
		\end{subfigure}
		\begin{subfigure}[b]{0.45\textwidth}
			\centering
			\includegraphics[width=\textwidth]{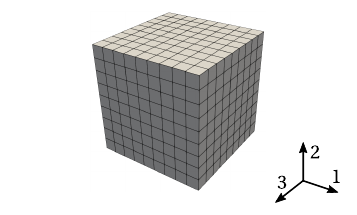}
			\caption{}
			\label{subfig:random_specimen_mesh}
		\end{subfigure}\vspace*{10pt}  \hfill
		\begin{subfigure}[b]{0.45\textwidth}
			\centering
			\includegraphics[width=\textwidth]{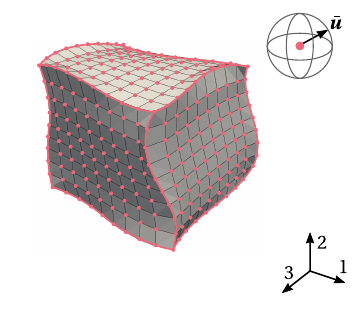}
			\caption{}
			\label{subfig:random_specimen_bc}
		\end{subfigure}
		\begin{subfigure}[b]{0.45\textwidth}
			\centering
			\includegraphics[width=\textwidth]{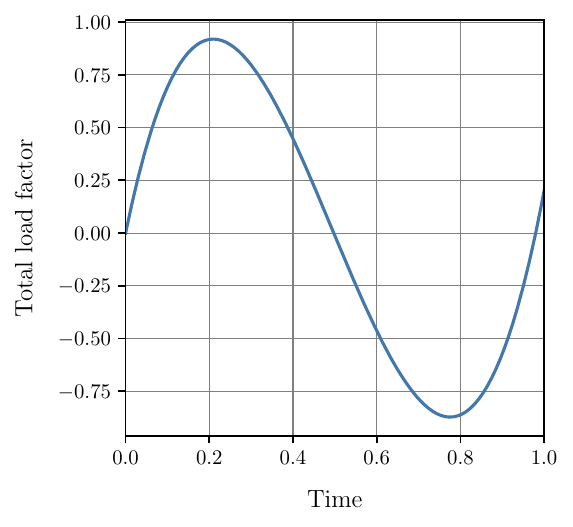}
			\caption{}
			\label{subfig:random_specimen_loading_factor}
		\end{subfigure}\hfill
		\caption{Randomly deformed material patch specimen: \subref{subfig:random_specimen_geometry} Geometry; \subref{subfig:random_specimen_mesh} Finite element mesh (665 8-node hexahedral finite elements); \subref{subfig:random_specimen_bc} Deformed configuration resulting from prescribed boundary displacements for $\lambda=0.92$; \subref{subfig:random_specimen_loading_factor} Total load factor history. Dimensions are displayed in millimiters.}
		\label{fig:random_specimen}
	\end{figure}

	An idealized randomly deformed material patch specimen is generated with SPDG (\ref{sec:spdg}) to showcase the global indirect discovery of a GRU material model from LZY model `ground-truth' data. By also including an internal regular grid of cubic voids, the purpose of such a specimen is to induce a broad converage of the strain-stress space and provide a rich data set to the physics-agnostic GRU material model. The specimen geometry is shown in Figure~\ref{subfig:random_specimen_geometry} together with the corresponding 8-node hexahedral finite element mesh discretization (see Figure~\ref{subfig:random_specimen_mesh}).
	
	The loading is applied by means of prescribed displacements at the specimen boundaries as described in \ref{sec:spdg}, leading to the deformed configuration illustrated in Figure~\ref{subfig:random_specimen_bc}. To promote an even higher strain field diversity, the random polynomial non-monotonic proportional loading scheme shown in Figure~\ref{subfig:random_specimen_loading_factor} is adopted. Given the prescribed displacement at a given boundary node, $\bar{\bm{u}}$, the corresponding displacement at time $t$, is given by $\bm{u}_{t}=\lambda(t) \, \bar{\bm{u}}$, $t=0, \, \dots, \, n_{t} - 1$. Note that a negative total load factor entails a similar displacement magnitude but applied in the opposite direction. The loading path is prescribed in a total of 200 increments.
	
	\section{Data set pruning algorithm \label{sec:pruning_algorithm}}
	This section provides a brief description of the data set pruning algorithm recently proposed by Li and coworkers \citep{li:2023a}. In this paper, the algorithm is employed to evaluate redundancy in the local strain-stress paths data set derived from a given specimen, and to explore its impact on the performance of a neural network material model.
	
	Given a data set of strain-stress paths, $D_{\mathrm{full}}$, the pruning algorithm operates by iteratively reducing the data set through the removal of paths exhibiting higher redundancy. The iterative process is shown in Table~\ref{tab:dataset_pruning_iterations}. Two types of iterations are distinguished here: (i) a model iteration, $k$, where the model $\mathcal{M}^{\, (k)}$, is trained\footnote{The data set $D_{t+v}$ denotes the joint training and validation data set, where the validation data set is solely used in the early stopping criterion during the training process.} on a data set $D_{t+v}^{\,(k)}$, and (ii) a pruning iteration, $k$, where the data set $D_{t+v}^{\,(k)}$ is obtained by pruning $D_{t+v}^{\, (k-1)}$. The algorithm starts by setting $D_{t+v}^{\,0} = D_{\mathrm{full}}$ and initializing the unused data set, $D_{\mathrm{unused}}^{\,0} = \emptyset$. Then, at each pruning iteration, $k$, $D_{t+v}^{\, (k-1)}$ is randomly split into training, validation and testing data sets, denoted as $D_{A}^{\, (k-1)}$, $D_{B}^{\, (k-1)}$ and $D_{C}^{\, (k-1)}$, respectively, and the pruning model, $\mathcal{M}^{*}$, is trained and tested on $D_{C}^{\, (k-1)}$. The testing paths with the lowest prediction error, $D_{C, \, \mathrm{low}}^{\, (k-1)}$, are pruned from the data set, $D_{t+v}^{\, (k)} = D_{t+v}^{\, (k-1)} - D_{C, \, \mathrm{low}}^{\, (k-1)}$, and added to the unused data set, $D_{\mathrm{unused}}^{\, (k)} = D_{\mathrm{unused}}^{\, (k-1)} + D_{C, \, \mathrm{low}}^{\, (k-1)}$. The pruning process terminates when given criteria are met, e.g., a minimum data set size is reached. The performance of the model $\mathcal{M}^{\, (k)}$ can be optionally evaluated throughout the process by training on $D_{t+v}^{\,(k)}$ and testing on any particular data set of interest (e.g., in-distribution, out-of-distribution, unused).
	
	\begin{table}[hbt]
		\caption{Data set pruning algorithm. The data set $D_{\mathrm{full}}$ is iteratively reduced by removing paths exhibiting higher redundancy based on the testing performance of the model $\mathcal{M}^{*}$ in each iteration.}
		\label{tab:dataset_pruning_iterations}
		\centering
		\setlength{\tabcolsep}{0.3cm}
		\renewcommand{\arraystretch}{1.5}
		\begin{tabular}{ccccc}
			\toprule
			\textbf{Iteration} & \textbf{Training} & \textbf{Model} & \textbf{Testing} & \textbf{Unused} \\ \midrule
			0 & $D_{t+v}^{\,0} = D_{\mathrm{full}}$ & $\mathcal{M}^{\, (0)}$  & Any & $D_{\mathrm{unused}}^{\,0} = \emptyset$ \\ \midrule
			\multirow{2}{*}{1} & $D_{A}^{\, (0)}$, $D_{B}^{\, (0)}$ & $\mathcal{M}^{*}$ & $D_{C}^{\, (0)}$ & $D_{\mathrm{unused}}^{\, (1)} = D_{\mathrm{unused}}^{\, (0)} + D_{C, \, \mathrm{low}}^{\, (0)}$ \\
			& $D_{t+v}^{\, (1)} = D_{t+v}^{\, (0)} - D_{C, \, \mathrm{low}}^{\, (0)}$ & $\mathcal{M}^{\, (1)}$ & Any & - \\ \midrule
			$\vdots$ & $\vdots$ & $\vdots$ & $\vdots$ & $\vdots$ \\ \midrule
			\multirow{2}{*}{$k$} & $D_{A}^{\, (k-1)}$, $D_{B}^{\, (k-1)}$ & $\mathcal{M}^{*}$ & $D_{C}^{\, (k-1)}$ & $D_{\mathrm{unused}}^{\, (k)} = D_{\mathrm{unused}}^{\, (k-1)} + D_{C, \, \mathrm{low}}^{\, (k-1)}$ \\
			& $D_{t+v}^{\, (k)} = D_{t+v}^{\, (k-1)} - D_{C, \, \mathrm{low}}^{\, (k-1)}$ & $\mathcal{M}^{\, (k)}$ & Any & - \\ \midrule
			$\vdots$ & $\vdots$ & $\vdots$ & $\vdots$ & $\vdots$ \\ \bottomrule
		\end{tabular}
	\end{table}
	
	\clearpage
	
	\section{Multi-layer gated recurrent unit (GRU) architecture \label{sec:gru_architecture}}
	
\begin{figure}[hbt]
		\centering
		\begin{subfigure}[b]{0.8\textwidth}
			\centering
			\includegraphics[width=\textwidth]{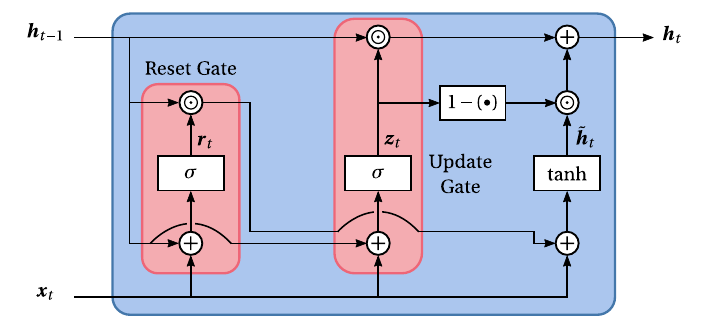}
			\caption{}
			\label{subfig:gru_layer}
		\end{subfigure}
		\begin{subfigure}[b]{0.8\textwidth}
			\centering
			\includegraphics[width=\textwidth]{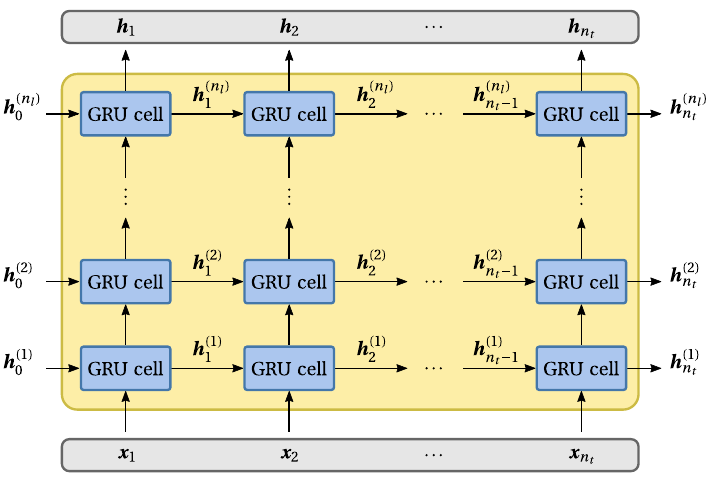}
			\caption{}
			\label{subfig:gru_multi_layer}
		\end{subfigure}
		\caption{Multi-layer gated recurrent unit (GRU) recurrent neural network model: \subref{subfig:gru_layer} GRU layer (or cell). The output of the reset gate determines the amount of old information that is used to generate new information, while the output of the update gate determines the amount of new information that updates the hidden (memory) state; \subref{subfig:gru_multi_layer} Multi-layer GRU architecture, consisting of a sequence of GRU layers, for which the output sequence is the output hidden state of the last layer.}
		\label{fig:gru_architecture}
	\end{figure}

	A multi-layer gated recurrent unit (GRU) \citep{cho:2014a, chung:2014a} is adopted as a recurrent neural network in both neural network and hybrid models architectures (see Figure~\ref{fig:gru_architecture}). Assume that the input sequence data is given by $\bm{x} \equiv [\bm{x}_{0}, \, \bm{x}_{1}, \, \dots, \, \bm{x}_{n_{t}}]$, where $\bm{x}_{t}$ is the $p$-feature input data at time $t$ and $n_{t}$ is the total sequence length. For each element $\bm{x}_{t}$, $t=0, \, \dots, \, n_{t}$, each GRU layer (see Figure~\ref{subfig:gru_layer}) computes the $h$-feature hidden state $\bm{h}_{t}$ by performing the following sequence of operations,
	\begin{align}
		\bm{z}_{t} &= \sigma \left( \bm{W}^{xz} \bm{x}_{t} + \bm{b}^{xz} + \bm{W}^{hz} \bm{h}_{t-1} + \bm{b}^{hz} \right) \, , \\
		\bm{r}_{t} &= \sigma \left( \bm{W}^{xr} \bm{x}_{t} + \bm{b}^{xr} + \bm{W}^{hr} \bm{h}_{t-1} + \bm{b}^{hr} \right) \, , \\
		\tilde{\bm{h}}_{t} &= \mathrm{tanh} \left( \bm{W}^{xh} \bm{x}_{t} + \bm{b}^{xh} + \bm{r}_{t} \odot \left( \bm{W}^{hh} \bm{h}_{t-1} + \bm{b}^{hh} \right) \right) \, , \\
		\bm{h}_{t} &= \bm{z}_{t} \odot \bm{h}_{t-1} + (\bm{1} - \bm{z}_{t}) \odot \tilde{\bm{h}}_{t} \, ,
	\end{align}
	where $\sigma$ denotes the sigmoid function, $\mathrm{tanh}$ is the hyperbolic tangent function, $\bm{z}_{t}$ is the $h$-feature output of the update gate, $\bm{r}_{t}$ is the $h$-feature output of the reset gate and $\tilde{\bm{h}}_{t}$ is the $h$-feature candidate hidden state. In turn, all instances of $\bm{W}$ and $\bm{b}$ denote either input-hidden or hidden-hidden weights and biases, i.e., the layer learnable parameters. In a multi-layer architecture with $n_{l}$ layers (see Figure~\ref{subfig:gru_multi_layer}), the input of any layer $l$, $l=2, \, \dots, \, n_{l}$, is the $h$-feature output $\bm{x}_{t}^{(l)} = \bm{h}_{t}^{(l-1)}$ of the previous layer, in which case $p=h$.
	
	In the context of this paper, the input is a sequence of 6-feature strain tensors, $\bm{\varepsilon} \equiv [\bm{\varepsilon}_{0}, \, \bm{\varepsilon}_{1}, \, \dots, \, \bm{\varepsilon}_{n_{t}}]$, while the output is a sequence of 6-feature stress tensors, $\bm{\sigma} \equiv [\bm{\sigma}_{0}, \, \bm{\sigma}_{1}, \, \dots, \, \bm{\sigma}_{n_{t}}]$. Given that the last GRU layer output has $h$-features, where $h$ is the number of hidden state features (or hidden layer size), a final linear layer is incorporated into the model such that
	\begin{equation}
		\bm{y}_{t} = \bm{W}^{hy} \bm{h}_{t} + \bm{b}^{hy} \, ,
	\end{equation}
	where $\bm{y} \equiv [\bm{y}_{0}, \, \bm{y}_{1}, \, \dots, \, \bm{y}_{n_{t}}]$ denotes the $q$-feature output sequence data. It is also remarked that the hidden state, $\bm{h}_{0}$, is initialized to zero in all layers.
	
\begin{figure}[hbt]
		\centering
		\includegraphics[width=0.5\textwidth]{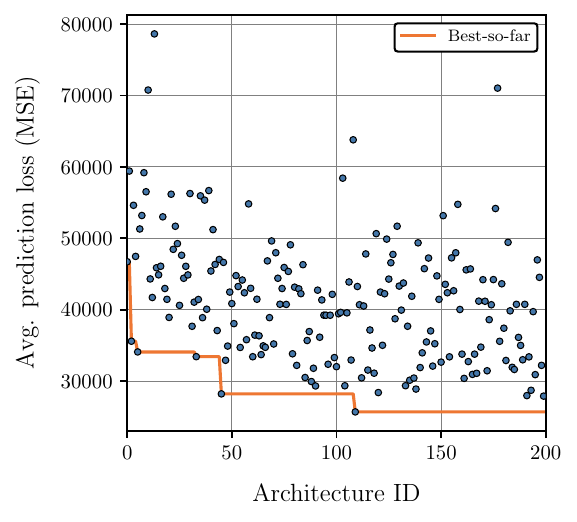}
		\caption{GRU material model hyperparameter optimization history with TPE (Tree-structured Parzen Estimator) algorithm. The hyperparameter search space includes the number of recurrent layers (1 to 6) and the hidden layer size (100 to 600), covering GRU models with a total number of trainable parameters between 33K and 12M, approximately. The objective function is the average testing loss (MSE) in a local synthetic data set of 512 random polynomial strain-stress paths generated with LZY model ($E=110$GPa, $\nu=0.33$, $a=1.0$, $b=0.05$, $c=1.0$, $d=0.5$, $s_{0}=900\,$MPa, $s_{1}=700\,$MPa, $s_{2}=0.5$).}
		\label{fig:gru_hyperparameter_opt_lzy}
	\end{figure}

	Having established the GRU model architecture, it remains to set the corresponding hyperparameters, namely the number of recurrent layers, $n_{l}$, and the hidden layer size, $h$. A hyperparameter optimization is performed using the Tree-structure Parzen Estimator (TPE) algorithm \citep{bergstra:2011a}. The search space is defined as $n_{l} \in [1, 6]$ and $h \in [100, 600]$, covering GRU models with a total number of learnable parameters between 33K and 12M, approximately. The objective function is the average testing loss (MSE) in a local synthetic data set of 512 random polynomial strain-stress paths generated with the LZY model ($E=110$GPa, $\nu=0.33$, $a=1.0$, $b=0.05$, $c=1.0$, $d=0.5$, $s_{0}=900\,$MPa, $s_{1}=700\,$MPa, $s_{2}=0.5$) after training the GRU model in a local synthetic data set of 320 random polynomial strain-stress paths. The hyperparameter optimization history is shown in Figure~\ref{fig:gru_hyperparameter_opt_lzy} and resulted in $n_{l}=3$ and $h=444$, yielding a GRU model with approximately 3M learnable parameters. This same architecture is used for both neural network and hybrid models.
	
	\clearpage
	
	\bibliography{mybibfile}

\begin{thebibliography}{87}
\expandafter\ifx\csname natexlab\endcsname\relax\def\natexlab#1{#1}\fi
\providecommand{\url}[1]{\texttt{#1}}
\providecommand{\href}[2]{#2}
\providecommand{\path}[1]{#1}
\providecommand{\DOIprefix}{doi:}
\providecommand{\ArXivprefix}{arXiv:}
\providecommand{\URLprefix}{URL: }
\providecommand{\Pubmedprefix}{pmid:}
\providecommand{\doi}[1]{\href{http://dx.doi.org/#1}{\path{#1}}}
\providecommand{\Pubmed}[1]{\href{pmid:#1}{\path{#1}}}
\providecommand{\bibinfo}[2]{#2}
\ifx\xfnm\relax \def\xfnm[#1]{\unskip,\space#1}\fi
\bibitem[{Lemaitre and Chaboche(1990)}]{lemaitre:1990a}
\bibinfo{author}{J.~Lemaitre}, \bibinfo{author}{J.~Chaboche},
  \bibinfo{title}{Mechanics of Solid Materials}, \bibinfo{publisher}{Cambridge
  University Press}, \bibinfo{year}{1990}.
\bibitem[{Simo and Hughes(2006)}]{simo:2006a}
\bibinfo{author}{J.~Simo}, \bibinfo{author}{T.~Hughes},
  \bibinfo{title}{Computational Inelasticity}, Interdisciplinary Applied
  Mathematics, \bibinfo{publisher}{Springer New York}, \bibinfo{year}{2006}.
\bibitem[{{de Souza Neto} et~al.(2008){de Souza Neto}, Peri, and
  Owen}]{desouzaneto:2008a}
\bibinfo{author}{E.~A. {de Souza Neto}}, \bibinfo{author}{D.~Peri},
  \bibinfo{author}{D.~R.~J. Owen}, \bibinfo{title}{Computational {{Methods}}
  for {{Plasticity}}}, \bibinfo{publisher}{{John Wiley \& Sons, Ltd}},
  \bibinfo{address}{{Chichester, UK}}, \bibinfo{year}{2008}.
\bibitem[{Bonet and Wood(2008)}]{bonet:2008a}
\bibinfo{author}{J.~Bonet}, \bibinfo{author}{R.~D. Wood},
  \bibinfo{title}{Nonlinear Continuum Mechanics for Finite Element Analysis},
  \bibinfo{edition}{2} ed., \bibinfo{publisher}{Cambridge University Press},
  \bibinfo{year}{2008}.
\bibitem[{Gurtin et~al.(2010)Gurtin, Fried, and Anand}]{gurtin:2010a}
\bibinfo{author}{M.~E. Gurtin}, \bibinfo{author}{E.~Fried},
  \bibinfo{author}{L.~Anand}, \bibinfo{title}{The Mechanics and Thermodynamics
  of Continua}, \bibinfo{publisher}{Cambridge University Press},
  \bibinfo{year}{2010}.
\bibitem[{Belytschko et~al.(2013)Belytschko, Liu, Moran, and
  Elkhodary}]{belytschko:2013a}
\bibinfo{author}{T.~Belytschko}, \bibinfo{author}{W.~Liu},
  \bibinfo{author}{B.~Moran}, \bibinfo{author}{K.~Elkhodary},
  \bibinfo{title}{Nonlinear Finite Elements for Continua and Structures}, No
  Longer used, \bibinfo{publisher}{Wiley}, \bibinfo{year}{2013}.
\bibitem[{Hill(1950)}]{hill:1950a}
\bibinfo{author}{R.~Hill}, \bibinfo{title}{The Mathematical Theory of
  Plasticity}, \bibinfo{publisher}{Oxford University Press},
  \bibinfo{year}{1950}.
\bibitem[{Truesdell and Noll(1965)}]{truesdell:1965a}
\bibinfo{author}{C.~Truesdell}, \bibinfo{author}{W.~Noll}, \bibinfo{title}{The
  Non-Linear Field Theories of Mechanics}, \bibinfo{publisher}{Springer Berlin
  Heidelberg}, \bibinfo{address}{Berlin, Heidelberg}, \bibinfo{year}{1965}, pp.
  \bibinfo{pages}{1--541}.
\bibitem[{Malvern(1969)}]{malvern:1969a}
\bibinfo{author}{L.~Malvern}, \bibinfo{title}{Introduction to the Mechanics of
  a Continuous Medium}, Prentice-Hall series in engineering of the physical
  sciences, \bibinfo{publisher}{Prentice-Hall}, \bibinfo{year}{1969}.
\bibitem[{Ogden(1984)}]{ogden:1984a}
\bibinfo{author}{R.~Ogden}, \bibinfo{title}{Non-linear Elastic Deformations},
  Ellis Horwood series in mathematics and its applications,
  \bibinfo{publisher}{E. Horwood}, \bibinfo{year}{1984}.
\bibitem[{Gurson(1977)}]{gurson:1977a}
\bibinfo{author}{A.~L. Gurson},
\newblock \bibinfo{title}{Continuum theory of ductile rupture by void
  nucleation and growth: Part i—yield criteria and flow rules for porous
  ductile media},
\newblock \bibinfo{journal}{Journal of Engineering Materials and Technology}
  \bibinfo{volume}{99} (\bibinfo{year}{1977}) \bibinfo{pages}{2--15}.
\bibitem[{Kothari and Anand(1998)}]{kothari:1998a}
\bibinfo{author}{M.~Kothari}, \bibinfo{author}{L.~Anand},
\newblock \bibinfo{title}{Elasto-viscoplastic constitutive equations for
  polycrystalline metals: Application to tantalum},
\newblock \bibinfo{journal}{Journal of the Mechanics and Physics of Solids}
  \bibinfo{volume}{46} (\bibinfo{year}{1998}) \bibinfo{pages}{51--83}.
\bibitem[{Holzapfel et~al.(2000)Holzapfel, Gasser, and Ogden}]{holzapfel:2000a}
\bibinfo{author}{G.~A. Holzapfel}, \bibinfo{author}{T.~C. Gasser},
  \bibinfo{author}{R.~W. Ogden},
\newblock \bibinfo{title}{A new constitutive framework for arterial wall
  mechanics and a comparative study of material models},
\newblock \bibinfo{journal}{J. Elast.} \bibinfo{volume}{61}
  (\bibinfo{year}{2000}) \bibinfo{pages}{1--48}.
\bibitem[{Anand and Gurtin(2003)}]{anand:2003a}
\bibinfo{author}{L.~Anand}, \bibinfo{author}{M.~E. Gurtin},
\newblock \bibinfo{title}{A theory of amorphous solids undergoing large
  deformations, with application to polymeric glasses},
\newblock \bibinfo{journal}{International Journal of Solids and Structures}
  \bibinfo{volume}{40} (\bibinfo{year}{2003}) \bibinfo{pages}{1465--1487}.
\bibitem[{Chen and Chen(1995)}]{chen:1995a}
\bibinfo{author}{T.~Chen}, \bibinfo{author}{H.~Chen},
\newblock \bibinfo{title}{Universal approximation to nonlinear operators by
  neural networks with arbitrary activation functions and its application to
  dynamical systems},
\newblock \bibinfo{journal}{IEEE Transactions on Neural Networks}
  \bibinfo{volume}{6} (\bibinfo{year}{1995}) \bibinfo{pages}{911--917}.
\bibitem[{Mozaffar et~al.(2019)Mozaffar, Bostanabad, Chen, Ehmann, Cao, and
  Bessa}]{mozaffar:2019}
\bibinfo{author}{M.~Mozaffar}, \bibinfo{author}{R.~Bostanabad},
  \bibinfo{author}{W.~Chen}, \bibinfo{author}{K.~Ehmann},
  \bibinfo{author}{J.~Cao}, \bibinfo{author}{M.~A. Bessa},
\newblock \bibinfo{title}{Deep learning predicts path-dependent plasticity},
\newblock \bibinfo{journal}{Proceedings of the National Academy of Sciences}
  \bibinfo{volume}{116} (\bibinfo{year}{2019}) \bibinfo{pages}{26414--26420}.
\bibitem[{Bessa et~al.(2017)Bessa, Bostanabad, Liu, Hu, Apley, Brinson, Chen,
  and Liu}]{bessa:2017a}
\bibinfo{author}{M.~Bessa}, \bibinfo{author}{R.~Bostanabad},
  \bibinfo{author}{Z.~Liu}, \bibinfo{author}{A.~Hu}, \bibinfo{author}{D.~W.
  Apley}, \bibinfo{author}{C.~Brinson}, \bibinfo{author}{W.~Chen},
  \bibinfo{author}{W.~Liu},
\newblock \bibinfo{title}{A framework for data-driven analysis of materials
  under uncertainty: Countering the curse of dimensionality},
\newblock \bibinfo{journal}{Computer Methods in Applied Mechanics and
  Engineering} \bibinfo{volume}{320} (\bibinfo{year}{2017})
  \bibinfo{pages}{633--667}.
\bibitem[{Bengio(2009)}]{bengio:2009a}
\bibinfo{author}{Y.~Bengio}, \bibinfo{title}{Learning Deep Architectures for
  AI}, \bibinfo{year}{2009}.
\bibitem[{Rolnick and Tegmark(2017)}]{rolnick:2017a}
\bibinfo{author}{D.~Rolnick}, \bibinfo{author}{M.~Tegmark},
\newblock \bibinfo{title}{The power of deeper networks for expressing natural
  functions},
\newblock \bibinfo{journal}{ArXiv} \bibinfo{volume}{abs/1705.05502}
  (\bibinfo{year}{2017}).
\bibitem[{Ghaboussi et~al.(1991)Ghaboussi, Garrett, and Wu}]{ghaboussi:1991a}
\bibinfo{author}{J.~Ghaboussi}, \bibinfo{author}{J.~H. Garrett},
  \bibinfo{author}{X.~Wu},
\newblock \bibinfo{title}{Knowledge‐based modeling of material behavior with
  neural networks},
\newblock \bibinfo{journal}{Journal of Engineering Mechanics}
  \bibinfo{volume}{117} (\bibinfo{year}{1991}) \bibinfo{pages}{132--153}.
\bibitem[{Theocaris and Panagiotopoulos(1993)}]{theocaris:1993a}
\bibinfo{author}{P.~Theocaris}, \bibinfo{author}{P.~Panagiotopoulos},
\newblock \bibinfo{title}{Neural networks for computing in fracture mechanics.
  methods and prospects of applications},
\newblock \bibinfo{journal}{Computer Methods in Applied Mechanics and
  Engineering} \bibinfo{volume}{106} (\bibinfo{year}{1993})
  \bibinfo{pages}{213--228}.
\bibitem[{Yagawa and Okuda(1996)}]{yagawa:1996a}
\bibinfo{author}{G.~Yagawa}, \bibinfo{author}{H.~Okuda},
\newblock \bibinfo{title}{Neural networks in computational mechanics},
\newblock \bibinfo{journal}{Archives of Computational Methods in Engineering}
  \bibinfo{volume}{3} (\bibinfo{year}{1996}) \bibinfo{pages}{435--512}.
\bibitem[{Waszczyszyn and Ziemiański(2001)}]{waszczyszyn:2001a}
\bibinfo{author}{Z.~Waszczyszyn}, \bibinfo{author}{L.~Ziemiański},
\newblock \bibinfo{title}{Neural networks in mechanics of structures and
  materials – new results and prospects of applications},
\newblock \bibinfo{journal}{Computers \& Structures} \bibinfo{volume}{79}
  (\bibinfo{year}{2001}) \bibinfo{pages}{2261--2276}.
\bibitem[{Lefik and Schrefler(2002)}]{lefik:2002a}
\bibinfo{author}{M.~Lefik}, \bibinfo{author}{B.~Schrefler},
\newblock \bibinfo{title}{One-dimensional model of cable-in-conduit
  superconductors under cyclic loading using artificial neural networks},
\newblock \bibinfo{journal}{Fusion Engineering and Design} \bibinfo{volume}{60}
  (\bibinfo{year}{2002}) \bibinfo{pages}{105--117}.
\bibitem[{Ghavamian and Simone(2019)}]{ghavamian:2019a}
\bibinfo{author}{F.~Ghavamian}, \bibinfo{author}{A.~Simone},
\newblock \bibinfo{title}{Accelerating multiscale finite element simulations of
  history-dependent materials using a recurrent neural network},
\newblock \bibinfo{journal}{Computer Methods in Applied Mechanics and
  Engineering} \bibinfo{volume}{357} (\bibinfo{year}{2019})
  \bibinfo{pages}{112594}.
\bibitem[{Wu et~al.(2020)Wu, Nguyen, Kilingar, and Noels}]{wu:2020a}
\bibinfo{author}{L.~Wu}, \bibinfo{author}{V.~D. Nguyen}, \bibinfo{author}{N.~G.
  Kilingar}, \bibinfo{author}{L.~Noels},
\newblock \bibinfo{title}{A recurrent neural network-accelerated multi-scale
  model for elasto-plastic heterogeneous materials subjected to random cyclic
  and non-proportional loading paths},
\newblock \bibinfo{journal}{Computer Methods in Applied Mechanics and
  Engineering} \bibinfo{volume}{369} (\bibinfo{year}{2020})
  \bibinfo{pages}{113234}.
\bibitem[{Fuhg et~al.(2025)Fuhg, {Anantha Padmanabha}, Bouklas, Bahmani, Sun,
  Vlassis, Flaschel, Carrara, and {De Lorenzis}}]{fuhg:2025a}
\bibinfo{author}{J.~Fuhg}, \bibinfo{author}{G.~{Anantha Padmanabha}},
  \bibinfo{author}{N.~Bouklas}, \bibinfo{author}{B.~Bahmani},
  \bibinfo{author}{W.~Sun}, \bibinfo{author}{N.~Vlassis},
  \bibinfo{author}{M.~Flaschel}, \bibinfo{author}{P.~Carrara},
  \bibinfo{author}{L.~{De Lorenzis}},
\newblock \bibinfo{title}{A review on data-driven constitutive laws for
  solids},
\newblock \bibinfo{journal}{Archives of Computational Methods in Engineering}
  \bibinfo{volume}{32} (\bibinfo{year}{2025}) \bibinfo{pages}{1841--1883}.
  \bibinfo{note}{Publisher Copyright: {\textcopyright} The Author(s) under
  exclusive licence to International Center for Numerical Methods in
  Engineering (CIMNE) 2024.}
\bibitem[{Lefik et~al.(2009)Lefik, Boso, and Schrefler}]{lefik:2009a}
\bibinfo{author}{M.~Lefik}, \bibinfo{author}{D.~Boso},
  \bibinfo{author}{B.~Schrefler},
\newblock \bibinfo{title}{Artificial neural networks in numerical modelling of
  composites},
\newblock \bibinfo{journal}{Computer Methods in Applied Mechanics and
  Engineering} \bibinfo{volume}{198} (\bibinfo{year}{2009})
  \bibinfo{pages}{1785--1804}. \bibinfo{note}{Advances in Simulation-Based
  Engineering Sciences – Honoring J. Tinsley Oden}.
\bibitem[{Unger and Könke(2009)}]{unger:2009a}
\bibinfo{author}{J.~F. Unger}, \bibinfo{author}{C.~Könke},
\newblock \bibinfo{title}{Neural networks as material models within a
  multiscale approach},
\newblock \bibinfo{journal}{Computers \& Structures} \bibinfo{volume}{87}
  (\bibinfo{year}{2009}) \bibinfo{pages}{1177--1186}. \bibinfo{note}{Civil-Comp
  Special Issue}.
\bibitem[{Le et~al.(2015)Le, Yvonnet, and He}]{le:2015a}
\bibinfo{author}{B.~A. Le}, \bibinfo{author}{J.~Yvonnet},
  \bibinfo{author}{Q.-C. He},
\newblock \bibinfo{title}{Computational homogenization of nonlinear elastic
  materials using neural networks},
\newblock \bibinfo{journal}{International Journal for Numerical Methods in
  Engineering} \bibinfo{volume}{104} (\bibinfo{year}{2015})
  \bibinfo{pages}{1061--1084}.
\bibitem[{Wang and Sun(2018)}]{wang:2018a}
\bibinfo{author}{K.~Wang}, \bibinfo{author}{W.~Sun},
\newblock \bibinfo{title}{A multiscale multi-permeability poroplasticity model
  linked by recursive homogenizations and deep learning},
\newblock \bibinfo{journal}{Computer Methods in Applied Mechanics and
  Engineering} \bibinfo{volume}{334} (\bibinfo{year}{2018})
  \bibinfo{pages}{337--380}.
\bibitem[{Liu et~al.(2019)Liu, Wu, and Koishi}]{liu:2019a}
\bibinfo{author}{Z.~Liu}, \bibinfo{author}{C.~Wu}, \bibinfo{author}{M.~Koishi},
\newblock \bibinfo{title}{A deep material network for multiscale topology
  learning and accelerated nonlinear modeling of heterogeneous materials},
\newblock \bibinfo{journal}{Computer Methods in Applied Mechanics and
  Engineering} \bibinfo{volume}{345} (\bibinfo{year}{2019})
  \bibinfo{pages}{1138--1168}.
\bibitem[{He and Chen(2022)}]{he:2022a}
\bibinfo{author}{X.~He}, \bibinfo{author}{J.-S. Chen},
\newblock \bibinfo{title}{Thermodynamically consistent machine-learned internal
  state variable approach for data-driven modeling of path-dependent
  materials},
\newblock \bibinfo{journal}{Computer Methods in Applied Mechanics and
  Engineering} \bibinfo{volume}{402} (\bibinfo{year}{2022})
  \bibinfo{pages}{115348}. \bibinfo{note}{A Special Issue in Honor of the
  Lifetime Achievements of J. Tinsley Oden}.
\bibitem[{Eghbalian et~al.(2023)Eghbalian, Pouragha, and Wan}]{eghbalian:2023a}
\bibinfo{author}{M.~Eghbalian}, \bibinfo{author}{M.~Pouragha},
  \bibinfo{author}{R.~Wan},
\newblock \bibinfo{title}{A physics-informed deep neural network for surrogate
  modeling in classical elasto-plasticity},
\newblock \bibinfo{journal}{Computers and Geotechnics} \bibinfo{volume}{159}
  (\bibinfo{year}{2023}) \bibinfo{pages}{105472}.
\bibitem[{Kalina et~al.(2025)Kalina, Brummund, Sun, and
  Kästner}]{kalina:2025a}
\bibinfo{author}{K.~A. Kalina}, \bibinfo{author}{J.~Brummund},
  \bibinfo{author}{W.~Sun}, \bibinfo{author}{M.~Kästner},
\newblock \bibinfo{title}{Neural networks meet anisotropic hyperelasticity: A
  framework based on generalized structure tensors and isotropic tensor
  functions},
\newblock \bibinfo{journal}{Computer Methods in Applied Mechanics and
  Engineering} \bibinfo{volume}{437} (\bibinfo{year}{2025})
  \bibinfo{pages}{117725}.
\bibitem[{Jadoon et~al.(2025)Jadoon, Meyer, and Fuhg}]{jadoon:2025a}
\bibinfo{author}{A.~A. Jadoon}, \bibinfo{author}{K.~A. Meyer},
  \bibinfo{author}{J.~N. Fuhg},
\newblock \bibinfo{title}{Automated model discovery of finite strain
  elastoplasticity from uniaxial experiments},
\newblock \bibinfo{journal}{Computer Methods in Applied Mechanics and
  Engineering} \bibinfo{volume}{435} (\bibinfo{year}{2025})
  \bibinfo{pages}{117653}.
\bibitem[{Kavanagh and Clough(1971)}]{kavanagh:1971a}
\bibinfo{author}{K.~T. Kavanagh}, \bibinfo{author}{R.~W. Clough},
\newblock \bibinfo{title}{Finite element applications in the characterization
  of elastic solids},
\newblock \bibinfo{journal}{International Journal of Solids and Structures}
  \bibinfo{volume}{7} (\bibinfo{year}{1971}) \bibinfo{pages}{11--23}.
\bibitem[{Chen et~al.(2024)Chen, Starman, Halilovič, Berglund, and
  Coppieters}]{chen:2024a}
\bibinfo{author}{B.~Chen}, \bibinfo{author}{B.~Starman},
  \bibinfo{author}{M.~Halilovič}, \bibinfo{author}{L.~A. Berglund},
  \bibinfo{author}{S.~Coppieters},
\newblock \bibinfo{title}{Finite element model updating for material model
  calibration: A review and guide to practice},
\newblock \bibinfo{journal}{Archives of Computational Methods in Engineering}
  (\bibinfo{year}{2024}).
\bibitem[{Bay et~al.(1999)Bay, Smith, Fyhrie, and Saad}]{bay:1999a}
\bibinfo{author}{B.~K. Bay}, \bibinfo{author}{T.~S. Smith},
  \bibinfo{author}{D.~P. Fyhrie}, \bibinfo{author}{M.~Saad},
\newblock \bibinfo{title}{Digital volume correlation: three-dimensional strain
  mapping using x-ray tomography},
\newblock \bibinfo{journal}{Experimental mechanics} \bibinfo{volume}{39}
  (\bibinfo{year}{1999}) \bibinfo{pages}{217--226}.
\bibitem[{Sutton et~al.(2009)Sutton, Orteu, and Schreier}]{sutton:2009a}
\bibinfo{author}{M.~Sutton}, \bibinfo{author}{J.~Orteu},
  \bibinfo{author}{H.~Schreier}, \bibinfo{title}{Image Correlation for Shape,
  Motion and Deformation Measurements: Basic Concepts,Theory and Applications},
  \bibinfo{publisher}{Springer US}, \bibinfo{year}{2009}.
\bibitem[{Pan et~al.(2009)Pan, Qian, Xie, and Asundi}]{pan:2009a}
\bibinfo{author}{B.~Pan}, \bibinfo{author}{K.~Qian}, \bibinfo{author}{H.~Xie},
  \bibinfo{author}{A.~Asundi},
\newblock \bibinfo{title}{Two-dimensional digital image correlation for
  in-plane displacement and strain measurement: a review},
\newblock \bibinfo{journal}{Measurement Science and Technology}
  \bibinfo{volume}{20} (\bibinfo{year}{2009}) \bibinfo{pages}{062001}.
\bibitem[{Buljac et~al.(2018)Buljac, {Trejo Navas}, Shakoor, Bouterf, Neggers,
  Bernacki, Bouchard, Morgeneyer, and Hild}]{buljac:2018a}
\bibinfo{author}{A.~Buljac}, \bibinfo{author}{V.-M. {Trejo Navas}},
  \bibinfo{author}{M.~Shakoor}, \bibinfo{author}{A.~Bouterf},
  \bibinfo{author}{J.~Neggers}, \bibinfo{author}{M.~Bernacki},
  \bibinfo{author}{P.-O. Bouchard}, \bibinfo{author}{T.~F. Morgeneyer},
  \bibinfo{author}{F.~Hild},
\newblock \bibinfo{title}{On the calibration of elastoplastic parameters at the
  microscale via x-ray microtomography and digital volume correlation for the
  simulation of ductile damage},
\newblock \bibinfo{journal}{European Journal of Mechanics - A/Solids}
  \bibinfo{volume}{72} (\bibinfo{year}{2018}) \bibinfo{pages}{287--297}.
\bibitem[{Wang et~al.(2024)Wang, Das, Joshi, Shaikeea, and
  Deshpande}]{wang:2024a}
\bibinfo{author}{Z.~Wang}, \bibinfo{author}{S.~Das},
  \bibinfo{author}{A.~Joshi}, \bibinfo{author}{A.~J.~D. Shaikeea},
  \bibinfo{author}{V.~S. Deshpande},
\newblock \bibinfo{title}{3d observations provide striking findings in rubber
  elasticity},
\newblock \bibinfo{journal}{Proceedings of the National Academy of Sciences}
  \bibinfo{volume}{121} (\bibinfo{year}{2024}) \bibinfo{pages}{e2404205121}.
\bibitem[{Choi and Kim(2006)}]{choi:2006a}
\bibinfo{author}{K.~Choi}, \bibinfo{author}{N.~Kim}, \bibinfo{title}{Structural
  Sensitivity Analysis and Optimization 1: Linear Systems}, Mechanical
  Engineering Series, \bibinfo{publisher}{Springer New York},
  \bibinfo{year}{2006}.
\bibitem[{Mr{\'o}z and Stavroulakis(2007)}]{mroz:2007a}
\bibinfo{author}{Z.~Mr{\'o}z}, \bibinfo{author}{G.~Stavroulakis},
  \bibinfo{title}{Parameter Identification of Materials and Structures}, CISM
  International Centre for Mechanical Sciences, \bibinfo{publisher}{Springer
  Vienna}, \bibinfo{year}{2007}.
\bibitem[{Friswell and Mottershead(2010)}]{friswell:2010a}
\bibinfo{author}{M.~Friswell}, \bibinfo{author}{J.~Mottershead},
  \bibinfo{title}{Finite Element Model Updating in Structural Dynamics}, Solid
  Mechanics and Its Applications, \bibinfo{publisher}{Springer Netherlands},
  \bibinfo{year}{2010}.
\bibitem[{Haftka et~al.(2013)Haftka, G{\"u}rdal, and Kamat}]{haftka:2013a}
\bibinfo{author}{R.~Haftka}, \bibinfo{author}{Z.~G{\"u}rdal},
  \bibinfo{author}{M.~Kamat}, \bibinfo{title}{Elements of Structural
  Optimization}, Solid Mechanics and Its Applications,
  \bibinfo{publisher}{Springer Netherlands}, \bibinfo{year}{2013}.
\bibitem[{Wengert(1964)}]{wengert:1964a}
\bibinfo{author}{R.~E. Wengert}, \bibinfo{title}{The Reverse Mode of
  Differentiation}, \bibinfo{type}{Technical Report}, University of California,
  \bibinfo{year}{1964}. \URLprefix
  \url{https://softlib.rice.edu/pub/CRPC-TRs/reports/CRPC-TR89003.pdf},
  \bibinfo{note}{accessed: 2025-05-07}.
\bibitem[{Griewank(1989)}]{griewank:1989a}
\bibinfo{author}{A.~Griewank}, \bibinfo{title}{Achieving Derivatives of
  High-Order by Automatic Differentiation}, \bibinfo{type}{Technical Report},
  CRPC Technical Report, \bibinfo{year}{1989}. \URLprefix
  \url{https://softlib.rice.edu/pub/CRPC-TRs/reports/CRPC-TR89003.pdf},
  \bibinfo{note}{accessed: 2025-05-07}.
\bibitem[{Baydin et~al.(2018)Baydin, Pearlmutter, Radul, and
  Siskind}]{baydin:2018a}
\bibinfo{author}{A.~G. Baydin}, \bibinfo{author}{B.~A. Pearlmutter},
  \bibinfo{author}{A.~A. Radul}, \bibinfo{author}{J.~M. Siskind},
\newblock \bibinfo{title}{Automatic differentiation in machine learning: a
  survey},
\newblock \bibinfo{journal}{Journal of machine learning research}
  \bibinfo{volume}{18} (\bibinfo{year}{2018}) \bibinfo{pages}{1--43}.
\bibitem[{Thakolkaran et~al.(2022)Thakolkaran, Joshi, Zheng, Flaschel, {De
  Lorenzis}, and Kumar}]{thakolkaran:2022a}
\bibinfo{author}{P.~Thakolkaran}, \bibinfo{author}{A.~Joshi},
  \bibinfo{author}{Y.~Zheng}, \bibinfo{author}{M.~Flaschel},
  \bibinfo{author}{L.~{De Lorenzis}}, \bibinfo{author}{S.~Kumar},
\newblock \bibinfo{title}{Nn-euclid: Deep-learning hyperelasticity without
  stress data},
\newblock \bibinfo{journal}{Journal of the Mechanics and Physics of Solids}
  \bibinfo{volume}{169} (\bibinfo{year}{2022}) \bibinfo{pages}{105076}.
\bibitem[{Wu et~al.(2025)Wu, Zhang, and Mao}]{wu:2025a}
\bibinfo{author}{X.~Wu}, \bibinfo{author}{Y.~Zhang}, \bibinfo{author}{S.~Mao},
\newblock \bibinfo{title}{Learning the physics-consistent material behavior
  from measurable data via pde-constrained optimization},
\newblock \bibinfo{journal}{Computer Methods in Applied Mechanics and
  Engineering} \bibinfo{volume}{437} (\bibinfo{year}{2025})
  \bibinfo{pages}{117748}.
\bibitem[{Akerson et~al.(2025)Akerson, Rajan, and Bhattacharya}]{akerson:2025a}
\bibinfo{author}{A.~Akerson}, \bibinfo{author}{A.~Rajan},
  \bibinfo{author}{K.~Bhattacharya},
\newblock \bibinfo{title}{Learning constitutive relations from experiments: 1.
  pde constrained optimization},
\newblock \bibinfo{journal}{Journal of the Mechanics and Physics of Solids}
  \bibinfo{volume}{201} (\bibinfo{year}{2025}) \bibinfo{pages}{106128}.
\bibitem[{Ghaboussi et~al.(1998)Ghaboussi, Pecknold, Zhang, and
  Haj-Ali}]{ghaboussi:1998a}
\bibinfo{author}{J.~Ghaboussi}, \bibinfo{author}{D.~A. Pecknold},
  \bibinfo{author}{M.~Zhang}, \bibinfo{author}{R.~M. Haj-Ali},
\newblock \bibinfo{title}{Autoprogressive training of neural network
  constitutive models},
\newblock \bibinfo{journal}{International Journal for Numerical Methods in
  Engineering} \bibinfo{volume}{42} (\bibinfo{year}{1998})
  \bibinfo{pages}{105--126}.
\bibitem[{Shin and Pande(2000)}]{shin:2000a}
\bibinfo{author}{H.~Shin}, \bibinfo{author}{G.~Pande},
\newblock \bibinfo{title}{On self-learning finite element codes based on
  monitored response of structures},
\newblock \bibinfo{journal}{Computers and Geotechnics} \bibinfo{volume}{27}
  (\bibinfo{year}{2000}) \bibinfo{pages}{161--178}.
\bibitem[{Lefik and Schrefler(2003)}]{lefik:2003a}
\bibinfo{author}{M.~Lefik}, \bibinfo{author}{B.~Schrefler},
\newblock \bibinfo{title}{Artificial neural network as an incremental
  non-linear constitutive model for a finite element code},
\newblock \bibinfo{journal}{Computer Methods in Applied Mechanics and
  Engineering} \bibinfo{volume}{192} (\bibinfo{year}{2003})
  \bibinfo{pages}{3265--3283}. \bibinfo{note}{Multiscale Computational
  Mechanics for Materials and Structures}.
\bibitem[{Gr{\'e}diac(1989)}]{grediac:1989a}
\bibinfo{author}{M.~Gr{\'e}diac},
\newblock \bibinfo{title}{Principe des travaux virtuels et identification},
\newblock \bibinfo{journal}{Comptes rendus de l'Acad{\'e}mie des sciences.
  S{\'e}rie 2, M{\'e}canique, Physique, Chimie, Sciences de l'univers, Sciences
  de la Terre} \bibinfo{volume}{309} (\bibinfo{year}{1989})
  \bibinfo{pages}{1--5}.
\bibitem[{Pierron and Gr{\'e}diac(2012)}]{pierron:2012a}
\bibinfo{author}{F.~Pierron}, \bibinfo{author}{M.~Gr{\'e}diac},
  \bibinfo{title}{The virtual fields method. Extracting constitutive mechanical
  parameters from full-field deformation measurements},
  \bibinfo{publisher}{Springer}, \bibinfo{year}{2012}.
\bibitem[{Lourenço et~al.(2024)Lourenço, Georgieva, Cueto, and
  Andrade-Campos}]{lourenco:2024a}
\bibinfo{author}{R.~Lourenço}, \bibinfo{author}{P.~Georgieva},
  \bibinfo{author}{E.~Cueto}, \bibinfo{author}{A.~Andrade-Campos},
\newblock \bibinfo{title}{An indirect training approach for implicit
  constitutive modelling using recurrent neural networks and the virtual fields
  method},
\newblock \bibinfo{journal}{Computer Methods in Applied Mechanics and
  Engineering} \bibinfo{volume}{425} (\bibinfo{year}{2024})
  \bibinfo{pages}{116961}.
\bibitem[{Sun et~al.(2025)Sun, Taguchi, Niki, Iizuka, and Yoneyama}]{sun:2025a}
\bibinfo{author}{D.~Sun}, \bibinfo{author}{S.~Taguchi},
  \bibinfo{author}{I.~Niki}, \bibinfo{author}{K.~Iizuka},
  \bibinfo{author}{S.~Yoneyama},
\newblock \bibinfo{title}{The virtual fields method for identifying
  viscoelastic properties based on stress-sensitivity virtual fields},
\newblock \bibinfo{journal}{Mechanics of Time-Dependent Materials}
  \bibinfo{volume}{29} (\bibinfo{year}{2025}) \bibinfo{pages}{43}.
\bibitem[{Kumar et~al.(2025)Kumar, Seidl, Granzow, Yang, and
  Fuhg}]{kumar:2025a}
\bibinfo{author}{S.~Kumar}, \bibinfo{author}{D.~T. Seidl},
  \bibinfo{author}{B.~N. Granzow}, \bibinfo{author}{J.~Yang},
  \bibinfo{author}{J.~N. Fuhg},
\newblock \bibinfo{title}{A comparative study of calibration techniques for
  finite strain elastoplasticity: Numerically-exact sensitivities for femu and
  vfm},
\newblock \bibinfo{journal}{arXiv preprint arXiv:2503.19782}
  (\bibinfo{year}{2025}).
\bibitem[{Dissanayake and Phan-Thien(1994)}]{dissanayake:1994a}
\bibinfo{author}{M.~W. M.~G. Dissanayake}, \bibinfo{author}{N.~Phan-Thien},
\newblock \bibinfo{title}{Neural-network-based approximations for solving
  partial differential equations},
\newblock \bibinfo{journal}{Communications in Numerical Methods in Engineering}
  \bibinfo{volume}{10} (\bibinfo{year}{1994}) \bibinfo{pages}{195--201}.
\bibitem[{Raissi et~al.(2019)Raissi, Perdikaris, and
  Karniadakis}]{raissi:2019a}
\bibinfo{author}{M.~Raissi}, \bibinfo{author}{P.~Perdikaris},
  \bibinfo{author}{G.~Karniadakis},
\newblock \bibinfo{title}{Physics-informed neural networks: A deep learning
  framework for solving forward and inverse problems involving nonlinear
  partial differential equations},
\newblock \bibinfo{journal}{Journal of Computational Physics}
  \bibinfo{volume}{378} (\bibinfo{year}{2019}) \bibinfo{pages}{686--707}.
\bibitem[{Jeong et~al.(2024)Jeong, Cho, Chung, and Kim}]{jeong:2024a}
\bibinfo{author}{I.~Jeong}, \bibinfo{author}{M.~Cho},
  \bibinfo{author}{H.~Chung}, \bibinfo{author}{D.-N. Kim},
\newblock \bibinfo{title}{Data-driven nonparametric identification of material
  behavior based on physics-informed neural network with full-field data},
\newblock \bibinfo{journal}{Computer Methods in Applied Mechanics and
  Engineering} \bibinfo{volume}{418} (\bibinfo{year}{2024})
  \bibinfo{pages}{116569}.
\bibitem[{Anton et~al.(2024)Anton, Tr{\"o}ger, Wessels, R{\"o}mer, Henkes, and
  Hartmann}]{anton:2024a}
\bibinfo{author}{D.~Anton}, \bibinfo{author}{J.-A. Tr{\"o}ger},
  \bibinfo{author}{H.~Wessels}, \bibinfo{author}{U.~R{\"o}mer},
  \bibinfo{author}{A.~Henkes}, \bibinfo{author}{S.~Hartmann},
\newblock \bibinfo{title}{Deterministic and statistical calibration of
  constitutive models from full-field data with parametric physics-informed
  neural networks},
\newblock \bibinfo{journal}{arXiv preprint arXiv:2405.18311}
  (\bibinfo{year}{2024}).
\bibitem[{Zhang and Bhattacharya(2024)}]{zhang:2024a}
\bibinfo{author}{Y.~Zhang}, \bibinfo{author}{K.~Bhattacharya},
\newblock \bibinfo{title}{Iterated learning and multiscale modeling of
  history-dependent architectured metamaterials},
\newblock \bibinfo{journal}{Mechanics of Materials} \bibinfo{volume}{197}
  (\bibinfo{year}{2024}) \bibinfo{pages}{105090}.
\bibitem[{von Mises(1913)}]{mises:1913}
\bibinfo{author}{R.~von Mises},
\newblock \bibinfo{title}{Mechanik der festen k{\"o}rper im
  plastisch-deformablen zustand},
\newblock \bibinfo{journal}{Nachrichten von der Gesellschaft der Wissenschaften
  zu G{\"o}ttingen, Mathematisch-Physikalische Klasse}  (\bibinfo{year}{1913})
  \bibinfo{pages}{582--592}.
\bibitem[{Drucker and Prager(1952)}]{drucker:1952a}
\bibinfo{author}{D.~C. Drucker}, \bibinfo{author}{W.~Prager},
\newblock \bibinfo{title}{Soil mechanics and plastic analysis or limit design},
\newblock \bibinfo{journal}{Quarterly of Applied Mathematics}
  \bibinfo{volume}{10} (\bibinfo{year}{1952}) \bibinfo{pages}{157--165}.
\bibitem[{Kingma and Ba(2017)}]{kingma:2014a}
\bibinfo{author}{D.~P. Kingma}, \bibinfo{author}{J.~Ba}, \bibinfo{title}{Adam:
  A method for stochastic optimization}, \bibinfo{year}{2017}.
  \href{http://arxiv.org/abs/1412.6980}{\tt arXiv:1412.6980}.
\bibitem[{Giton et~al.(2006)Giton, Caro-Bretelle, and
  Ienny}]{giton2006hyperelastic}
\bibinfo{author}{M.~Giton}, \bibinfo{author}{A.~Caro-Bretelle},
  \bibinfo{author}{P.~Ienny},
\newblock \bibinfo{title}{Hyperelastic behaviour identification by a forward
  problem resolution: Application to a tear test of a silicone-rubber},
\newblock \bibinfo{journal}{Strain} \bibinfo{volume}{42}
  (\bibinfo{year}{2006}).
\bibitem[{Avril et~al.(2008)Avril, Bonnet, Bretelle, Gr{\'e}diac, Hild, Ienny,
  Latourte, Lemosse, Pagano, Pagnacco et~al.}]{avril2008overview}
\bibinfo{author}{S.~Avril}, \bibinfo{author}{M.~Bonnet}, \bibinfo{author}{A.-S.
  Bretelle}, \bibinfo{author}{M.~Gr{\'e}diac}, \bibinfo{author}{F.~Hild},
  \bibinfo{author}{P.~Ienny}, \bibinfo{author}{F.~Latourte},
  \bibinfo{author}{D.~Lemosse}, \bibinfo{author}{S.~Pagano},
  \bibinfo{author}{E.~Pagnacco}, et~al.,
\newblock \bibinfo{title}{Overview of identification methods of mechanical
  parameters based on full-field measurements},
\newblock \bibinfo{journal}{Experimental Mechanics} \bibinfo{volume}{48}
  (\bibinfo{year}{2008}) \bibinfo{pages}{381--402}.
\bibitem[{Bhatia et~al.(2025)Bhatia, Koehler, and Thuerey}]{bhatia:2025a}
\bibinfo{author}{K.~Bhatia}, \bibinfo{author}{F.~Koehler},
  \bibinfo{author}{N.~Thuerey},
\newblock \bibinfo{title}{Prdp: Progressively refined differentiable physics},
\newblock \bibinfo{journal}{arXiv preprint arXiv:2502.19611}
  (\bibinfo{year}{2025}).
\bibitem[{Bergstra et~al.(2011)Bergstra, Bardenet, Bengio, and
  K\'{e}gl}]{bergstra:2011a}
\bibinfo{author}{J.~Bergstra}, \bibinfo{author}{R.~Bardenet},
  \bibinfo{author}{Y.~Bengio}, \bibinfo{author}{B.~K\'{e}gl},
\newblock \bibinfo{title}{Algorithms for hyper-parameter optimization},
\newblock in: \bibinfo{editor}{J.~Shawe-Taylor}, \bibinfo{editor}{R.~Zemel},
  \bibinfo{editor}{P.~Bartlett}, \bibinfo{editor}{F.~Pereira},
  \bibinfo{editor}{K.~Weinberger} (Eds.), \bibinfo{booktitle}{Advances in
  Neural Information Processing Systems}, volume~\bibinfo{volume}{24},
  \bibinfo{publisher}{Curran Associates, Inc.}, \bibinfo{year}{2011}, pp.
  \bibinfo{pages}{2546--2554}.
\bibitem[{Yi et~al.(2024)Yi, Cheng, and Bessa}]{yi2024practical}
\bibinfo{author}{J.~Yi}, \bibinfo{author}{J.~Cheng}, \bibinfo{author}{M.~A.
  Bessa},
\newblock \bibinfo{title}{Practical multi-fidelity machine learning: fusion of
  deterministic and bayesian models},
\newblock \bibinfo{journal}{arXiv preprint arXiv:2407.15110}
  (\bibinfo{year}{2024}).
\bibitem[{Ihuaenyi et~al.(2024)Ihuaenyi, Luo, Li, and Zhu}]{ihuaneyi:2024a}
\bibinfo{author}{R.~C. Ihuaenyi}, \bibinfo{author}{J.~Luo},
  \bibinfo{author}{W.~Li}, \bibinfo{author}{J.~Zhu},
\newblock \bibinfo{title}{Seeking the most informative design of test specimens
  for learning constitutive models},
\newblock \bibinfo{journal}{Extreme Mechanics Letters} \bibinfo{volume}{69}
  (\bibinfo{year}{2024}) \bibinfo{pages}{102169}.
\bibitem[{Tung and Li(2024)}]{tung:2024a}
\bibinfo{author}{C.-H. Tung}, \bibinfo{author}{J.~Li},
\newblock \bibinfo{title}{The anti-dogbone: Evaluating and designing optimal
  tensile specimens for deep learning of constitutive relations},
\newblock \bibinfo{journal}{Extreme Mechanics Letters} \bibinfo{volume}{69}
  (\bibinfo{year}{2024}) \bibinfo{pages}{102157}.
\bibitem[{Jia et~al.(2025)Jia, Li, and Zhang}]{jia:2025a}
\bibinfo{author}{Y.~Jia}, \bibinfo{author}{W.~Li}, \bibinfo{author}{X.~S.
  Zhang},
\newblock \bibinfo{title}{Multimaterial topology optimization of elastoplastic
  composite structures},
\newblock \bibinfo{journal}{Journal of the Mechanics and Physics of Solids}
  \bibinfo{volume}{196} (\bibinfo{year}{2025}) \bibinfo{pages}{106018}.
\bibitem[{Vijayakumaran et~al.(2025)Vijayakumaran, Russ, Paulino, and
  Bessa}]{vijayakumaran:2025a}
\bibinfo{author}{H.~Vijayakumaran}, \bibinfo{author}{J.~B. Russ},
  \bibinfo{author}{G.~H. Paulino}, \bibinfo{author}{M.~A. Bessa},
\newblock \bibinfo{title}{Consistent machine learning for topology optimization
  with microstructure-dependent neural network material models},
\newblock \bibinfo{journal}{Journal of the Mechanics and Physics of Solids}
  \bibinfo{volume}{196} (\bibinfo{year}{2025}) \bibinfo{pages}{106015}.
\bibitem[{Dekhovich et~al.(2023)Dekhovich, Turan, Yi, and
  Bessa}]{dekhovich:2023a}
\bibinfo{author}{A.~Dekhovich}, \bibinfo{author}{O.~T. Turan},
  \bibinfo{author}{J.~Yi}, \bibinfo{author}{M.~A. Bessa},
\newblock \bibinfo{title}{Cooperative data-driven modeling},
\newblock \bibinfo{journal}{Computer Methods in Applied Mechanics and
  Engineering} \bibinfo{volume}{417} (\bibinfo{year}{2023})
  \bibinfo{pages}{116432}.
\bibitem[{Yi and Bessa(2025)}]{yi2025cooperative}
\bibinfo{author}{J.~Yi}, \bibinfo{author}{M.~A. Bessa},
\newblock \bibinfo{title}{Cooperative bayesian and variance networks
  disentangle aleatoric and epistemic uncertainties},
\newblock \bibinfo{journal}{arXiv preprint arXiv:2505.02743}
  (\bibinfo{year}{2025}).
\bibitem[{Coscia et~al.(2025)Coscia, Welling, Demo, and Rozza}]{coscia:2025a}
\bibinfo{author}{D.~Coscia}, \bibinfo{author}{M.~Welling},
  \bibinfo{author}{N.~Demo}, \bibinfo{author}{G.~Rozza}, \bibinfo{title}{Barnn:
  A bayesian autoregressive and recurrent neural network},
  \bibinfo{year}{2025}. \href{http://arxiv.org/abs/2501.18665}{\tt
  arXiv:2501.18665}.
\bibitem[{Lou et~al.(2022)Lou, Zhang, Zhang, and Yoon}]{lou:2022a}
\bibinfo{author}{Y.~Lou}, \bibinfo{author}{C.~Zhang},
  \bibinfo{author}{S.~Zhang}, \bibinfo{author}{J.~W. Yoon},
\newblock \bibinfo{title}{A general yield function with differential and
  anisotropic hardening for strength modelling under various stress states with
  non-associated flow rule},
\newblock \bibinfo{journal}{International Journal of Plasticity}
  \bibinfo{volume}{158} (\bibinfo{year}{2022}) \bibinfo{pages}{103414}.
\bibitem[{Cazacu and Barlat(2004)}]{barlat:2004a}
\bibinfo{author}{O.~Cazacu}, \bibinfo{author}{F.~Barlat},
\newblock \bibinfo{title}{A criterion for description of anisotropy and yield
  differential effects in pressure-insensitive metals},
\newblock \bibinfo{journal}{International Journal of Plasticity}
  \bibinfo{volume}{20} (\bibinfo{year}{2004}) \bibinfo{pages}{2027--2045}.
  \bibinfo{note}{Daniel C. Drucker Memorial Issue}.
\bibitem[{Zhang and Lou(2023)}]{lou:2023a}
\bibinfo{author}{C.~Zhang}, \bibinfo{author}{Y.~Lou},
\newblock \bibinfo{title}{Characterization and modelling of evolving plasticity
  behaviour up to fracture for fcc and bcc metals},
\newblock \bibinfo{journal}{Journal of Materials Processing Technology}
  \bibinfo{volume}{317} (\bibinfo{year}{2023}) \bibinfo{pages}{117997}.
\bibitem[{Li et~al.(2023)Li, Persaud, Choudhary, DeCost, Greenwood, and
  Hattrick-Simpers}]{li:2023a}
\bibinfo{author}{K.~Li}, \bibinfo{author}{D.~Persaud},
  \bibinfo{author}{K.~Choudhary}, \bibinfo{author}{B.~DeCost},
  \bibinfo{author}{M.~Greenwood}, \bibinfo{author}{J.~Hattrick-Simpers},
\newblock \bibinfo{title}{Exploiting redundancy in large materials datasets for
  efficient machine learning with less data},
\newblock \bibinfo{journal}{Nature Communications} \bibinfo{volume}{14}
  (\bibinfo{year}{2023}).
\bibitem[{Cho et~al.(2014)Cho, Van~Merri{\"e}nboer, Gulcehre, Bahdanau,
  Bougares, Schwenk, and Bengio}]{cho:2014a}
\bibinfo{author}{K.~Cho}, \bibinfo{author}{B.~Van~Merri{\"e}nboer},
  \bibinfo{author}{C.~Gulcehre}, \bibinfo{author}{D.~Bahdanau},
  \bibinfo{author}{F.~Bougares}, \bibinfo{author}{H.~Schwenk},
  \bibinfo{author}{Y.~Bengio},
\newblock \bibinfo{title}{Learning phrase representations using rnn
  encoder-decoder for statistical machine translation},
\newblock \bibinfo{journal}{arXiv preprint arXiv:1406.1078}
  (\bibinfo{year}{2014}).
\bibitem[{Chung et~al.(2014)Chung, Gulcehre, Cho, and Bengio}]{chung:2014a}
\bibinfo{author}{J.~Chung}, \bibinfo{author}{C.~Gulcehre},
  \bibinfo{author}{K.~Cho}, \bibinfo{author}{Y.~Bengio},
\newblock \bibinfo{title}{Empirical evaluation of gated recurrent neural
  networks on sequence modeling},
\newblock \bibinfo{journal}{arXiv preprint arXiv:1412.3555}
  (\bibinfo{year}{2014}).

\end{thebibliography}
	
\end{document}